\documentclass[11pt]{article}
\usepackage[utf8]{inputenc}

\usepackage[margin=1in]{geometry}

\usepackage{lmodern, setspace}
\usepackage{amsthm}
\usepackage{caption}
\usepackage{subcaption}
\usepackage{thm-restate}
\usepackage{extarrows}
\usepackage[normalem]{ulem}

\usepackage{graphics}

\usepackage{amsmath, amssymb, bm, xcolor}
\usepackage{xspace}

\DeclareMathOperator{\argmin}{argmin}

\DeclareMathOperator{\evec}{\bf e}

\DeclareMathOperator{\clusters}{\sf CLR}
\DeclareMathOperator{\diam}{\sf diam}
\DeclareMathOperator{\dist}{\sf dist}
\DeclareMathOperator{\families}{\sf FML}

\DeclareMathOperator{\cH}{\cal H}
\DeclareMathOperator{\cD}{\cal D}

\DeclareMathOperator{\cQ}{\cal Q}
\DeclareMathOperator{\cA}{\cal A}

\DeclareMathOperator{\cF}{\cal F}

\DeclareMathOperator{\bd}{\bf d}
\DeclareMathOperator{\bu}{\bf u}

\DeclareMathOperator{\cP}{\cal P}
\DeclareMathOperator{\bq}{\bf q}

\DeclareMathOperator{\bw}{\bf w}

\DeclareMathOperator{\bx}{\bf x}
\DeclareMathOperator{\bt}{\bf t}
\DeclareMathOperator{\bv}{\bf v}
\DeclareMathOperator{\by}{\bf y}

\DeclareMathOperator{\bmu}{\bm{\mu}}

\DeclareMathOperator{\bzero}{\bf 0}
\DeclareMathOperator{\bsigma}{\bm{\sigma}}
\newcommand*{\trsp}{^{\mkern-1.5mu\top}}

\newtheorem{theorem}{Theorem}
\newtheorem{lemma}{Lemma}

\newtheorem{algorithm}{Algorithm}
\newtheorem{definition}{Definition}

\newtheorem{claim}{Claim}
\newtheorem{remark}{Remark}

\usepackage{enumerate}

\usepackage{algorithm, algpseudocode}

\usepackage[colorlinks=true, linkcolor=blue, citecolor=blue]{hyperref}

\usepackage[color=green!40]{todonotes}

\usepackage[style=authoryear, sorting=nyt]{biblatex}
\DeclareCiteCommand{\cite}
  {\usebibmacro{prenote}}
  {\bibhyperref{%
     \printnames{labelname}%
     \setunit{\nameyeardelim}%
     \printfield{year}%
   }}
  {\multicitedelim}
  {\usebibmacro{postnote}}

\usepackage{lipsum}

\begin{document}
\title{Distributionally Robust Newsvendor on a Metric}
\author{
  Ayoub Foussoul\\
  {\footnotesize Industrial Engineering and Operations Research, Columbia University, New York, USA, 
  \href{mailto:af3209@columbia.edu}{af3209@columbia.edu}}
  \and
  Vineet Goyal\\
  {\footnotesize Industrial Engineering and Operations Research, Columbia University, New York, USA, 
  \href{mailto:vg2277@columbia.edu}{vg2277@columbia.edu}}
}
\date{}
\maketitle

\begin{abstract}
We consider a fundamental generalization of the classical newsvendor problem where the seller needs to decide on the inventory of a product jointly for multiple locations on a metric as well as a fulfillment policy to satisfy the uncertain demand that arises sequentially over time after the inventory decisions have been made. To address the distributional-ambiguity, we consider a distributionally robust setting where the decision-maker only knows the mean and variance of the demand, and the goal is to make inventory and fulfillment decisions to minimize the worst-case expected inventory and fulfillment cost.
We design an approximation policy for the problem with theoretical guarantees on its performance. Our policy generalizes the classical solution of~\cite{scarf1957min}, maintaining its simplicity and interpretability: it identifies a hierarchical set of clusters, assigns a “virtual” underage cost to each cluster, then makes sure that each cluster holds at least the inventory suggested by Scarf’s solution if the cluster behaved as a single point with “virtual” underage cost. As demand arrives sequentially, our policy fulfills orders from nearby clusters, minimizing fulfillment costs, while balancing inventory consumption across the clusters to avoid depleting any single one. We show that the policy achieves a poly-logarithmic approximation. To the best of our knowledge, this is the first algorithm with provable performance guarantees. Furthermore, our numerical experiments show that the policy performs well in practice.
\end{abstract}

\section{Introduction and Motivation}
An e-commerce seller has multiple warehouses and receives uncertain demand that arrives sequentially to multiple demand locations. When a demand request realizes, it can be fulfilled from any of the warehouses. The seller needs to decide about an initial inventory allocation at the warehouses and a demand fulfillment policy for the sequentially arriving demand to minimize their total expected fulfillment, overage, and underage costs. One straightforward approach is to fulfill each request from the nearest warehouse, in which case each warehouse plans its inventory independently. However, this prevents risk pooling over demand uncertainty between warehouses. In particular, by adopting a more flexible fulfillment strategy and jointly planning inventory across the warehouses, the seller can leverage the advantages of risk pooling over demand uncertainty and reduce costs. We explore in this work the benefits of risk pooling of demand uncertainty when warehouses and demand locations reside in a metric.

In practice, only partial information is typically available about the demand distribution. For instance, it is often much more difficult to estimate the exact demand distribution compared to estimating simpler statistics, such as the first and second moments. In this work, we study the problem when only the first and second moments of the demand are known, and demand across locations is uncorrelated. This arises in settings where locations operate under uncorrelated local conditions, with different customer bases, market dynamics, or regional preferences. We adopt a distributionally robust model in which the goal of the seller is to choose an initial inventory allocation and sequential demand fulfillment policy that minimize their expected cost under the worst case demand distribution and worst case sequence of arrivals of demand. We refer to this problem as the {Distributionally Robust Newsvendor on a Metric} problem~\eqref{eq:online}.

The Distributionally Robust Newsvendor on a Metric~\eqref{eq:online} is a fundamental generalization of the classical distributionally robust newsvendor problem of~\cite{scarf1957min}. In his seminal paper,~\cite{scarf1957min} studies the problem for a single location and gives a simple closed-form solution. 
In particular,~\cite{scarf1957min} shows that the optimal initial inventory level is given by,
\begin{align}\label{eq:scarf}
    q = \mu + \frac{\sigma}{2} \left(\sqrt{\frac{b}{h}} - \sqrt{\frac{h}{b}}\right),
\end{align}
or $q=0$, whichever yields a lower cost, where $\mu$ and $\sigma$ denote the mean and standard deviation respectively of the unknown demand distribution and $b$ and $h$ denote the per-unit underage and overage respectively. In~\eqref{eq:scarf}, the firm orders the expected demand $\mu$ plus a correction term that depends on the relative magnitude of the underage cost $b$ and overage cost $h$, and is scaled by the demand uncertainty (represented by $\sigma$).

The multi-location problem has also been studied before in~\cite{govindarajan2021distribution} in an offline setting when demand arrives all at once and is fulfilled using a minimum cost fulfillment. We refer to this problem as the Offline Distributionally Robust Newsvendor on a Metric~\eqref{eq:offline}. In~\cite{govindarajan2021distribution}, the authors give an SDP heuristic, based on a dual reformulation of the problem. However, no interpretable solutions with theoretical guarantees on their performance are known for this fundamental problem. 
In this work, we bridge this gap and design an approximation solution with theoretical guarantees on its performance. Moreover, our solution generalizes Scarf's solution for a single location maintaining its simplicity and interpretability and show a good performance in practice. Then, building upon this result, we design an approximation policy for the Distributionally Robust Newsvendor on a Metric problem~\eqref{eq:online} when demand arrives sequentially over time.
Below is a more detailed exposition of our contributions.

\subsection{Our Main Contributions}

We adopt a two step approach. First, we derive an approximation solution for the offline problem~\eqref{eq:offline} when demand arrives all at once. We refer to this solution as the {\em Generalized-Scarf-on-a-Metric} (GSM) solution. Then, building on this result, we derive an approximation policy for our problem~\eqref{eq:online} where demand realizes sequentially over time. We refer to this policy as the {\em Hierarchical Balance} (HB) policy.

\vspace{3mm}{\noindent \bf The { Generalized-Scarf-on-a-Metric} (GSM) solution.} We give an approximation solution for the offline problem~\eqref{eq:offline} when demand arrives all at once. We refer to our solution as the {\em Generalized-Scarf-on-a-Metric} (GSM) solution.

Our solution is a generalization of Scarf's solution maintaining its simplicity and interpretability, and can be described as follows: Our solution identifies a set of clusters of locations $\Gamma$ and assigns each cluster $C \in \Gamma$ a ``virtual" underage cost $b_C$. The solution then holds the minimum amount of inventory necessary to ensure that each cluster $C \in \Gamma$ contains a inventory of at least
$$
\sum_{i \in C} \mu_i + \frac{\sqrt{\sum_{i \in C} \sigma_i^2}}{2}  \cdot \left(\sqrt{\frac{b_C}{h}}-\sqrt{\frac{h}{b_C}}\right),
$$
and that the entire set of locations holds an inventory of at least,
$$
\sum_{i \in [n]} \mu_i + \frac{\sqrt{\sum_{i \in [n]} \sigma_i^2}}{2}\cdot \left(\sqrt{\frac{b}{h}}-\sqrt{\frac{h}{b}}\right),
$$
where $n$ denotes the number of locations, $\mu_i$ and $\sigma_i$ the mean and standard deviation respectively of demand at location $i$, and $b$ and $h$ denote the underage and overage cost respectively. Note that the first expression corresponds exactly to the ordering quantity given by Scarf solution~\eqref{eq:scarf} with mean and standard deviation corresponding to the total demand within cluster $C$, subject to the ``virtual" underage cost $b_C$ and original overage cost $h$. Similarly, the second expression corresponds to the ordering quantity given by Scarf solution~\eqref{eq:scarf} with mean and standard deviation corresponding to the total demand, subject to underage cost $b$ and overage cost $h$. In fact, our solution is a generalization of Scarf's original solution for a single point.

We focus in this paper on the setting where the per-unit underage cost $b$ is larger than the
per-unit overage cost $h$. This typically holds in practice as underages cause customer
dissatisfaction, lost revenue, and reputational harm, while overages lead to storage, insurance, and
depreciation costs, which are generally less immediate and severe.
In particular, we have the following theorem.
\begin{theorem}(informal) Suppose $b>h$. The Generalized-Scarf-on-a-Metric (GSM) solution is a polylogarithmic approximation to~\eqref{eq:offline}. Furthermore, it is a constant approximation for the special case of uniform metric.
\end{theorem}

To provide some intuition behind the structure of our solution, consider the special case where the points are separated by a fixed distance $\ell$ (uniform metric). In this case, the solution considers each point as a separate cluster with a ``virtual" underage cost of $\ell - h$, and ensures a inventory inside each point that is at least the inventory given by Scarf's solution with virtual underage cost $\ell-h$ and original overage cost $h$. This stems, roughly speaking, from the fact that a point incurs a cost of $\ell$ to acquire inventory from another point with excess inventory, while the other point saves $h$ by no longer holding the excess inventory. This setup protects against fulfillment costs. Additionally, our solution ensures a total inventory that is at least the inventory given by Scarf's solution if the points were combined into a single point. This protects against the overall underage costs. We show that ensuring that both of these conditions are met leads to a constant approximation. Now, in general metric spaces, we introduce special partitions of the metric space, which we refer to as Hierarchically Well-Separated Partitions and which allow us to ``mimic" uniform metrics.

\vspace{3mm}{\noindent \bf The Hierarchical Balance (HB) policy.}
Building upon the previous result, 
we present an approximation policy for the online problem~\eqref{eq:online} where we observe the demand in a sequential manner and need to decide where it must be fulfilled from at the time it materializes. In particular, we present a sequential fulfillment policy that makes sequential assignment decisions starting from the inventory solution obtained using our (GSM) solution and that gives a polylogarithmic approximation. Specifically, our policy fulfills demand from clusters containing the demand in increasing order of diameter. Within each cluster, inventory is drawn equally from each of its subclusters, and equally from each of the subclusters of each subcluster, and so on. We refer to our policy as the {\em Hierarchical Balance (HB)} policy. We have the following theorem:
\begin{theorem}(informal) Suppose $b>h$. The Hierarchical Balance (HB) policy is a  polylogarithmic approximation to~\eqref{eq:online}.
\end{theorem}

Intuitively, we are fulfilling the demand from clusters in increasing order of diameter to ensure small fulfillment costs. Within each cluster, we make sure to draw inventory evenly from its subclusters and evenly from their subclusters and so on, to prevent any single cluster from being depleted excessively, as we remain uncertain about where future demand will arise. To prove our result we use a potential function argument that allows us to link the cost from an online fulfillment of demand to the cost from a special offline fulfillment of demand.

We note that our Hierarchical Balance (HB) policy closely resembles heuristics employed by major e-commerce businesses. For instance, in 2022, Amazon adopted a regional fulfillment approach in the U.S., dividing the country into eight regions. Customer orders are initially fulfilled locally within the designated region, and only if unavailable, the order is fulfilled from other regions~\cite{amazon}. This approach mirrors our policy when two clustering levels are considered, and our results therefore provide a theoretical backing for the practical effectiveness of such hierarchical inventory strategies.

We also remark that the clustering used in our policy depends solely on the underlying metric space and is independent of the demand distribution. As a result, the same clustering can be applied to multiple products simultaneously, and our results remain valid in multi-product settings.

\subsection{Related Work}

The offline multi-location newsvendor problem is extensively studied in the literature when the demand distribution is known exactly, with a primary focus on examining the effects of centralization and inventory pooling on the performance. The problem was first introduced by~\cite{eppen1979note}, who studied the problem under normally distributed demand and no fulfillment costs, highlighting the benefits of centralization in risk pooling. The results of~\cite{eppen1979note} were subsequently extended to more general distributions and fulfillment costs~\cite{chen1989effects,lin2001effects,corbett2006generalization}. Further work explored the impact of demand variability on the benefits of risk pooling~\cite{berman2011benefits,gerchak2003relation},~\cite{yang2021multilocation} analyzed the problem in a risk-averse framework, and distinguished between the effects of centralization and inventory pooling, and \cite{herer2023asymptotic} studied the asymptotic performance as the number of locations increases.

In practice, the demand distribution is often unknown. Distributionally robust optimization (DRO) gives a good framework to address such uncertainty, where the objective is to optimize against the worst-case distribution within a specified ambiguity set. Various ambiguity sets have been studied in the literature, which can be generally categorized into three types. The first category includes sets containing distributions supported on single points, reducing the (DRO) problem to a robust optimization problem \cite{ben1998robust,bertsimas2011theory,el2021optimality}. The second category consists of statistically based ambiguity sets, which include distributions near a nominal (typically empirical) distribution, with proximity measured using statistical distances such as Wasserstein metrics \cite{pflug2007ambiguity,mohajerin2018data}, $\zeta$-structures \cite{zhao2018data}, or $\Phi$-divergences \cite{jiang2016data}. The third category comprises moment-based ambiguity sets, such as the one considered in this work, which include distributions that satisfy specific moment constraints \cite{delage2010distributionally,zymler2013worst,wiesemann2014distributionally}.

The offline multi-location newsvendor problem has been studied within a distributionally robust framework. For a single location,~\cite{scarf1957min} considered the problem under a moment-based ambiguity set and provided a complete solution with a closed form (see~\cite{gallego1993distribution} for a simpler derivation of the same result). For multiple locations, the problem has been explored under different ambiguity sets. \cite{chou2006robust} considered the problem under a chance-constrained ambiguity set and proposed an approximation based on linear decision rules. \cite{govindarajan2021distribution} considered the problem under a moment-based ambiguity set similar to ours, and gave an SDP heuristic. More recently, \cite{li2022distributionally} studied the problem under a scenario-based ambiguity set and proposed a solution through reformulation as a scalable linear program. A similar inventory allocation then demand fulfillment on a metric problem has been studied by~\cite{el2021power} in a two-stage robust framework.

The idea of approximating metric spaces with families of well-separated clusters, where algorithms are easier to design and analyze, has proved to be useful in a variety of contexts. For instance, a well-known technique in designing algorithms for (integral) online bipartite matching problems is the use of randomized embeddings of metric spaces into better clustered metric spaces such as the so-called Hierarchically Well-Separated Trees (HSTs)~\cite{bartal1998approximating,fakcharoenphol2003tight,bansal2007log,meyerson2006randomized,gupta2012online,gupta2019stochastic}. We note in passing that randomized embeddings are not applicable in a distributionally robust framework, where the objective function is non-linear with respect to the metric distances. Approximating metric spaces by well-separated families of clusters has also been effective in designing distributed graph algorithms \cite{linial1993low,chang2020deterministic,chang2021near} and approximating NP-hard graph optimization problems \cite{even2000divide}.

\section{Our Model and Notation}
In the Distributionally Robust Newsvendor on a Metric problem~\eqref{eq:online}, we are given a metric space $(X,\ell)$ consisting of $|X| = n$ locations, where $\ell_{ij}$ denotes the distance between locations $i \in X$ and $j \in X$. Each location represents both a warehouse and a demand location. At the beginning of the selling season, a firm determines the inventory level $q_i \in \mathbb{R}_+$ at each location $i \in X$. Then, an adversary samples a demand vector $\bd = (d_i)_{i \in X} \in \mathbb{R}_+^X$ from an unknown distribution $\cD$ of mean $\bmu=(\mu_i)_{i \in X}$ and diagonal covariance matrix $\Sigma = (\Sigma_{ij})_{i\in X, \; j \in X}$ where $\mu_i>0$ and $\Sigma_{ii} > 0$ are the mean and variance at location $i \in X$ respectively. The adversary then reveals the demand sequentially, in parts, in the worst possible sequence. As each portion of demand arrives, the firm must fulfill it
immediately from the available inventory thus far. The total cost incurred includes a per-unit cost equal to the distance $\ell_{ij}$ between locations $i$ and $j$ for shipped inventory between $i$ and $j$, and a per-unit underage (resp. overage) cost $b$ (resp. $h$) for unfulfilled demand (resp. remaining inventory) at the end of the selling season. The firm’s objective is to determine both an initial inventory allocation
and an online fulfillment policy that minimize the expected cost under the worst-case demand
distribution and the worst-case sequence of demand arrivals.   In particular, let $\cF(\bmu, \Sigma)$ denote the set of all distributions over demand vectors of mean $\bmu$ and covariance matrix $\Sigma$. For a demand vector $\bd \geq \bzero$, let 
$\Pi(\bd)= \bigcup_{T=1}^{\infty} \left\{(\bd^1, \dots, \bd^T) \in (\mathbb{R}_+^X)^T\;\left|\; \sum_{t=1}^T \bd^t = \bd\right.\right\}$
denote the set of all possible arrival sequences of $\bd$. The firm solves the problem,
\begin{align}
\label{eq:online}
\tag{DRNM}
    \inf_{\bq \geq \bzero, \; \cA}\; \left\{ \sup_{\cD \in \cF(\bmu, \Sigma)} \mathbb{E}_{\bd \sim \cD}\; \left(\sup_{\pi \in \Pi(\bd)}C_{\cA}(\bq, \pi)\right)\right\},
\end{align}
where if we let $x^{\cA}_{ij}(\bq, \pi)$~\footnote{In particular, $x^{\cA}_{ii}(\bq, \pi)$ denotes the amount of inventory used to fulfill demand in-place.} denote the amount of inventory shipped between locations $i$ and $j$ by the online fulfillment policy $\cA$ under the initial inventory allocation $\bq$ and demand arrival sequence $\pi$, the cost function $C_{\cA}(\bq, \pi)$ is given by,
\begin{align*}
    C_{\cA}(\bq, \pi) = h\sum_{i \in X}\left(q_i - \sum_{j\in X} x^{\cA}_{ij}(\bq, \pi)\right) + b \sum_{j \in X} \left(d_j-\sum_{i \in X} x^{\cA}_{ij}(\bq, \pi)\right) + \sum_{i \in X}\sum_{j \in X} \ell_{ij} x^{\cA}_{ij}(\bq, \pi).
\end{align*}
Here $h\sum_{i \in X}\left(q_i - \sum_{j\in X} x^{\cA}_{ij}(\bq, \pi)\right)$ is the total overage cost, $b \sum_{j \in X} \left(d_j-\sum_{i \in X} x^{\cA}_{ij}(\bq, \pi)\right)$ is the total underage cost, and $\sum_{i \in X}\sum_{j \in X} \ell_{ij} x^{\cA}_{ij}(\bq, \pi)$ is the total fulfillment cost. We denote the optimal value of~\eqref{eq:online} by $z_{\textsf{DRNM}}$. We focus in this paper on the setting where the per-unit underage cost $b$ is larger than the per-unit overage cost $h$, i.e., $b \geq h$. 

The Offline Distributionally Robust Newsvendor on a Metric problem~\eqref{eq:offline} where demand arrives all at once and gets fulfilled using a minimum cost fulfillment is given formally by the distributionally robust problem,
\begin{align}
\label{eq:offline}
\tag{ODRNM}
     \inf_{\bq \geq \bzero}\quad \left\{\sup_{\cD \in \cF(\bmu, \Sigma)} \mathbb{E}_{\bd \sim \cD}\; C(\bq, \bd)\right\},
\end{align}
where the cost function $C(\bq,\bd)$ is given by
$$
\begin{array}{rcl}
\min\limits_{\mathbf{x}} & \quad & h\sum\limits_{i \in X}\Bigl(q_i - \sum\limits_{j \in X} x_{ij}\Bigr)
+ b\sum\limits_{j \in X}\Bigl(d_j - \sum\limits_{i \in X} x_{ij}\Bigr)
+ \sum\limits_{i \in X}\sum\limits_{j \in X} \ell_{ij}\, x_{ij} \\[1ex]
\text{s.t.} & \quad & \sum\limits_{j \in X} x_{ij} \le q_i,\quad \forall\, i \in X, \\
           & \quad & \sum\limits_{i \in X} x_{ij} \le d_j,\quad \forall\, j \in X, \\
           & \quad & x_{ij} \ge 0,\quad \forall\, i,j \in X.
\end{array}
$$
Here, $x_{ij}$ denotes the amount of inventory shipped from location $i$ to location $j$. The total underage, overage, and fulfillment cost is minimized while ensuring that the amount of inventory shipped out of location $i \in X$ does not exceed the inventory at $i$, and that the amount of inventory shipped into location $j \in X$ does not exceed the demand at $j$. We denote the optimal objective of~\eqref{eq:offline} by $z_{\textsf{ODRNM}}$.

\begin{remark}[Generalizations]
    We note that our model subsumes the more general model where demand locations are different from warehouses (see Appendix~\ref{apx:generallocations}). For simplicity and to ease readability we restrict ourselves to the simplified model where demand locations and warehouses are the same in the remainder of the paper.

    We also note that our algorithm is robust to misspecification in the demand mean and variance in the sense that if upper-bound estimates of the demand mean and variance at each location (e.g., obtained from high probability upper-confidence bounds estimated from data) are used instead of the true parameters, our approximation bounds still hold, with an extra error that vanishes as the estimates become more accurate. The details are provided in Appendix~\ref{apx:misspecification}.
\end{remark}

\paragraph{Notation.} In the remainder of the paper, we denote by $\kappa = \frac{\max_{\substack{i,j \in X:i\neq j}} \ell_{ij}}{\min_{\substack{i,j \in X:i\neq j}} \ell_{ij}}$ the aspect ratio of $(X,\ell)$, by $\diam(S)= \max_{\substack{i,j \in S}} \ell_{ij}$ the diameter of a subset $S \subset X$ , and by $\dist(S,T) = \min_{i \in S, j \in T} \ell_{ij}$ the distance between two subsets $S, T \subset X$. Let $\bsigma=(\sigma_i)_{i \in X}$ denote the vector of standard deviations of the distributions $\cF(\bmu, \Sigma)$, let $\bsigma^2=(\sigma^2_i)_{i \in X}$ denote the vector of variances. Given a subset $S \subset X$ and a vector $\bt \in \mathbb{R}^X$, we denote by $\bt_S=(t_i)_{i \in S}$ the sub-vector of $\bt$ corresponding to the components $S$, and we denote $t_S = \sum_{i \in S} t_i$, except for $\bt = \bsigma$ where we denote $\sigma_S = \sqrt{\sum_{i \in S} \sigma^2_i}$ instead. Similarly, given a square matrix $M \in \mathbb{R}^{X \times X}$ and a subset $S \subset X$, we denote by $M_S$ the square submatrix of $M$ corresponding to the rows and columns of indices in $S$. For a subset $S \subset X$, we denote by $\cF_{S}(\bmu, \Sigma)$ the set of distributions that are the restrictions of the distributions $\cF(\bmu, \Sigma)$ to $S$, that is, if $\cD \in \cF_{S}(\bmu, \Sigma)$ then $\bd \sim \cD$ is supported on $\mathbb{R}_+^{S}$, has mean $\bmu_S$ and covariance matrix $\Sigma_{S}$. Note that in particular $\cF_{X}(\bmu, \Sigma)=\cF(\bmu, \Sigma)$. Let $\nu = \max_{i \in X} \frac{\sigma_i}{\mu_i}$ denote the maximum coefficient of variation. Let $(x)^+$ denote the positive part of $x \in \mathbb{R}$. As the reader may have already noticed, vectors are denoted by bold symbols throughout the paper, while scalars and matrices are represented using regular notation.

\section{The Generalized-Scarf-on-a-Metric (GSM) Solution}\label{sec:policy}

We give in this section an approximation solution for the offline problem~\eqref{eq:offline} where demand realizes all at once. Our solution consists of an initial inventory allocation $\bq^{\textsf{GSM}}$ and we refer to it as the {\em Generalized-Scarf-on-a-Metric} (GSM) solution. We begin with a series of preliminary results that simplify our problem. 

\subsection{Preliminaries.}
The following proposition shows that one can assume without loss of generality that $b+h \geq \diam(X)$. 

\begin{restatable}{proposition}{restatebhl}\label{prop:bhl}
    It can be assumed without loss of generality that $b+h \geq \ell_{ij}$ for every $i,j \in X$. Also, given $\bq, \bd \geq \bzero$, there exists an optimal solution of $C(\bq, \bd)$ that fulfills all the demand or consumes all of the inventory.
\end{restatable}

Proposition~\ref{prop:bhl} allows us to simplify the expression of the cost function $C(\bq, \bd)$ for $\bq, \bd \geq \bzero$ as
$$
C(\bq,\bd)
= h\left(q_X - d_X\right)^+ + b\left(d_X - q_X\right)^+ + \Theta(\bq, \bd),
$$
where $\Theta(\bq, \bd)$ is equal to
$$
\inf_{\bx \geq \bzero} \left\{ \sum_{i,j} \ell_{ij} x_{ij} \;\left|\; \sum_{j} x_{ij} \leq q_i, \; \sum_i x_{ij} = d_j, \; \forall i,j \in X\right.\right\}
$$
if $d_X \leq q_X$, and to
$$
\inf_{\bx \geq \bzero} \left\{ \sum_{i,j} \ell_{ij} x_{ij} \;\left|\; \sum_{j} x_{ij} = q_i, \; \sum_i x_{ij} \leq d_j, \; \forall i,j \in X\right.\right\}
$$
if $d_X > q_X$.
Here $(q_X-d_X)^+$ and $(d_X-q_X)^+$ are the total inventory overage and demand underage respectively, after the maximum possible demand is fulfiled, and $\Theta(\bq, \bd)$ is the fulfillment cost.

Next, we show that it is sufficient to consider inventory vectors $\bq$ that order at least the mean $\mu_i$ at each location $i \in X$. But before, let us show the following lemma that we use in proving this fact, and that we use in multiple other occasions in the paper. The lemma gives a lower-bound on the probability that a demand vector $\bd$ sampled from a distribution of $\cF_S(\bmu, \Sigma)$ for some $S \subset X$ exceeds its mean point-wise by a multiple of its variance vector. 

\begin{restatable}{lemma}{restateChebyshev}\label{lem:chebyshev}
    For every $\alpha \geq 1$ and $S \subset X$, there exists a distribution $\cD \in \cF_S(\bmu, \Sigma)$ such that, 
    $$\mathbb{P}_{\bd \sim \cD}\left( \bd \geq \bmu_S + \frac{\alpha}{\sigma_S}\bsigma_S^2\right) = \Omega\left(\frac{1}{1+\alpha^2}\right).$$
\end{restatable}
The constants in the above expression depends solely on the maximum coefficient of variation $\nu$ which is assumed here to be constant. Note that in practice, the coefficient of variation of demand typically ranges between 0 and 1,
with a coefficient larger than 1 indicating a highly variable and unpredictable demand. We note that for distributions supported on $\mathbb{R}^{|S|}$ (i.e., which can also take negative values), the bound in Lemma~\ref{lem:chebyshev} can be obtained from the multi-dimensional Chebyshev inequality of~\cite{marshall1960multivariate} (see also~\cite{bertsimas2005optimal}), which gives a lower bound on the probability that a vector with given first and second moment belongs to a given convex set. The main challenge here is to construct a distribution of non-negative support that (approximately) achieves such bound.

Next, the following lemma shows that we can assume without loss of generality that an inventory of at least $\mu_i$ is ordered at each location $i \in X$.

\begin{restatable}{lemma}{restateAtLeastMu}\label{lem:AtLeastMu}
    There exists an $O(1)$-approximation solution $\bq$ of \eqref{eq:offline} such that $\bq \geq \bmu$.
\end{restatable}

To prove Lemma~\ref{lem:AtLeastMu}, we consider an optimal solution $\bq^*$ of~\eqref{eq:offline}, then we construct a modified solution $\tilde{\bq}^*$ that is somehow the closest vector to $\bq^*$ that ``fills up" the locations having inventory less than the mean. Next, we use a special triangle inequality to upper-bound the cost $C(\tilde{\bq}^*, \bd)$ in function of the cost $C(\bq^*, \bd)$, a cost that accounts for the shipping between $\bq^*$ and $\tilde{\bq}^*$, and a cost that accounts for any extra units that $\tilde{\bq}^*$ might contain compared to $\bq^*$. Then, using Lemma~\ref{lem:chebyshev}, we show that the two later costs do not exceed $O(1) \cdot z_{\textsf{ODRNM}}$.

Given Lemma~\ref{lem:AtLeastMu}, we restrict the search to only initial inventory allocation vectors such that $\bq \geq \bmu$, while losing a constant factor in our approximation bound. The objective is henceforth to find the {\em safety stock} to order above the mean at each location. In the following, we let $\cQ = \{\bq \;|\; \bq \geq \bmu\}$. 

\subsection{High Level Idea}

The high level idea behind our solution is as follows: First, consider the special case or a uniform metric where the points are separated by a fixed distance $\ell$ ($> h$), our solution associates to each point $i$ a ``virtual" underage cost of $b_i = \ell - h$ and ensures an inventory of at least $\mu_i + \frac{\sigma_i}{2} \cdot \left(\sqrt{\frac{\ell-h}{h}}-\sqrt{\frac{h}{\ell-h}}\right)$ at $i$, as given by Scarf's solution~\eqref{eq:scarf} with underage cost $b_i$ and overage cost $h$. This stems, roughly speaking, from the fact that a point incurs a cost of $\ell$ to acquire inventory from another point with excess inventory, while the other point saves $h$ by no longer holding the excess inventory. Having such inventory inside each cluster ensures protecting against fulfillment costs. Additionally, our solution ensures a total inventory of at least
$
\mu_X + \frac{\sigma_X}{2} \cdot \left(\sqrt{\frac{\ell-h}{h}}-\sqrt{\frac{h}{\ell-h}}\right)
$ 
overall, as given by Scarf's solution~\eqref{eq:scarf} if all the points were combined into a single point with underage cost $b$ and overage cost $h$. Having such inventory overall protecting against the overall underage costs. In general metric spaces, we introduce special partitions of the metric space, which we refer to as Hierarchically Well-Separated Partitions and which cover the metric space using a ``small number" of families of distant clusters (subsets of
locations) in a hierarchical manner, such that each family ``mimics" a uniform metric. We now introduce {\em Well-Separated Hierarchical Partitions}.

\subsection{Well-Separated Hierarchical Partitions.}

We now introduce a special type of partitions of metric spaces that we call {\em Well-Separated Hierarchical Partitions}. Let us first define the notion of a {\em Well-Separated Partition.}

\begin{definition}[Well-Separated Partition]
    Let $\alpha, \beta \geq 1$, and let $\Delta >0$. An $(\alpha, \beta)$-Well-Separated Partition of $(X, \ell)$ of margin $\Delta$ is a set of $K \leq \beta$ families of clusters $\cF_1, \dots, \cF_K$, such that,
    \begin{itemize}
        \item Each family of clusters $\cF_k = \{C_1, \dots, C_{m_k}\}$ where $C_i \subset X$ for every $i \in [m_k]$, is well separated in the sense that the distance between any pair of clusters is at least $\Delta$, i.e.,
        $$
        \dist(C_i, C_j) > \Delta, \quad \forall i\neq j \in [m_k].
        $$
        \item The diameter of every cluster is at most $\alpha\cdot \Delta$, i.e,
        $$
        \diam(C) < \alpha \cdot \Delta, \quad \forall C \in \bigcup_{k \in K} \cF_k.
        $$
        \item The set of all clusters $\bigcup_{k \in K} \cF_k$ forms a partition of $X$.
    \end{itemize}
\end{definition}

\begin{figure}
    \centering
    \includegraphics[width=1\linewidth]{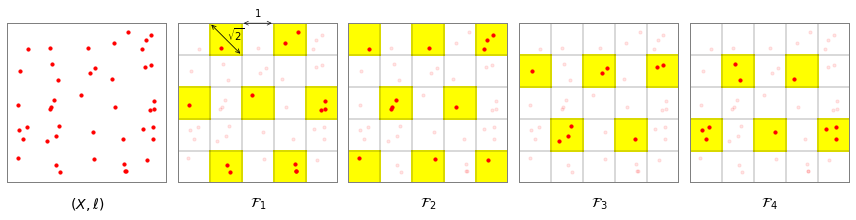}
    \caption{Example of a $(\sqrt{2}, 4)$-well-separated partition of margin $1$ of the metric space $(X, \ell)$ given by the red points in the Euclidean plane. The partition is given by $K=4$ families $\cF_1, \cF_2, \cF_3, \cF_4$ where each family is given by the points inside the yellow squares. The union of the families is the whole space $(X, \ell)$, the clusters within each family are at least a distance of $1$ away from each other, and the diameter of each cluster is at most $\sqrt{2} \times 1 = \sqrt{2}$.}
    \label{fig:wsp}
\end{figure}

\noindent See Figure~\ref{fig:wsp} for an example. We now define {\em Well-Separated Hierarchical Partitions}.

\begin{definition}[Well-Separated Hierarchical Partition]\label{def:wshp}
    Let $\alpha, \beta \geq 1$ and $\gamma > 1$.
    Let 
    $$
    \Delta_1 = \frac{1}{\alpha}\max\left\{\min_{i,j \in X: i \neq j} \ell_{ij}, \frac{\max_{i,j \in X: i \neq j} \ell_{ij}}{n}\right\}$$ 
    and
    $$ 
    R = \left\lceil \frac{\log\left(\max_{i,j \in X: i \neq j} \ell_{ij}/\Delta_1\right)}{\log(\gamma)} \right\rceil + 1
    $$
    and for every $r=1, \dots, R$, 
    $$
    \Delta_{r} = \gamma^{r-1} \Delta_{1}.$$
    
    An $(\alpha, \beta, \gamma)$-Well-Separated Hierarchical Partition of $X$ is a sequence $\cH = \cP_1,\dots, \cP_{R}$ of $R$ well-separated partitions of $X$ such that,
    \begin{itemize}
        \item For every $r \in \{1, \dots, R\}$, $\cP_r$ is an $(\alpha, \beta)$-well-separated partition of margin $\Delta_{r}$.
        \item For every $r \in \{1, \dots, R-1\}$, $\cP_{r+1}$ is a coarsening of $\cP_{r}$, that is, every cluster $C$ of some family of $\cP_r$ is contained in some cluster $C'$ of some family of $\cP_{r+1}$.
    \end{itemize}
\end{definition}

In the following, given an $(\alpha, \beta, \gamma)$-well-separated partition $\cH = \cP_1,\dots,\cP_{R}$ of $(X,\ell)$, 
and given $r \in \{1, \dots, R\}$, we denote by $\families(\cP_r)$ the set of families of clusters of $\cP_r$. Next, we denote by $$\clusters(\cH)= \bigcup_{r=1}^{R}\bigcup_{\cF \in \families(\cP_r)} \cF$$ the set of all clusters in $\cH$, and denote by $$\clusters(\cP_r) = \bigcup_{\cF \in \families(\cP_r)} \cF$$ the set of all clusters in $\cP_r$. Then for $r \in \{2, \dots, R\}$, and a cluster $C \in \clusters(\cP_r)$, we denote by $$\clusters(C) = \{C' \;|\; C' \in \clusters(\cP_{r-1}), \; C' \subset C\}$$ the set of clusters in $\cP_{r-1}$ that make $C$ (recall that $C$ at level $r$ is partitioned into smaller clusters at level $r-1$).

By definition of $\Delta_1$, there are two possible cases. First, the case when $\Delta_1 = \frac{1}{\alpha} \min_{i,j \in X: i \neq j} \ell_{ij}$. In this case, the clusters of the first partition $\clusters(\cP_1)$ are necessarily singletons. In fact, the diameter of any such cluster is strictly less than $\alpha \Delta_1 = \min_{i,j \in X: i \neq j}\ell_{ij}$ and hence the cluster must be a singleton. Second, the case when $\Delta_1 = \frac{1}{n\alpha} \max_{i,j \in X: i \neq j} \ell_{ij}$. In this case, the clusters of the first partition $\clusters(\cP_1)$ might have more than one element but have a diameter that is strictly less than $\max_{i,j \in X: i \neq j} \ell_{ij}/n$, which is, as we will see later, a small enough diameter for our purposes.

Next, by definition of $R$, $\Delta_R \geq \diam(X)$. This implies that the last partition $\cP_R$ contains a single family with a single cluster given by the whole space $X$. This is because there cannot exist two clusters inside $X$ separated by strictly more than $\diam(X)$.

Hence, a well-separated-hierarchical partition starts with clusters $\clusters(\cP_1)$ that are either singletons or that have a small diameter and ends with a single cluster in $\clusters(\cP_R)$ given by the whole set $X$.

The following lemma constructs examples of well-separated hierarchical partitions of small parameters for uniform metric spaces, subspaces of an Euclidean space $\mathbb{R}^d$ of small dimension $d$, and any general metric space. 

\begin{restatable}{lemma}{restatewshp}
\label{lem:wshp}
If $(X, \ell)$ is,
\begin{itemize}
    \item A uniform metric space (i.e., such that $\ell_{ij}=\lambda$ for every $i \neq j$), then, $(X, \ell)$ admits an $(\alpha,\beta,\gamma)$-well-separated hierarchical partition for every $\alpha > 1$, $\beta=1$ and $\gamma > \alpha$.
    \item A subspace of the Euclidean space $\mathbb{R}^d$, then, $(X, \ell)$ admits an $(\alpha,\beta,\gamma)$-well-separated hierarchical partition for every $\alpha > \sqrt{d}$, $\beta=2^d$ and an even integer $\gamma \geq 2$.
    \item A general metric space, then $(X, \ell)$ admits an $(\alpha, \beta, \gamma)$-well-separated hierarchical partition for every $\alpha = 6\log(n)+1$, $\beta = \log(n)$, $\gamma > 12\log(n)+2$.
\end{itemize}
All of the above well-separated hierarchical partitions can be constructed in polynomial time.
\end{restatable}

\subsection{The Generalized-Scarf-on-a-Metric (GSM) solution}\label{sec:GSM}
We are now ready to present our solution. Let $\cH=\cP_1,\dots,\cP_{R}$ be an $(\alpha,\beta,\gamma)$-well-separated hierarchical partition of $(X, \ell)$ for some parameters $\alpha,\beta, \gamma$.  Given $r \in \{1, \dots, R-1\}$ and a cluster $C \in \clusters(\cP_r)$, let $\diam_{p}(C)$ denote the diameter of the parent cluster of $C$ in $\cP_{r+1}$ (i.e., cluster of $\cP_{r+1}$ containing $C$). Let $$
\Gamma = \left\{C \in \clusters(\cH)\;|\; \diam_{p}(C) \geq 2h \right\}
$$
Then for every $C \in \Gamma$, let
$$
b_C = \diam_p(C)-h
$$
Our initial inventory allocation $\bq^{\textsf{GSM}}$ is given by an optimal solution of the following covering LP,
\begin{align}
\label{eq:gsm}
    \min_{\bq} \;\;
    & h (q_X-\mu_X) \notag\\
    s.t. \;\; 
    & q_{C} - \mu_C \geq \frac{\sigma_C}{2}\left(\sqrt{\frac{b_C}{h}} - \sqrt{\frac{h}{b_C}}\right), \;\; \forall C \in \Gamma \notag\tag{GSM}\\
    & q_i - \mu_i \geq \frac{\sigma_i}{\sum_{i\in X} \sigma_i} \frac{\sigma_X}{2}\left(\sqrt{\frac{b}{h}} - \sqrt{\frac{h}{b}}\right), \;\; \forall i \in X \notag
\end{align}
In~\eqref{eq:gsm}, we are minimizing the total safety stock across all locations while making sure in the first set of constraints that each cluster $C \in \Gamma$ contains a safety stock as given by the Scarf solution~\eqref{eq:scarf} with mean and standard deviation corresponding to the total demand within $C$ subject to a per-unit overage cost $h$ and per-unit underage cost $b_C$. The second set of constraints make sure that the entire set $X$ contains a safety stock $\frac{\sigma_X}{2}\left(\sqrt{\frac{b}{h}} - \sqrt{\frac{h}{b}}\right)$ as given by the Scarf solution~\eqref{eq:scarf} with mean and standard deviation corresponding to the total demand within $X$ subject to a per-unit overage cost $h$ and per-unit underage cost $b$. For technical reasons, we make sure to divide this requirement between locations proportional to the standard deviation. In fact, this only is necessary in the case when the clusters $\cP_1$ of the first level are not individual locations, and as we will see, such division makes sure that the total expected fulfillment costs between the locations inside the clusters of $\cP_1$ (of diameter at most $\diam(X)/n$ in this case) is bounded by a fraction of the total underage cost, and hence negligible. 

We show next that our solution gives a good approximation to~\eqref{eq:offline}. In particular, we show the following lemma,
\begin{restatable}{lemma}{restatemain1}
    \label{lem:main1}
    The initial inventory allocation $\bq^{\textsf{GSM}}$ is such that,
    $$\sup_{\cD \in \cF(\bmu, \Sigma)} \mathbb{E}_{\bd\sim \cD}\; C(\bq^{\textsf{GSM}}, \bd) \leq O\left(\alpha\gamma\beta R\right) \cdot z_{\textsf{ODRNM}}$$
\end{restatable}

Lemma~\ref{lem:wshp} and Lemma~\ref{lem:main1} imply our first main theorem stated formally as follows,
\setcounter{theorem}{0}
\begin{theorem}\label{thm:main1}
    There exists an initial inventory solution which is an
    \begin{itemize}
        \item $O(1)$ approximation of~\eqref{eq:offline} when $(X, \ell)$ is a uniform metric space (i.e., same distance between every two points).
        \item $O(\log(n))$ approximation of~\eqref{eq:offline} when $(X, \ell)$ is a subset of an Euclidean space $\mathbb{R}^d$ of constant dimension $d$.
        \item $O(\log^4(n))$ approximation of~\eqref{eq:offline} when $(X, \ell)$ is any general metric space.
    \end{itemize}
\end{theorem}
The next section is devoted to the proof of Lemma~\ref{lem:main1}. 
We give an outline of the proof in the main part and differ the full proofs to the Appendix.

\section{Analysis of the Generalized-Scarf-on-a-Metric (GSM) Solution}\label{sec:GSMAnalysis}

The proof consists of two parts. First, we upper-bound the objective value of $\bq^{\textsf{GSM}}$. Then we show an (approximately) similar lower-bound for $z_{\textsf{ODRNM}}$.

\vspace{3mm}{\noindent \bf Upper-Bound of the Objective Value of $\bq^{\textsf{GSM}}$.} We begin by defining the notion of a hierarchical fulfillment of demand $\bd \geq \bzero$ from inventory $\bq \in \cQ$ with respect to the well-separated hierarchical partition $\cH$, illustrated in Algorithm~\ref{alg:hier-fulfill}.
\algtext*{EndIf}
\algtext*{EndFor}
\algtext*{EndWhile}
\begin{algorithm}
\caption{Hierarchical Fulfillment}\label{alg:hier-fulfill}
\begin{algorithmic}[1]
\State \textbf{Input:} Inventory vector $\mathbf q\in\mathcal Q$, demand vector $\mathbf d\ge\mathbf0$, well-separated hierarchical partition  $\mathcal H=\cP_1,\dots,\cP_L$.
\State \textbf{Output:} Hierarchical fulfillment $\bx \in\mathbb R_{+}^{X\times X}$ of demand $\bd$ from inventory $\bq$ w.r.t. $\cH$.
\State Initialize $\mathbf q_{\rm rem}\gets\mathbf q$, $\mathbf d_{\rm rem}\gets\mathbf d$, $\bx\gets0$.
\For{$i \in X$} \Comment{Fulfill maximum demand locally}
  \State $x_{ii}\gets\min\{q_{{\rm rem},i},\,d_{{\rm rem},i}\}$
  \State $q_{{\rm rem},i}\gets q_{{\rm rem},i}-x_{ii}$; \quad $d_{{\rm rem},i}\gets d_{{\rm rem},i}-x_{ii}$
\EndFor
\For{$\ell=1,\dots,L$}
  \For{each cluster $C\in\clusters(P_{\ell})$} \Comment{Fulfill a maximum demand within each cluster of $\cP_\ell$}
    \State Let $S=\{i\in C:q_{{\rm rem},i}>0\}$ and $D=\{j\in C:d_{{\rm rem},j}>0\}$
    \If{$S\neq\emptyset$ \textbf{and} $D\neq\emptyset$}
      \State Let $\mathbf{f}$ be an optimal solution of the LP, $$\displaystyle\max_{f_{ij}\ge0}\left\{\sum_{i\in S}\sum_{j\in D}f_{ij}\quad\text{s.t.}\;\sum_{j\in D}f_{ij}\le q_{{\rm rem},i}\;\forall i\in S,\;\sum_{i\in S}f_{ij}\le d_{{\rm rem},j}\;\forall j\in D\right\}$$
      \ForAll{$i\in S,\;j\in D$}
        \State $x_{ij}\gets x_{ij}+f_{ij}$; \quad $q_{{\rm rem},i}\gets q_{{\rm rem},i}-f_{ij}$; \quad $d_{{\rm rem},j}\gets d_{{\rm rem},j}-f_{ij}$
      \EndFor
    \EndIf
  \EndFor
\EndFor
\State \textbf{Return} $\bx$
\end{algorithmic}
\end{algorithm}
In particular, for a given inventory vector $\bq \in \cQ$ and demand vector $\bd \geq \bzero$, the hierarchical fulfillment of demand $\bd \geq \bzero$ from inventory $\bq \in \cQ$ with respect to $\cH$ is the solution that fulfills a maximum demand inside each location from the inventory at the location, then for the remaining demand, fulfills (arbitrarily) a maximum demand inside each cluster $C \in \clusters(\cP_{1})$ from the inventory inside $C$, then for the remaining demand, fulfills (arbitrarily) a maximum demand inside each cluster $C \in \clusters(\cP_{2})$ from the inventory inside $C$, and so on.
Let
\begin{align*}
    C^H(\bq, \bd) &= 
    h(q_X-d_X)^+ + b(d_X-q_X)^+ \\
    & \quad +\sum_{r=1}^{R-1} \sum_{C \in \clusters(\cP_{r})} \diam_p(C) \cdot (d_C - q_C)^+ + \frac{\diam(X)}{n}\sum_{i \in X}(d_i-q_i)^+
\end{align*}
The following lemma gives an upper-bound of the cost function $C(\bq, \bd)$.

\begin{restatable}{lemma}{restateUBHierarchicalMatching}\label{lem:UBHierarchicalMatching}
    For every $\bq \in \cQ$ and $\bd \geq \bzero$,
    $$
    C(\bq, \bd) \leq C^H(\bq, \bd).
    $$
\end{restatable}
\noindent The proof of Lemma~\ref{lem:UBHierarchicalMatching} relies on showing that $C^H(\bq, \bd)$ is an upper-bound of the cost one would occur if instead of using an optimal fulfillment of the demand, uses a hierarchical fulfillment with respect to $\cH$.

Next, we show the following upper-bound on the expected value of $C^H(\bq^{\textsf{GSM}}, \bd)$ under the worst case distribution of demand.

\begin{restatable}{lemma}{restateUB}\label{lem:UB}
The initial inventory allocation $\bq^{\textsf{GSM}}$ is such that,
    \begin{align*}
        &\sup_{\cD \in \cF(\bmu, \Sigma)} \mathbb{E}_{\bd\sim \cD}\; C^H(\bq^{\textsf{GSM}}, \bd) \leq 2\sigma_X \sqrt{bh} + \sum_{C \in \Gamma}\sigma_C \sqrt{b_Ch} + \sum_{C \in \clusters(H) \setminus \Gamma} \frac{\diam_p(C)}{2}\sigma_C
    \end{align*}
\end{restatable}
\noindent In Lemma~\ref{lem:UB}, we upper-bound $\sup_{\cD \in \cF(\bmu, \Sigma)} \mathbb{E}_{\bd\sim \cD}\; C^H(\bq^{\textsf{GSM}}, \bd)$ by a function of $\bq^{\textsf{GSM}}$, then use the fact that $\bq^{\textsf{GSM}}$ is an optimal solution of~\eqref{eq:gsm} to bound the resulting function.

From Lemma~\ref{lem:UBHierarchicalMatching} and~\ref{lem:UB}, we conclude that,
\begin{align}
    \label{eq:UB}
    &\sup_{\cD \in \cF(\bmu, \Sigma)} \mathbb{E}_{\bd\sim \cD}\; C(\bq^{\textsf{GSM}}, \bd) \leq 2\sigma_X \sqrt{bh} + \sum_{C \in \Gamma}\sigma_C \sqrt{b_Ch} +\sum_{C \in \clusters(H) \setminus \Gamma} \frac{\diam_p(C)}{2} \sigma_C
\end{align}

Next, we give a lower-bound of $z_{\textsf{ODRNM}}$.

\vspace{3mm}{\noindent \bf Lower-Bound of $z_{\textsf{ODRNM}}$.}
In the following, for every $r \in \{1, \dots, R-1\}$, for every cluster $C \in \clusters(\cP_r)$, denote by 
$$\overline{C} = \left\{i \in X \;|\; \exists j \in C, \; s.t. \;\ell_{ij} < \frac{\Delta_r}{2}\right\}$$ the set of all locations within distance of at most $\frac{\Delta_r}{2}$ of $C$. Given a inventory and demand vectors $\bq \in \cQ$ and $\bd \geq \bzero$, let
$$\overline{q}_C = q_C + \sum_{i \in \overline{C} \setminus C} (q_i - d_i)^+$$ be the amount of inventory around the locations of $C$ at a distance of at most $\frac{\Delta_r}{2}$ of $C$ after subtracting the ``in-place" demand from each point of $\overline{C} \setminus C$. The following lemma lower-bounds the cost $C(\bq, \bd)$.

\begin{restatable}{lemma}{restateLBqd}\label{lem:LBqd}
    Let $\bq \in \cQ$ and $\bd \geq \bzero$. Then,
    \begin{align*}
        C(\bq, \bd) 
        &\geq h\left(q_X - d_X\right) + (b+h) \left(d_{X} - q_X\right)^+
    \end{align*}
    and for every $r \in \{1, \dots, R-1\}$ and $\cF \in \families(\cP_r)$,
    \begin{align*}
        C(\bq, \bd) 
        &\geq h\left(q_X - d_X\right) + \frac{\Delta_r}{2} \sum_{C \in \cF} \left(d_{C} - \overline{q}_{C}\right)^+
    \end{align*}
\end{restatable}
The proof of the second inequality in Lemma~\ref{lem:LBqd} follows, roughly speaking, from the fact that for each cluster $C$, a least an amount of demand $(d_C-\overline{q}_C)^+$ given by the deviation of demand inside $C$ from the inventory around $C$ at a distance $\Delta_r/2$, needs to travel a distance of at least $\Delta_r/2$ to get fulfilled. 

We now show the following lemma which given two subsets $S \subset T$ of $X$, gives a lower-bound on the expected value of the deviation of the total demand within $S$ from a given quantity $Q$ subject to the demand being at least $\mu_i$ in every location of $T$.

\begin{restatable}{lemma}{restateExpBound}\label{lem:ExpBound}
    Given subsets $S \subset T \subset X$ and a quantity $Q \geq \mu_S$, it holds that
    \begin{align*}
        \sup_{\bd \in \cF_T(\bmu, \Sigma)} &\mathbb{E}_{\bd}\left((d_S - Q)^+ \cdot \mathbf{1}_{\{d_i \geq \mu_i, \forall i \in T\}}\right)= 
        \Omega \left(\sqrt{\sigma_S^2 + \left(Q - \mu_S\right)^2} - \left(Q - \mu_S\right)\right)
    \end{align*}
\end{restatable}

The proof of Lemma~\ref{lem:ExpBound} is via an explicit construction of a near worst-case distribution using Lemma~\ref{lem:chebyshev}. Combining Lemma~\ref{lem:LBqd} and Lemma~\ref{lem:ExpBound}, we prove the following lower-bound for the expected cost $C(\bq, \bd)$ under worst-case distribution of demand for every $\bq \in \cQ$,

\begin{restatable}{lemma}{restateLBzStar}\label{lem:LBzStar}
    Let $\bq \in \cQ$,
    \begin{align*}
        &\sigma_X \sqrt{bh} + \sum_{C \in \Gamma}\sigma_C \sqrt{b_Ch} +\sum_{C \in \clusters(H) \setminus \Gamma} \sigma_C\diam_p(C) \leq O(\gamma\alpha\beta R) \cdot\sup_{\cD \in \cF(\bmu, \Sigma)} \mathbb{E}_{\bd\sim \cD}\; C(\bq, \bd) 
    \end{align*}
\end{restatable}

Finally, by taking, in Lemma~\ref{lem:LBzStar}, $\bq$ to be the optimal solution of~\eqref{eq:offline} in $\cQ$, we get the lower-bound,

\begin{align}\label{eq:LB}
   & \sigma_X \sqrt{bh} + \sum_{C \in \Gamma}\sigma_C \sqrt{b_Ch} + \sum_{C \in \clusters(H) \setminus \Gamma} \sigma_C\diam_p(C) 
    \leq O(\gamma\alpha\beta R) \cdot z_{\textsf{ODRNM}}.
\end{align}
Lemma~\ref{lem:main1} follows from the upper-and lower-bounds~\eqref{eq:UB} and~\eqref{eq:LB}.

\section{The Hierarchical Balance (HB) Policy}
We now present an approximation policy for the online problem~\eqref{eq:online} where we observe the demand in a sequential manner and need to decide where it must be fulfilled from at the time it materializes. Our policy consists of an initial inventory allocation $\bq^{\textsf{HB}}$ and an online fulfillment policy $\cA^{\textsf{HB}}$. We refer to our policy as the {\em Hierarchical Balance} (HB) policy.

Let $\cH=\cP_1,\dots,\cP_{R}$ be an $(\alpha,\beta,\gamma)$-well-separated hierarchical partition of $(X, \ell)$ for some parameters $\alpha,\beta, \gamma$, such that $\gamma > \log(n)\alpha$. Given $r \in \{1, \dots, R-1\}$ and a cluster $C \in \clusters(\cP_r)$, let $\diam_{p}(C)$ denote the diameter of the parent cluster of $C$ in $\cP_{r+1}$. Let $$
\Gamma = \{C \in \clusters(\cH)\;|\; \diam_{p}(C) \geq 2h\}
$$
Then for every $C \in \Gamma$, let
$$
b_C = \diam_p(C)-h
$$

\vspace{3mm}{\noindent \bf Initial Inventory Allocation $\bq^{\textsf{HB}}$.} Our initial inventory allocation $\bq^{\textsf{HB}}$ is given by an optimal solution of the covering LP
\begin{align}
\label{eq:hfs}
    \min_{\bq} \;\;
    & h (q_X-\mu_X) \notag\\
    s.t. \;\; 
    & q_{C} - \mu_C \geq \frac{\sigma_C}{2}\left(\sqrt{\frac{b_C}{h}} - \sqrt{\frac{h}{b_C}}\right), \;\; \forall C \in \Gamma \notag\tag{HB}\\
    & q_i - \mu_i \geq \frac{\sigma_i}{\sum_{i\in X} \sigma_i} \frac{\sigma_X}{2}\left(\sqrt{\frac{b}{h}} - \sqrt{\frac{h}{b}}\right), \;\; \forall i \in X. \notag
\end{align}

\vspace{0mm}
{\noindent \bf Online fulfillment Policy $\cA^{\textsf{HB}}$.} 
The pseudo-code for the online fulfillment policy $\cA^{\textsf{HB}}$ is given in Algorithm~\ref{alg:online}. In this algorithm, after a part $\bd^t$ of the demand arrives. We fulfill it incrementally, one portion at a time, in an arbitrary order until all demand is fulfilled or no inventory remains (the while loop). In particular, at each step, we fulfill a small enough portion of demand $\delta d$ at some location $i$, which we will henceforth refer to as an infinitesimal portion, in the following way:
\begin{itemize}
    \item If there is inventory remaining at location $i$, the infinitesimal portion is fulfilled from this inventory in-place.
    \item Otherwise, let $r \in \{1, \dots, R\}$ be the smallest integer for which there exists remaining inventory in the cluster of $\cP_{r}$ to which $i$ belongs, say cluster $C$. Let $C_1, \dots, C_k$ be the clusters in $\clusters(C)$ that still have remaining inventory. We take an equal amount of $\frac{\delta d}{k}$ inventory from each one of these clusters. Within each cluster, an equal amount of inventory is taken from each subcluster that still has available inventory, and so on. If a cluster of $\cP_1$ is reached, take an equal amount of inventory from each of the locations that still have available inventory inside the cluster.
\end{itemize}
The infinitesimal portion $\delta d$ is chosen to be the minimum between the remaining demand at location $i$ and the maximum demand that can be fulfilled using the above procedure before one of the locations runs out of inventory. If we run out of inventory before fulfilling all the demand, all the subsequent demand will remain unfulfilled and an underage cost $b$ per unit is paid. Similarly, if at the end of the selling season we still have inventory, an overage cost $h$ per unit is paid.
\algtext*{EndIf}
\algtext*{EndFor}
\algtext*{EndWhile}
\begin{algorithm}
\caption{Online Fulfillment Policy $\cA^{\textsf{HB}}$}\label{alg:online}
\begin{algorithmic}[1]
\setlength{\itemsep}{0.3em}
\State Initialize current inventory vector $\bq^{\textsf{curr}} \gets \bq^{\textsf{HB}}$

\For{every time step $t$ before the end of the selling season}
    \State Receive demand vector $\mathbf{d}^t$
    \State Initialize current demand vector $\bd^{\textsf{curr}} \gets \bd^{t}$.
    \While{$d^{\textsf{curr}}_X > 0$ and $q_X^{\textsf{curr}} > 0$} \label{linestep}
        \State Let $i \in X$ such that $d^{\textsf{curr}}_i > 0$.
        \If{$q^{\textsf{curr}}_i > 0$} \label{lineif}
            \State Let $\delta d \gets \min(d^{\textsf{curr}}_i,\; q^{\textsf{curr}}_i)$
            \State Send $\delta d$ inventory from $i$ to $i$.
            \State Update $d_i^{\textsf{curr}} \gets d_i^{\textsf{curr}} - \delta d$ and $q_i^{\textsf{curr}} \gets q_i^{\textsf{curr}} - \delta d$
        \Else
            \State Let $r \gets \min \{ r' \in [R] \mid \exists \, C \in \clusters(\cP_{r'}) \text{ s.t. } i \in C, q^{\textsf{curr}}_C > 0\}$
            \State Let $C \in \clusters(\cP_r)$ \text{ s.t. } $i \in C,\; q^{\textsf{curr}}_C > 0$
            
            \For{$C' \in \cup_{r'=2}^r \clusters(\cP_{r'})$ \text{ s.t. } $C' \subset C$}
                \State Let $k_{C'} \gets |\{C'' \in \clusters(C') |\; q^{\textsf{curr}}_{C''} > 0\}|$
            \EndFor
            \For{$C' \in \clusters(\cP_{1})$ \text{ s.t. } $C' \subset C$}
                \State Let $k_{C'} \gets |\{j \in C' |\; q^{\textsf{curr}}_{j} > 0\}|$
            \EndFor
            
            \For{$j \in C$ s.t. $q^{\textsf{curr}}_j > 0$}
                \State Let $C_j^r,C_j^{r-1}, \dots, C^1_j$ be the clusters of levels $r, r-1, \dots, 1$ to which $j$ belongs
                \State Let $k_j = k_{C_j^r}\times k_{C_j^{r-1}} \times \dots \times k_{C_j^1}$
            \EndFor
            \State Let $\delta d = \min\{d^{\textsf{curr}}_i, \; \min_{j \in C} k_j \cdot q^{\textsf{curr}}_j\}$
            \For{$j \in C$ s.t. $q^{\textsf{curr}}_j > 0$}
                \State Send $\delta d/k_j$ inventory from $j$ to $i$
                \State Update $d_i^{\textsf{curr}} \gets d_i^{\textsf{curr}} - \delta d/k_j$ and $q_j^{\textsf{curr}} \gets q_j^{\textsf{curr}} - \delta d /k_j$
            \EndFor
        \EndIf
    \EndWhile
    \If{$d^{\textsf{curr}}_X > 0$} \Comment{We pay underage cost if remaining demand}
    \State Pay underage cost $b \cdot d^{\textsf{curr}}_X$
    \EndIf
\EndFor
\If{$q^{\textsf{curr}}_X > 0$} \Comment{We pay overage cost if remaining inventory}
    \State Pay overage cost $h \cdot q^{\textsf{curr}}_X$
\EndIf
\end{algorithmic}
\end{algorithm}

It is a priori unclear that $\cA^{\textsf{HB}}$ can be implemented in polynomial time. The following lemma shows that this is indeed the case.

\begin{restatable}{lemma}{restateInputSize}\label{lem:InputSize}
    $\cA^{\textsf{HB}}$ can be implemented in polynomial time.
\end{restatable}

We show next that our policy gives a good approximation to~\eqref{eq:online}. In particular, we show the following lemma,
\begin{restatable}{lemma}{restateOnline}\label{lem:main2}
    The initial inventory allocation $\bq^{\textsf{HB}}$ and the online fulfillment policy $\cA^{\textsf{HB}}$ are such that,
    \begin{align*}
        \sup_{\cD \in \cF(\bmu, \Sigma)} &\mathbb{E}_{\bd \sim \cD}\; \left(\sup_{\pi \in \Pi(\bd)}C_{\cA^{\textsf{HB}}}(\bq^{\textsf{HB}}, \pi)\right)
        \leq O(\alpha\beta\gamma R \log(n)) \cdot z_{\textsf{ODRNM}}  \leq O(\alpha\beta\gamma R \log(n)) \cdot z_{\textsf{DRNM}}
    \end{align*}
\end{restatable}

The above lemma shows that our approximation bounds are not only with respect to the optimum of~\eqref{eq:online}, but also with respect to the optimum of the offline problem~\eqref{eq:offline}. 
Lemma~\ref{lem:wshp} and Lemma~\ref{lem:main2} imply our second main theorem stated formally as follows,
\begin{theorem}\label{thm:main2}
    There exists an initial inventory allocation and an online fulfillment policy that is an
    \begin{itemize}
        \item $O(\log^2(n))$ approximation of~\eqref{eq:online} when $(X, \ell)$ is a uniform metric space.
        \item $O(\log^3(n))$ approximation of~\eqref{eq:online} when $(X, \ell)$ is a subset of an Euclidean space $\mathbb{R}^d$ of constant dimension $d$.
        \item $O(\log^6(n))$ approximation of~\eqref{eq:online} when $(X, \ell)$ is any general metric space.
    \end{itemize}
\end{theorem}
The next section is devoted to the proof of Lemma~\ref{lem:main2}. We give an outline of the proof in the main part and differ the full proofs to the Appendix.

\section{Analysis of the Hierarchical Balance (HB) Policy}
\label{sec:HBanalysis}
Recall that,
\begin{align*}
    C^H(\bq, \bd) &= 
    h(q_X-d_X)^+ + b(d_X-q_X)^+ 
    \\
    &\quad +\sum_{r=1}^{R-1} \sum_{C \in \clusters(\cP_{r})} \diam_p(C) \cdot (d_C - q_C)^+ + \frac{\diam(X)}{n}\sum_{i \in X}(d_i-q_i)^+.
\end{align*}
for every $\bq$ and $\bd$. Define,
\begin{align*}
    C^H_{\cA^{\textsf{HB}}}(\bq, \pi) 
    &= h\sum_{i \in X}\left(q_i - \sum_{j\in X} x^{\cA^{\textsf{HB}}}_{ij}(\bq, \pi)\right) + b \sum_{j \in X} \left(d_j-\sum_{i \in X} x^{\cA^{\textsf{HB}}}_{ij}(\bq, \pi)\right) + \sum_{i \in X}\sum_{j \in X} \ell^H_{ij} x^{\cA^{\textsf{HB}}}_{ij}(\bq, \pi)
\end{align*}
for every $\bq$ and $\pi$, where for every $i \neq j \in X$, and letting $r \in \{1, \dots, R\}$ be the smallest integer such that $i$ and $j$ belong to the same cluster $C \in \clusters(\cP_r)$, we set $\ell^H_{ij} = \diam(C)$. We set $\ell^H_{ii} = 0$ for every $i \in X$. Note that $\ell^H_{ij} \geq \ell_{ij}$ for every $i, j \in X$. Hence, 
$
C_{\cA^{\textsf{HB}}}(\bq, \pi) \leq C^H_{\cA^{\textsf{HB}}}(\bq, \pi).
$
for every $\bq$ and $\pi$.

The following lemma shows that it is sufficient to compare $C^H_{\cA^{\textsf{HB}}}(\bq^{\textsf{HB}}, \pi)$ and $C^H(\bq^{\textsf{HB}}, \bd)$, for demand vectors of same total quantity as $\bq^{\textsf{HB}}$.

\begin{restatable}{lemma}{restateReductionToC}\label{lem:reductionToC}
    Assume
    $$
    C^H_{\cA^{\textsf{HB}}}(\bq^{\textsf{HB}}, \pi) \leq O(\log(n)) \cdot C^H(\bq^{\textsf{HB}}, \bd) 
    $$
    for every $\bd \geq \bzero$ such that $d_X = q^{\textsf{HB}}_X$ and $\pi \in \Pi(\bd)$, then,
    $\bq^{\textsf{HB}}$ and $\cA^{\textsf{HB}}$ achieve the approximation bound in Lemma~\ref{lem:main2}.
\end{restatable}

We now fix $\bd \geq \bzero$ such that $d_X = q^{\textsf{HB}}_X$ and $\pi \in \Pi(\bd)$, and prove that $$C^H_{\cA^{\textsf{HB}}}(\bq^{\textsf{HB}}, \pi) \leq O(\log(n)) \cdot C^H(\bq^{\textsf{HB}}, \bd).$$ Note that when $d_X = q^{\textsf{HB}}_X$, there is no overage or underage and the whole demand will be fulfilled by the online algorithm.

Consider some cluster $C \in \clusters(\cH)$. When an infinitesimal portion of demand $\delta d$ is considered, there are two cases, either the portion arrives inside $C$ and in which case it is either entirely fulfilled from inventory inside $C$ or entirely fulfilled from inventory outside $C$ by construction of the algorithm, or the portion arrives outside $C$ in which case it might happen that a part of it gets fulfilled from inventory inside $C$. Let $\delta d^1, \dots, \delta d^L$ denote all of the infinitesimal demand portions or parts of infinitesimal demand portions that get fulfilled using inventory inside $C$, in the order they were considered by the online algorithm, and suppose they belong to location $i^1, \dots, i^L$ respectively. Note that because all demand is fulfilled with no remaining inventory (recall that $d_X = q^{\textsf{HB}}_X$), the total demand given by all the portions $\sum_{l=1}^L \delta d^l$ is equal to the total inventory inside $C$, i.e., $q^{\textsf{HB}}_C$.

If $C \in \cP_r$ for $r\geq 2$, define $V_C=\{C(0)\} \cup \clusters(C)$ as the set of clusters of $\cP_{r-1}$ that form $C$ plus a ``virtual cluster" $C(0)$. If $C \in \cP_1$, define $V_C=\{C(0)\} \cup C$ as the set of locations inside $C$ plus a ``virtual cluster" $C(0)$. The virtual cluster is introduced for simplicity, it will be considered in the following as a cluster inside $C$ with $q=0$ inventory inside it and such that whenever some $\delta d^l$ arrives outside of $C$ and is fulfilled from inventory inside $C$, we say that it arrived inside $C(0)$ and fulfilled by inventory outside $C(0)$.

For every $C' \in V_C$, let 
$
\theta^{\textsf{online}, C}(C') 
$
denote the total amount of demand from $\delta d^1, \dots, \delta d^{L}$ that arrived inside $C'$ and that the online algorithm fulfilled using inventory outside of $C'$ (but inside $C$ by definition of $\delta d^1, \dots, \delta d^{L}$).

Next, consider an offline algorithm that fulfills the portions $\delta d^{1}, \dots, \delta d^{L}$ inside $C$ in a way that fulfills a maximum amount of demand to the remaining inventory inside each $C' \in V_C$, then fulfills the remaining demand arbitrary inside $C$. Then let,
$
\theta^{\textsf{offline}, C}(C') 
$
denote the amount of demand from $\delta d^{1}, \dots, \delta d^{L}$ that this offline algorithm fulfills outside of $C' \in V_C$. Note that,
$
\theta^{\textsf{offline}, C}(C') = \left(\sum_{l: i^l \in C'} \delta d^l - q^{\textsf{HB}}_{C'}\right)^+.
$

The following lemma upper-bounds the total amount of demand that the online algorithm fulfills from outside the cluster/location it arrived to, as a function of the total amount the offline algorithm fulfills from outside of the cluster/location it arrived to. 
\begin{restatable}{lemma}{restateComparePotential}\label{lem:comparepotential}
The following holds,
    $$
    \sum_{C' \in V_C} \theta^{\textsf{online}, C}(C') \leq 3\ln(n) \sum_{C' \in V_C} \theta^{\textsf{offline}, C}(C')
    $$
\end{restatable}

To prove Lemma~\ref{lem:comparepotential}, we define for every $C' \in V_C$ and $l \in \{1, \dots, L+1\}$,
$
\theta^{\textsf{online}, C}_{[1:l-1]}(C') 
$
as the total amount of demand from $\delta d^1, \dots, \delta d^{l-1}$ that arrived inside $C'$ and that the online algorithm fulfilled using inventory outside of $C'$, where by convention $
\theta^{\textsf{online}, C}_{[1:0]}(C') = 0
$.

Next, consider the inventory remaining after the online algorithm fulfilled $\delta d^1, \dots, \delta d^{l-1}$. Consider an offline algorithm that fulfills the remaining portions $\delta d^{l}, \dots, \delta d^{L}$ inside $C$ using this remaining inventory in a way that fulfills a maximum amount of demand to the remaining inventory inside each $C' \in V_C$, then fulfills the remaining demand arbitrary inside $C$. We define for every $C' \in V_C$ and $l \in \{1, \dots, L+1\}$,
$
\theta^{\textsf{offline}, C}_{[l:L]}(C') 
$
as the amount of demand from $\delta d^{l}, \dots, \delta d^{L}$ that the offline algorithm fulfilled outside of $C' \in V_C$, where by convention $
\theta^{\textsf{offline}, C}_{[L+1:L]}(C') = 0
$.

Let $C(0), C(1), \dots, C(U)$ denote the order by which the elements of $V_C$ run out of inventory (the virtual cluster $C(0)$ has no inventory and hence its first in the list). We consider the following potential function
$$
\Phi_l =  \sum_{u=0}^U \theta^{\textsf{online}, C}_{[1:l-1]}(C(u)) + \sum_{u=0}^U 3\ln(U-u+1) \cdot \theta^{\textsf{offline}, C}_{[l:L]}(C(u)),
$$
for every $l \in \{1,\dots, L+1\}$, and show that $\Phi_l$ is non-increasing in $l$. The lemma follows from the fact that 
$$
    \Phi_{L+1} = \sum_{C' \in V_C} \theta^{\textsf{online}, C}_{[1:L]}(C') = \sum_{C' \in V_C} \theta^{\textsf{online}, C}(C') \leq \Phi_{1}
$$
and
$$
\Phi_{1} \leq 3\ln(n) \sum_{C' \in V_C} \theta^{\textsf{offline}, C}_{[1:L]}(C') = 3\ln(n) \sum_{C' \in V_C} \theta^{\textsf{offline}, C}(C').
$$
Now, the cost $C^H_{\cA^{\textsf{HB}}}(\bq^{\textsf{HB}}, \pi)$ can be upper-bounded in function of the $\theta^{\textsf{online}, C}(C')$'s as follows,
\begin{restatable}{lemma}{restateOnlineUB}\label{lem:onlineUB}
The cost $C^H_{\cA^{\textsf{HB}}}(\bq^{\textsf{HB}}, \pi)$ can be upper-bounded as follows,
    $$
    C^H_{\cA^{\textsf{HB}}}(\bq^{\textsf{HB}}, \pi) \leq \sum_{C \in \clusters(\cH)} \diam(C) \sum_{C' \in V_C}\theta^{\textsf{online}, C}(C')
 $$
\end{restatable}

Then, the cost $C^H(\bq^{\textsf{HB}}, \bd)$ can be lower-bounded as a function of the $\theta^{\textsf{online}, C}(C')$'s and $\theta^{\textsf{offline}, C}(C')$'s as follows,
\begin{restatable}{lemma}{restateOnlineLB}\label{lem:onlineLB}
The cost $C^H(\bq^{\textsf{HB}}, \bd)$ can be lower-bounded as follows,
\begin{align*}
    C^H(\bq^{\textsf{HB}}, \bd) & \geq \frac{1}{2}\sum_{C \in \clusters(\cH)} \diam(C) \sum_{C' \in V_C}\theta^{\textsf{offline}, C}(C')
    \\
    &\quad - \frac{1}{21\log(n)}\sum_{C \in \clusters(\cH)} \diam(C) \sum_{C' \in V_C}\theta^{\textsf{online}, C}(C')
\end{align*}
\end{restatable}

Finally from Lemma~\ref{lem:comparepotential} and Lemma~\ref{lem:onlineLB}, we conclude that,
$$
C^H(\bq^{\textsf{HB}}, \bd)\geq \Omega\left(\frac{1}{\log(n)}\right)\sum_{C \in \clusters(\cH)} \diam(C) \sum_{C' \in V_C}\theta^{\textsf{online}, C}(C')
$$
combined with Lemma~\ref{lem:onlineUB}, we get,
\begin{align*}
    C^H(\bq^{\textsf{HB}}, \bd) 
    &= \Omega\left(\frac{1}{\log(n)}\right)C^H_{\cA^{\textsf{HB}}}(\bq^{\textsf{HB}}, \pi)
\end{align*}
which concludes the proof.

\section{Numerical Experiments}
In this section, we evaluate the numerical performance of our policy. Our results demonstrate that our policy not only provides strong theoretical guarantees but also performs effectively in practice. We conduct two sets of experiments: the first assesses the performance of our Generalized-Scarf-on-a-Metric (GSM) solution for the offline problem, and the second evaluates the performance of the online fulfillment policy $\cA^{\textsf{HB}}$ of Hierarchical Balance (HB).

\subsection{The Generalized-Scarf-on-a-Metric (GSM) Solution} \label{sec:num}

To evaluate the numerical performance of our Generalized-Scarf-on-a-Metric (GSM) Solution, we adopt an approach similar to that used in \cite{govindarajan2021distribution}
to evaluate their heuristic and consider special metric spaces where the cost function $C(\bq, \bd)$ has a simple expression in terms of $\bq$ and $\bd$ and which are a good approximations of general metric spaces in practice, and consider the slightly modified problem where demand is allowed to take negative values. These adjustments enable an exact reformulation of the problem as an exponentially sized SDP (that can still be solved for a small number of locations). We then compare our approximation solution to an optimal solution of the SDP. We note that no algorithms (even inefficient ones) are known for solving the problem exactly otherwise.

Specifically, we consider metric spaces where $X$ is given by the leaves of a rooted tree $T$, and $\ell_{ij}$ represents the (shortest path) distance between leaves $i$ and $j$ in the tree. Let $L_1, \dots, L_K$ denote the vertices at levels $1, 2, \dots, K$ of the tree $T$, with level $1$ representing the leaves and level $K$ the root. The edge weights of $T$ have the following structure: the weight between leaves and their parents at level $2$ is $c \lambda$ for some $c > 0$ and $\lambda > 1$; between nodes at level $2$ and their parents at level $3$, it is $c (\lambda^2 - \lambda)$; between nodes at level $3$ and their parents at level $4$, it is $c (\lambda^3 - \lambda^2)$; and so on.
In practice, locations $X$ often exhibit hierarchically clustered structures similar to these tree metrics. For example, locations might be divided into cities that are very far apart, within which neighborhoods are far apart but not as distant, and so on.
It is also known that any metric space can be randomly embedded into trees of this form with low distortion~\cite{fakcharoenphol2003tight}, further suggesting that a good performance on these metric spaces suggests a good performance on general metric spaces.

The metric spaces constructed above admit a natural well-separated hierarchical partition defined as follows: let $r_0$ be the smallest level bellow which the nodes of $L_{r}$, for $r \leq r_0$ all have a single child. Then, let $\cH =\cP_{1}, \dots, \cP_{K-r_0+1}$, where each $\cP_{r}$ consists of a single family with clusters $\{C_{v}\}_{v \in L_{r+r_0-1}}$, where $C_{v}$ is the set of leaves that are children of $v \in L_{r+r_0-1}$. This is indeed an $(\alpha, \beta, \gamma)$-well-separated hierarchical partition for $\alpha = (1+\epsilon)^{K-r_0}$, $\beta=1$ and $\gamma=(1+\epsilon)\lambda$, for every $\epsilon>0$ small enough. In fact, (i) it is clearly a coarsening family of partitions. (ii) for $n$ large enough $\Delta_1 = \min_{i \neq j \in X} \ell_{ij}/\alpha = 2c\lambda^{r_0}/(1+\epsilon)^{K-r_0}$ and hence $\Delta_{r} = 2c\lambda^{r_0+r-1}/(1+\epsilon)^{K-r_0}$ and $R=\lceil\log(\alpha \lambda^{K-r_0})/\log(\gamma)\rceil + 1 = K-r_0+1$. (iii) each $\cP_r$ is an $(\alpha, \beta)$-well separated partition of margin $\Delta_r$ as the distance between any two clusters is $2c\lambda^{r+r_0-1}$ which is strictly larger than $\Delta_r$, and the diameter of each cluster is at most $2c\lambda^{r+r_0-2}$ which is strictly less than $\Delta_r$. We consider these well separated hierarchical partitions in our experiments (for simplicity we let $\epsilon \rightarrow 0$) and let 
$\bq^{\textsf{GSM}}$ denote the corresponding Generalized-Scarf-on-a-Metric (GSM) solution.

We compare our solution $\bq^{\textsf{GSM}}$ to an optimal solution of~\eqref{SDP} which we denote by $\bq^{\textsf{SDP}}$. As mentioned earlier,~\eqref{SDP} can be obtained by leveraging the explicit expression of the cost function $C(\bq, \bd)$ and assuming demand can take negative values. The derivation is given in Appendix~\ref{apx:derivation-sdp} for completeness.
\begin{align}\label{SDP}
    \min_{\substack{Y, \by, y_0\\\bq \geq \bzero}} \quad 
    &\sum_{i,j \in X} (\mu_i \mu_j + \Sigma_{ij}) Y_{ij} + \sum_{i \in X} \mu_i y_i + y_0 \notag\\
    s.t. \quad 
    & 
    \begin{bmatrix}
        Y & \frac{\by - f( \bm{\epsilon}_1, \dots, \bm{\epsilon}_K)}{2}\\
        \frac{\by\trsp - f( \bm{\epsilon}_1, \dots, \bm{\epsilon}_K)\trsp}{2} & y_0 + f( \bm{\epsilon}_1, \dots, \bm{\epsilon}_K)\trsp \bq
    \end{bmatrix}\succeq 0, \tag{SDP}
    \\
    &\quad \forall \bm{\epsilon}_1 \in \Pi_1, \dots, \forall \bm{\epsilon}_K \in \Pi_K \notag.
\end{align}

To compare the performance of the SDP solution $\bq^{\textsf{SDP}}$ and that of our solution $\bq^{\textsf{GSM}}$, we compare their (sample) average 
cost across multiple choices of the ground-truth demand distribution, specifically Normal, Log-Normal, and Gamma. In particular, we solve the following sample average approximation LP,
$$
\begin{array}{rcl}
\min\limits_{\{\bx^{\bd}\}_{\bd \in D}} & \quad & \displaystyle\sum_{\bd \in D} \Biggl(
h\sum\limits_{i \in X}\Bigl(q_i - \sum\limits_{j \in X} x^{\bd}_{ij}\Bigr)
+ b\sum\limits_{j \in X}\Bigl(d_j - \sum\limits_{i \in X} x^{\bd}_{ij}\Bigr)
+ \sum\limits_{i \in X}\sum\limits_{j \in X} \ell_{ij}\, x^{\bd}_{ij}
\Biggr) \\[1ex]
\text{s.t.} & \quad & \displaystyle\sum\limits_{j \in X} x^{\bd}_{ij} \le q_i,\quad \forall\, i \in X,\ \forall\, \bd \in D, \\[1ex]
           & \quad & \displaystyle\sum\limits_{i \in X} x^{\bd}_{ij} \le d_j,\quad \forall\, j \in X,\ \forall\, \bd \in D, \\[1ex]
           & \quad & x^{\bd}_{ij} \ge 0,\quad \forall\, i,j \in X,\ \forall\, \bd \in D.
\end{array}
$$

\noindent for $\bq \in \{\bq^{\textsf{SDP}}, \bq^{\textsf{GSM}}\}$ and for $D \in \{D^{\textsf{normal}}, D^{\textsf{log-normal}}, D^{\textsf{gamma}}\}$ where each $D^{\textsf{normal}}, D^{\textsf{log-normal}}, D^{\textsf{gamma}}$ is a set of $m$ samples from a Normal, Log-Normal and Gamma distribution respectively of mean $\bmu$ and covariance $\Sigma$. This setup ensures that we evaluate our solutions against demand distributions with different tail behaviors.

\begin{figure}[htbp]
\centering
\subcaptionbox{$K=2$\label{fig:offline3}}{%
  \includegraphics[width=0.3\textwidth, trim=0 0 0 19, clip]{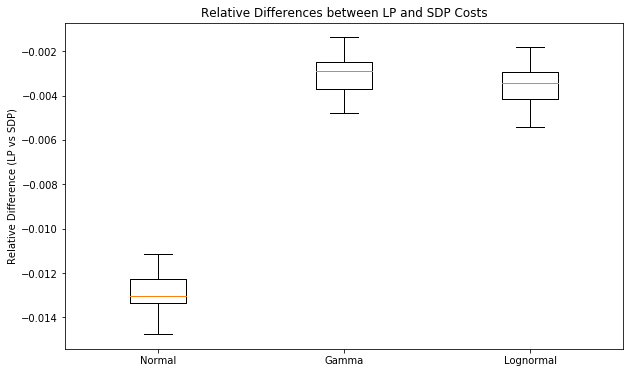}%
}
\hfill
\subcaptionbox{$K=3$\label{fig:offline2}}{%
  \includegraphics[width=0.3\textwidth, trim=0 0 0 19, clip]{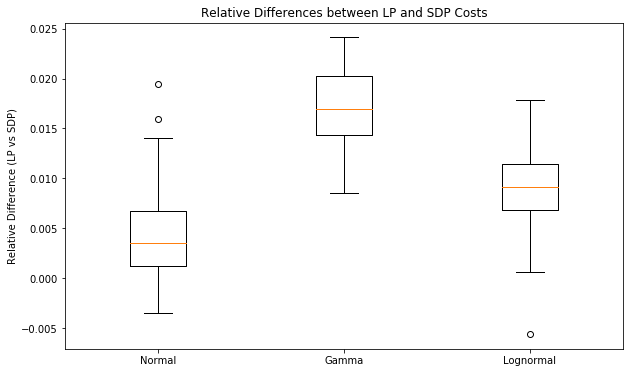}%
}
\hfill
\subcaptionbox{$K=4$\label{fig:offline1}}{%
  \includegraphics[width=0.3\textwidth, trim=0 0 0 19, clip]{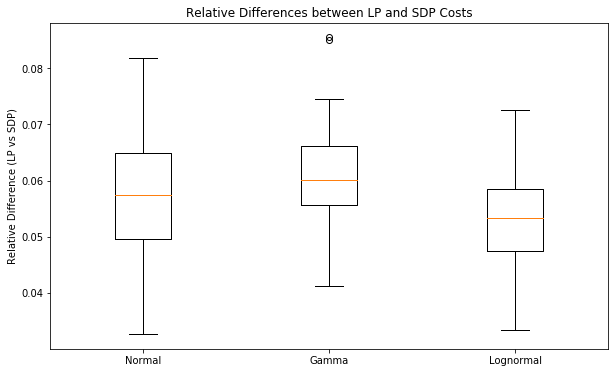}%
}
\caption{The relative difference between the optimal cost of the SAA linear program under the SDP solution $\bq^{\textsf{SDP}}$ and the LP solution $\bq^{\textsf{GSM}}$ for different ground-truth demand distributions.}
\label{fig:offline}
\end{figure}

\begin{figure}[htbp]
\centering
\subcaptionbox{$K=2$\label{fig:tree_offline_3}}{%
  \includegraphics[width=0.3\textwidth, trim=0 0 0 35, clip]{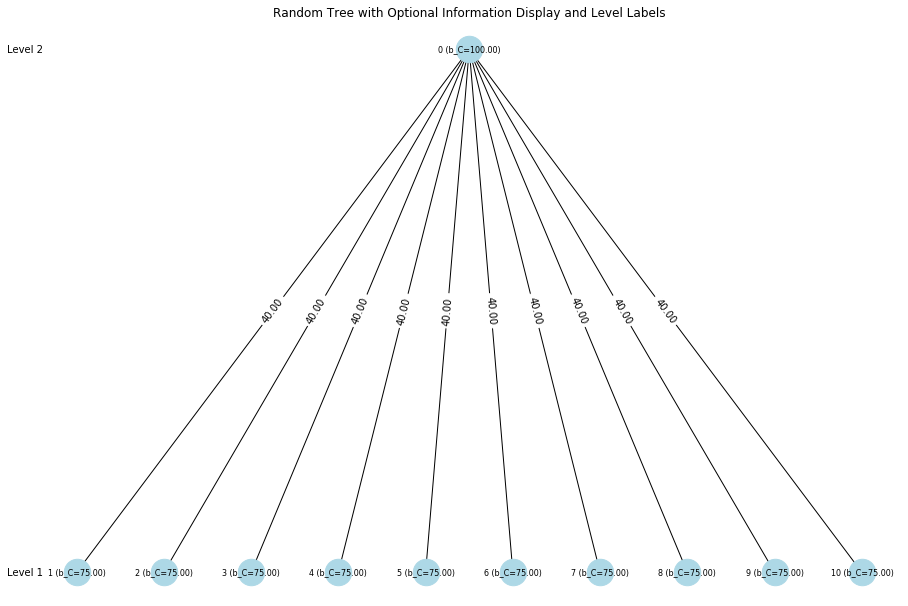}%
}
\hfill
\subcaptionbox{$K=3$\label{fig:tree_offline_2}}{%
  \includegraphics[width=0.3\textwidth, trim=0 0 0 35, clip]{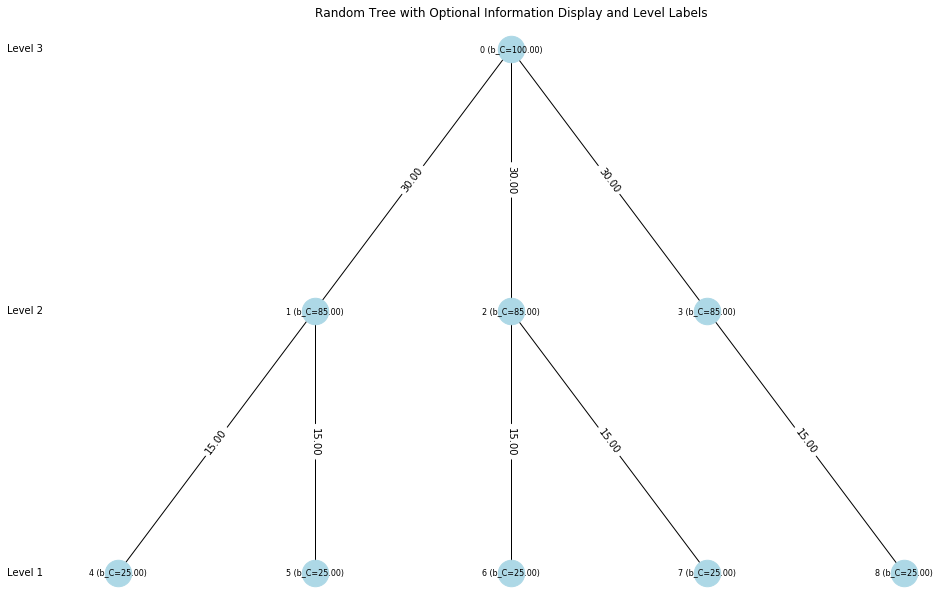}%
}
\hfill
\subcaptionbox{$K=4$\label{fig:tree_offline_1}}{%
  \includegraphics[width=0.3\textwidth, trim=0 0 0 35, clip]{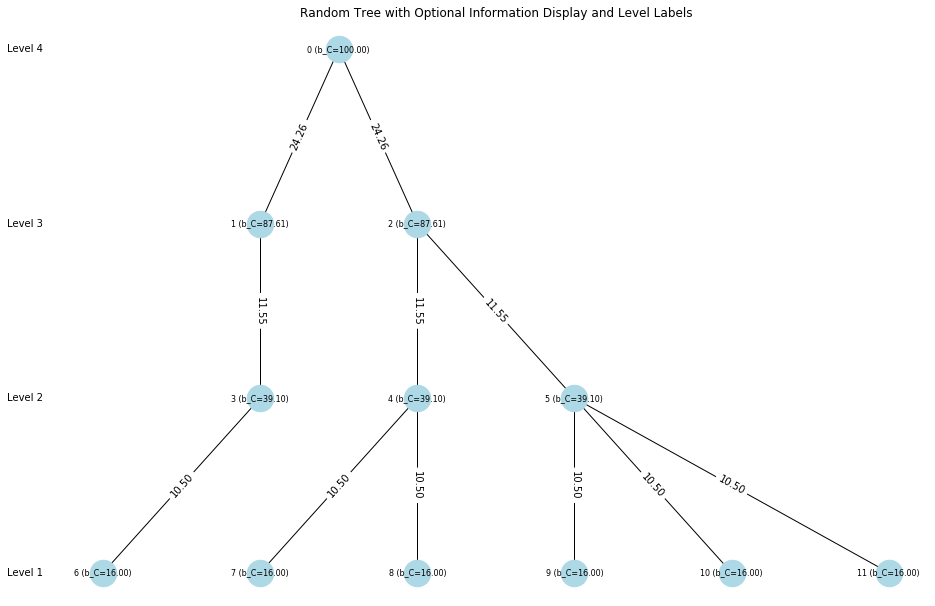}%
}
\caption{Trees used to construct the metric space in each experiment.}
\label{fig:tree_offline}
\end{figure}

\vspace{3mm}{\noindent \bf Experimental Setup.}
We set $b = 100$, $h = 5$, and $m = 1000$, and conduct our tests on three different metric spaces given by tree metrics constructed as described above, with the number of levels $K = 2$, $3$, and $4$. Note that when $K = 2$, the resulting metric space corresponds simply to a uniform metric space where the distance between every two locations is the same. Figure~\ref{fig:offline} shows the relative difference $\frac{{\textsf{cost}}(\bq^{\textsf{GSM}}) - {\textsf{cost}}(\bq^{\textsf{SDP}})}{ {\textsf{cost}}(\bq^{\textsf{SDP}})}$ obtained from $N = 50$ repetitions of the experiments. Figure~\ref{fig:tree_offline} shows the trees used to construct the metric space in each experiment. The numbers on the edges indicate the edge weights. The mean values $\bmu$ are randomly chosen from $[200, 1500]$, and for each $i \in X$, the standard deviation $\sigma_i$ is sampled randomly not too large compared to the mean (uniformly sampled from $[0.3, 0.8] \times \mu_i$).

\vspace{3mm}{\noindent \bf Results.} We observe from Figure~\ref{fig:offline} that our solution $\bq^{\textsf{GSM}}$ performs well across all the instances we consider and under all the three different distributions, with a relative gap of at most $6\%$ compared to the SDP solution. In the uniform metric space, our solution is slightly better (though not by much, less than $0.3\%$). These results suggest that our solution, despite its simplicity, effectively captures the key structure of the problem. As such, it provides a strong candidate to use in practice for the offline problem, offering both solid performance and interpretability.

\subsection{The Online Fulfillment Policy}

We evaluate the numerical performance of the online fulfillment policy $\cA^{\textsf{HB}}$ of Hierarchical Balance (HB) by comparing it to fulfilling demand offline. We show that even against this ``strong" benchmark, our policy still performs well. As there is no longer a need for a simple closed-form expression for the cost function $C(\bq, \bd)$, we can now consider more general metric spaces, in particular, we will consider points on the Euclidean plane. We compute the initial inventory allocation $\bq^{\textsf{HB}}$ and then compare the (sample) average cost of fulfilling demand offline and online, under different ground-truth distributions (Normal, Log-Normal, and Gamma).

\begin{figure}[htbp]
  \centering
  \subcaptionbox{$n=10$\label{fig:online10}}{
    \includegraphics[width=0.45\textwidth, trim=0 0 0 19, clip]{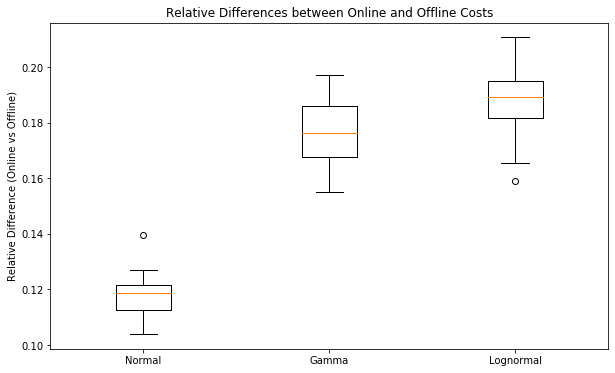}
  }
  \hfill
  \subcaptionbox{$n=15$\label{fig:online15}}{
    \includegraphics[width=0.45\textwidth, trim=0 0 0 19, clip]{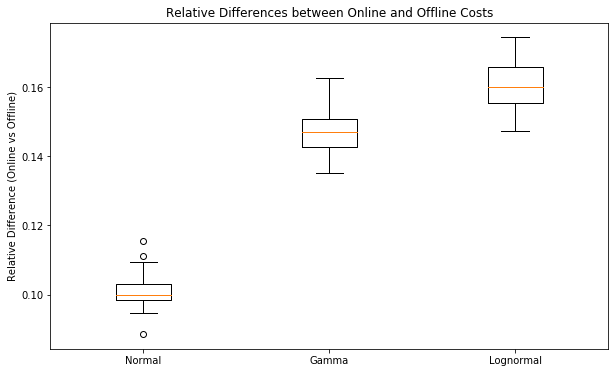}
  }\\[1ex]
  \subcaptionbox{$n=20$\label{fig:online20}}{
    \includegraphics[width=0.45\textwidth, trim=0 0 0 19, clip]{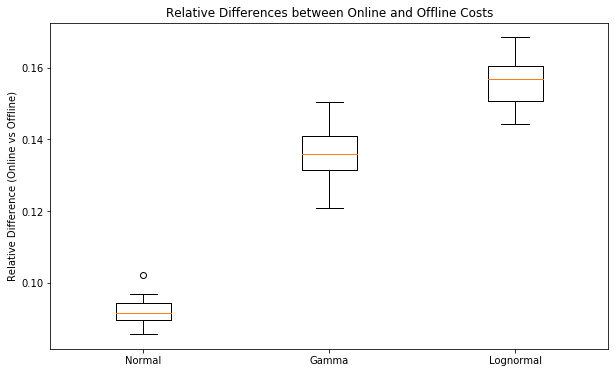}
  }
  \hfill
  \subcaptionbox{$n=25$\label{fig:online25}}{%
    \includegraphics[width=0.45\textwidth, trim=0 0 0 19, clip]{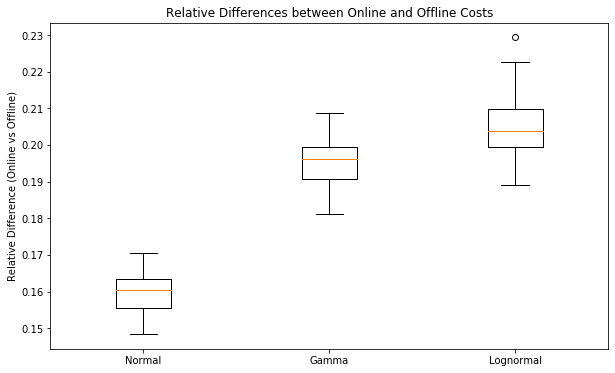}
  }
  \caption{The relative difference between the (sample) average cost of fulfilling demand offline and online when demand arrives one location at a time, for different values of $n$ and across various ground-truth demand distributions.}
  \label{fig:online}
\end{figure}

\vspace{3mm}{\noindent \bf Experimental Setup.} 
We choose $b = 100$, $h = 5$, and $m = 1000$, and generate $n = 10, 15, 20, 25$ random points in the Euclidean plane. We use the Well-Separated Hierarchical Partitions as constructed in Lemma~\ref{lem:wshp}. Figure~\ref{fig:online} shows the relative difference $\frac{{\textsf{cost}}(\mathrm{online}) - {\textsf{cost}}(\mathrm{offline})}{ {\textsf{cost}}(\mathrm{offline})}$ obtained from $N = 50$ repetitions of the experiments. The mean values $\bmu$ are randomly sampled from $[200, 1500]$, and for each $i \in X$, the standard deviation $\sigma_i$ is sampled randomly not too large compared to the mean (uniformly sampled from $[0.3, 0.8] \times \mu_i$).

\vspace{3mm}{\noindent \bf Results.} We observe from Figure~\ref{fig:online} that even against the strong benchmark of the offline optimum, the relative gap of our online fulfillment policy does not exceed $20\%$ across all instances considered. These results suggest that our online fulfillment policy is a strong candidate for practical use in real-time demand fulfillment.

\printbibliography

\newpage
{\Huge \noindent Appendix}

\section{Generalization to Demand Locations Different from Warehouses.}
\label{apx:generallocations}

We note that our simplified model subsumes the more general model where demand locations are different from warehouses: Given an instance of the general model, one can form an instance of the simplified model where the demand mean and variance at a warehouse are the total mean and variance from all demand locations for which it is the closest warehouse, and where demand requests arriving at a demand location are assumed to arrive at the nearest warehouse for that location. Given $\bq^{\sf HB}$ and $\cA^{\sf HB}$—the inventory solution and fulfillment policy from our Hierarchical Balance (HB) policy for the simplified model instance—in the general model instance, $\bq^{\sf HB}$ is used as the inventory solution; for the fulfillment policy, demand arriving at a demand location is routed to its nearest warehouse, and $\cA^{\sf HB}$ is then used to fulfill this demand. Following a similar analysis to that of our policy in the simplified model, it can be shown that this modified algorithm for the general model achieves the same approximation guarantees.

% Using triangle inequality and noting that all demand distributions we construct in our analysis rely on Lemma~\ref{lem:chebyshev} which still holds under the general model 

% if the set of distributions $\cF_{S}(\bmu, \Sigma)$ where $\mu_i$ (resp. $\Sigma_{ii}$) is the total mean (resp. variance) of demand from locations whose nearest warehouse is $i$ is replaced with the subset of distributions we get from routing the demand at each demand location to its nearest warehouse in the more general model, the same approximation guarantees can be proven under the general model by following the same analysis of our algorithms.

\section{Robustness to misspecification in the demand mean and variance.}\label{apx:misspecification}

Assume upper-bound estimates $\hat{\bmu}$ and $\hat{\bsigma}$ of the true demand mean $\bmu$ and standard deviations $\bsigma$ are used to compute our solutions. For example, one might use high probability upper-confidence bounds estimated from data. The following lemma relates the loss in optimality to the error in the estimates.

\begin{lemma}
    Let $\hat{\bq}^{\sf GSM}$ (resp. $\hat{\bq}^{\sf HB}$) be the inventory solution given by the Generalized-Scarf-on-a-Metric (GSM) solution (resp. the Hierarchical Balance policy (HB)) under upper-bound estimates $\hat{\bmu}$ and $\hat{\bsigma}$ of the true demand mean $\bmu$ and standard deviations $\bsigma$, and let ${\bq}^{\sf GSM}$ (resp. ${\bq}^{\sf HB}$) the solution under the true parameters. Assume $\frac{\hat{\mu_i}}{\mu_i}, \frac{\hat{\sigma_i}}{\sigma_i} \leq 1+\delta$ for every $i \in X$.
    Then,
    $$
    \sup_{\cD \in \cF(\bmu, \Sigma)} \mathbb{E}_{\bd \sim \cD}\; \left(C(\hat{\bq}^{\textsf{GSM}}, \bd)\right) \leq O((1+\delta)\alpha\beta\gamma R)\cdot z_{\sf ODRNM} + hq^{\textsf{GSM}}_X\cdot \delta,
    $$
    resp.
    $$
    \sup_{\cD \in \cF(\bmu, \Sigma)} \mathbb{E}_{\bd \sim \cD}\; \left(\sup_{\pi \in \Pi(\bd)}C_{\cA^{\textsf{HB}}}(\hat{\bq}^{\textsf{HB}}, \pi)\right) 
    \leq 
     O((1+\delta)\alpha\beta\gamma R \log(n))\cdot z_{\sf DRNM} + O(\log(n))hq^{\textsf{HB}}_X\cdot \delta.
    $$
\end{lemma}
In particular, in both~\eqref{eq:offline} and~\eqref{eq:online}, the extra loss from using the upper-bound estimates $\hat{\bmu}$ and $\hat{\bsigma}$ instead of the true parameters $\bmu$ and $\bsigma$ vanishes as the estimates become more accurate.

\begin{proof}
    First, recall that by Lemma~\ref{lem:UBHierarchicalMatching},
    $$
    \sup_{\cD \in \cF(\bmu, \Sigma)} \mathbb{E}_{\bd \sim \cD}\; \left(
        C(\hat{\bq}^{\textsf{GSM}}, \bd)
        \right)
        \leq 
        \sup_{\cD \in \cF(\bmu, \Sigma)} \mathbb{E}_{\bd \sim \cD}\; \left(
        C^H(\hat{\bq}^{\textsf{GSM}}, \bd)
        \right).
    $$
    Also, following the lines of the analysis in Section~\ref{sec:HBanalysis}, it can be shown that,
    $$
    \sup_{\cD \in \cF(\bmu, \Sigma)} \mathbb{E}_{\bd \sim \cD}\; \left(\sup_{\pi \in \Pi(\bd)}
        C_{\cA^{\textsf{HB}}}(\hat{\bq}^{\textsf{HB}}, \pi)
        \right)
        \leq 
        O(\log(n))\sup_{\cD \in \cF(\bmu, \Sigma)} \mathbb{E}_{\bd \sim \cD}\; \left(
        C^H(\hat{\bq}^{\textsf{HB}}, \bd)
        \right)
    $$
    It is therefore sufficient to show that,
    $$
    \sup_{\cD \in \cF(\bmu, \Sigma)} \mathbb{E}_{\bd \sim \cD}\; \left(
        C^H(\hat{\bq}^{\textsf{GSM}}, \bd)
        \right) \leq O((1+\delta)\alpha\beta\gamma R) \cdot z_{\sf ODRNM} + hq^{\sf GSM}_X \cdot \delta.
    $$
    First, note that for every $C \in \Gamma$,
    \begin{align}
        \label{eq:qhat1}\hat{q}^{\textsf{GSM}}_{C} &\geq 
        \hat{\mu}_C + \frac{\hat{\sigma}_C}{2}\left(\sqrt{\frac{b_C}{h}} - \sqrt{\frac{h}{b_C}}\right) \geq 
        \mu_C + \frac{\sigma_C}{2}\left(\sqrt{\frac{b_C}{h}} - \sqrt{\frac{h}{b_C}}\right),
    \end{align}
    and for every $i \in X$,
    \begin{align}
    \label{eq:qhat2}
      \hat{q}^{\sf GSM}_i 
      &\geq \hat{\mu}_i + \frac{\hat{\sigma}_i}{\sum_{i\in X} \hat{\sigma}_i} \frac{\hat{\sigma}_X}{2}\left(\sqrt{\frac{b}{h}} - \sqrt{\frac{h}{b}}\right) \geq 
      \mu_i + \frac{1}{1+\delta}\frac{\sigma_i}{\sum_{i\in X} \sigma_i} \frac{{\sigma}_X}{2}\left(\sqrt{\frac{b}{h}} - \sqrt{\frac{h}{b}}\right).
    \end{align}
    Similarly,
    \begin{align*}
        {q}^{\textsf{GSM}}_{C} & \geq 
        \mu_C + \frac{\sigma_C}{2}\left(\sqrt{\frac{b_C}{h}} - \sqrt{\frac{h}{b_C}}\right) 
        \geq \frac{1}{1+\delta}\left(\hat{\mu}_C + \frac{\hat{\sigma}_C}{2}\left(\sqrt{\frac{b_C}{h}} - \sqrt{\frac{h}{b_C}}\right)\right)
    \end{align*}
    and for every $i \in X$,
    \begin{align*}
      q^{\sf GSM}_i 
      & \geq 
      \mu_i + \frac{\sigma_i}{\sum_{i\in X} \sigma_i} \frac{{\sigma}_X}{2}\left(\sqrt{\frac{b}{h}} - \sqrt{\frac{h}{b}}\right) \geq 
      \frac{1}{1+\delta}\left(\hat{\mu}_i + \frac{\hat{\sigma}_i}{\sum_{i\in X} \hat{\sigma}_i} \frac{\hat{\sigma}_X}{2}\left(\sqrt{\frac{b}{h}} - \sqrt{\frac{h}{b}}\right)\right).
    \end{align*}
    In particular, $(1+\delta){\bq}^{\textsf{GSM}}$ is feasible for~\eqref{eq:gsm} under the estimates $\hat{\bmu}$ and $\hat{\bsigma}$. 

    Next, in our proof of Lemma~\ref{lem:UB} we show that,
    $$
    \sup_{\cD \in \cF(\bmu, \Sigma)} \mathbb{E}_{\bd \sim \cD}\; \left(
        C^H(\hat{\bq}^{\textsf{GSM}}, \bd)
        \right) \leq f(\hat{\bq}^{\textsf{GSM}}) + g(\hat{\bq}^{\textsf{GSM}}) + h(\hat{\bq}^{\textsf{GSM}}) + \sum_{r=1}^{R-1} s_r(\hat{\bq}^{\textsf{GSM}}),
    $$
    where $f,g,h,s_r$ are defined under the true parameters $\bmu$ and $\bsigma$. 
    
    We have, $f(\hat{\bq}^{\textsf{GSM}}) = f({\bq}^{\textsf{GSM}}) + h(\hat{q}^{\textsf{GSM}}_X - {q}^{\textsf{GSM}}_X)
    $ Hence, 
    \begin{align*}
        f(\hat{\bq}^{\sf GSM}) &\leq \sigma_X\left(\sqrt{bh} - h\sqrt{\frac{h}{b}}\right) + \sum_{C \in \Gamma} \frac{\sigma_C}{2}\left(\sqrt{b_Ch} - h\sqrt{\frac{h}{b_C}}\right) + h(\hat{q}_X^{\textsf{GSM}} - {q}_X^{\textsf{GSM}})\\
        & \leq \sigma_X\left(\sqrt{bh} - h\sqrt{\frac{h}{b}}\right) + \sum_{C \in \Gamma} \frac{\sigma_C}{2}\left(\sqrt{b_Ch} - h\sqrt{\frac{h}{b_C}}\right) + hq_X^{\sf GSM}\cdot \delta.\end{align*}
        The first inequality was shown in the proof of Lemma~\ref{lem:UB} and the second inequality follows from the feasibility of $(1+\delta)\bq^{\sf GSM}$ in \eqref{eq:gsm} under the estimates which implies that $h((1+\delta)q^{\sf GSM}_X - \hat{\mu}_X) \geq h(\hat{q}^{\sf GSM}_X - \hat{\mu}_X)$.
        Next, following the same lines of the proof of Lemma~\ref{lem:UB} and using the inequalities~\eqref{eq:qhat1} and~\eqref{eq:qhat2} we can upper-bound $g(\hat{\bq}^{\textsf{GSM}}) $, $h(\hat{\bq}^{\textsf{GSM}}) $, and $\sum_{r=1}^{R-1} s_r(\hat{\bq}^{\textsf{GSM}})$ as follows,
    \begin{align*}
    g(\hat{\bq}^{\textsf{GSM}}) 
    \leq  \frac{\sigma_X}{2} \left(\sqrt{bh} + h\sqrt{\frac{h}{b}}\right)\end{align*}
    and,
    \begin{align*}
    h(\hat{\bq}^{\textsf{GSM}}) 
    &= \frac{\diam(X)}{2n}  \sum_{i \in X} \left(\sqrt{\sigma_i^2 + \left(\hat{q}^{\textsf{GSM}}_i - \mu_i\right)^2} - \left(\hat{q}^{\textsf{GSM}}_i - \mu_i\right)\right)\\
    &\leq \frac{b+h}{2n}  \sum_{i \in X} \left(\sqrt{\sigma_i^2 + \left(\frac{1}{1+\delta}\frac{\sigma_i}{\sum_i \sigma_i}\frac{\sigma_X}{2}\left(\sqrt{\frac{b}{h}} - \sqrt{\frac{h}{b}}\right)\right)^2} - \left(\frac{1}{1+\delta}\frac{\sigma_i}{\sum_i \sigma_i}\frac{\sigma_X}{2}\left(\sqrt{\frac{b}{h}} - \sqrt{\frac{h}{b}}\right)\right)\right)\\
    &\leq (1+\delta)\frac{b+h}{2n}  \sum_{i \in X} \left(\sqrt{\sigma_i^2 + \left(\frac{\sigma_i}{\sum_i \sigma_i}\frac{\sigma_X}{2}\left(\sqrt{\frac{b}{h}} - \sqrt{\frac{h}{b}}\right)\right)^2} - \left(\frac{\sigma_i}{\sum_i \sigma_i}\frac{\sigma_X}{2}\left(\sqrt{\frac{b}{h}} - \sqrt{\frac{h}{b}}\right)\right)\right)\\
    &\leq (1+\delta)\frac{b+h}{2n} \left(\sqrt{n\sigma_X^2 + \left(\frac{\sigma_X}{2}\left(\sqrt{\frac{b}{h}} - \sqrt{\frac{h}{b}}\right)\right)^2} - \left(\frac{\sigma_X}{2}\left(\sqrt{\frac{b}{h}} - \sqrt{\frac{h}{b}}\right)\right)\right),
\end{align*}
where the second inequality follows from the fact that for every $a, x > 0$ and $\kappa \leq 1$,
\begin{align*}
    \sqrt{a+(\kappa x)^2} - \kappa x 
    = \frac{a}{\sqrt{a+(\kappa x)^2} + \kappa x }
    \leq \frac{a}{\sqrt{\kappa^2 a+(\kappa x)^2} + \kappa x }
    = \frac{1}{\kappa}\frac{a}{\sqrt{a+x^2} + x }
\end{align*}
and,
$$
\sum_{r=1}^{R-1} s_r(\hat{\bq}^{\textsf{GSM}}) \leq \sum_{C \in \Gamma} \frac{\sigma_C}{2}\left(\sqrt{b_C h} + h \sqrt{\frac{h}{b_C}}\right) + \sum_{C \in \clusters(\cH) \setminus \Gamma} \frac{\diam_p(C)}{2}  \sigma_C
$$
Implying that,
\begin{align*}
    \sup_{\cD \in \cF(\bmu, \Sigma)} \mathbb{E}_{\bd \sim \cD}\left(C^H(\hat{\bq}^{\textsf{GSM}}, \bd)\right)
    &\leq (1+\delta)\left(2\sigma_X \sqrt{bh} + \sum_{C \in \Gamma}\sigma_C \sqrt{b_Ch} + \sum_{C \in \clusters(H) \setminus \Gamma} \frac{\diam_p(C)}{2}\sigma_C\right) + hq_X^{\sf GSM}\cdot \delta
\end{align*}
The desired inequality follows from the above inequality and the lower-bound~\eqref{eq:LB}.    
\end{proof}

\section{The Generalized-Scarf-on-a-Metric (GSM) Solution: Omitted Proofs}

\subsection{Proof of Proposition~\ref{prop:bhl}}

    Given vectors $\bq \geq \bzero$ and $\bd \geq \bzero$, consider the linear program $\tilde{C}(\bq, \bd)$ we get from $C(\bq, \bd)$ by replacing $\ell_{ij}$ by $\min\{\ell_{ij}, b+h\}$. If $x^*$ is an optimal solution of $C(\bq, \bd)$, we must have $x^*_{ij}=0$ whenever $\ell_{ij} > b+h$ (if $\ell_{ij} > b+h$ and $x^*_{ij} >0$, then changing $x^*_{ij}$ to $0$ still gives a feasible solution with strictly smaller cost), and hence $x^*$ is feasible for $\tilde{C}(\bq, \bd)$ with same cost. Similarly, if we let $\tilde{x}^*$ be an optimal solution of $\tilde{C}(\bq, \bd)$, then we can assume without loss of generality that $\tilde{x}^*_{ij}=0$ whenever $\ell_{ij} > b+h$, i.e., whenever $\min\{\ell_{ij}, b+h\} = b+h$ (if $\ell_{ij} > b+h$, the objective coefficient of $\tilde{x}^*_{ij}$ in $\tilde{C}(\bq, \bd)$ is $-h-b+\min\{\ell_{ij}, b+h\} = 0$), which again implies that $\tilde{x}^*$ is feasible for $C(\bq, \bd)$ with same objective. Therefore $C(\bq, \bd)= \tilde{C}(\bq, \bd)$. Let $\tilde{\ell}_{ij} = \min\{\ell_{ij}, b+h\}$ for every $i,j \in X$, then $(X, \tilde{\ell})$ is still a metric space. Hence, by considering this new truncated metric space instead, we can assume without loss of generality that $\ell_{ij} \leq b+h$, as this does not change the objective value of a feasible solution $\bq \geq \bzero$ of~\eqref{eq:offline}.
    
    Now, under this assumption, we can assume without loss of generality that for every $\bq \geq \bzero$ and $\bd \geq \bzero$ the optimal solution of $C(\bq, \bd)$ fulfills all the demand or consumes all of the inventory, in fact, if there is a quantity of inventory $\delta q$ left at $i$ and a quantity of unfulfilled demand $\delta d$ at $j$, one can fulfill $\min\{\delta d, \delta q\}$ of the demand from the inventory at $i$, losing $\ell_{ij}\cdot\min\{\delta d, \delta q\} \leq (b+h)\cdot\min\{\delta d, \delta q\} $ in transportation cost, and gaining $(b+h)\cdot\min\{\delta d, \delta q\}$ in overage and underage cost, and hence get a solution of at least as good objective value.

\subsection{Proof of Lemma~\ref{lem:chebyshev}}

    Let $ \lambda \geq 1$,
    $
    \epsilon = \min\left\{\frac{1}{\sqrt{2}}, \frac{\lambda}{\sqrt{2}(1+\sqrt{1+\lambda^2})}\cdot \frac{1}{\nu}\right\},$
    $\eta = \nu \sqrt{2(1-\epsilon^2)}$, and finally let $m = |S|$. Consider the discrete distribution in Table~\ref{table:dist}
    over $3m+1$ vectors $\bv_0$, $\{\bv_i\}_{i \in S}$, $\{\bw_i\}_{i \in S}$ and $\{\bt_i\}_{i \in S}$ with respective probabilities $p_0$, $\{p_i\}_{i \in S}$, $\{q_i\}_{i \in S}$ and $\{r_i\}_{i \in S}$. We show that this is a probability distribution, of mean $\bzero$ and second moment $I_S$, the identity matrix of $\mathbb{R}^{S \times S}$.

    \begin{table}[h!]
    \centering
    \begin{tabular}{|c|c|}
    \hline
    \text{Probability} & \text{Vector} \\ 
    \hline
     $p_0=\frac{1}{2(1+\lambda^2)}$ & $\bv_0 = \frac{\lambda\epsilon\sqrt{2}}{\sigma_S}\bsigma_S$ \\[1mm]
    $p_i = \frac{\lambda^2}{2(1+\lambda^2)}\left(\frac{\sigma_i}{\sigma_S}\right)^2$ & $\bv_i =\frac{\epsilon\sqrt{2(1+\lambda^2)}}{\lambda}\frac{\sigma_S}{\sigma_i} \cdot {\evec_i} - \frac{\sqrt{2}(1+\sqrt{1+\lambda^2})}{\lambda}\frac{\epsilon}{\sigma_S}{ \bsigma_S}$\\[1mm]
    $q_i = \frac{1}{2m(1+m\eta^2)}$ & $\bw_i =m\eta\sqrt{2(1-\epsilon^2)}\evec_i$\\[1mm]
    $r_i = \frac{\eta^2}{2(1+m\eta^2)}$ & $\bt_i =-\frac{\sqrt{2(1-\epsilon^2)}}{\eta}\evec_i$\\[1mm]
    \hline
    \end{tabular}
    \\
    \vspace{3mm}
    \caption{Distribution taking $\bv_0$ w.p. $p_0$, $\bv_i, \bw_i, \bt_i$ w.p. $p_i, q_i, r_i$ respectively for every $i \in S$. $\evec_i$ is the indicator vector of $i \in S$.}\label{table:dist}
    \end{table}
    
    \vspace{3mm}{\noindent \bf Probability mass.} It is easy to check that,
    $p_0 + \sum_{i \in S}p_i + \sum_{i \in S} q_i + \sum_{i \in S} r_i = 1.$
    
    \vspace{3mm}{\noindent \bf Mean.} We have,
    \begin{align*}
        p_0 \bv_0 + \sum_{i \in S} p_i \bv_i 
        &=
        \frac{1}{\sqrt{2}(1+\lambda^2)}\frac{\lambda\epsilon}{\sigma_S}\bsigma_S + \sum_{i\in S} \frac{\lambda^2}{\sqrt{2}(1+\lambda^2)}\left(\frac{\sigma_i}{\sigma_S}\right)^2\frac{\epsilon\sqrt{(1+\lambda^2)}}{\lambda\sigma_i} 
         \sigma_S \cdot {\evec_i} \\
         &- \sum_{i \in S} \frac{\lambda^2}{\sqrt{2}(1+\lambda^2)}\left(\frac{\sigma_i}{\sigma_S}\right)^2\frac{1+\sqrt{(1+\lambda^2)}}{\lambda}\frac{\epsilon}{\sigma_S} \bsigma_S\\
         &= \frac{\lambda\epsilon}{\sqrt{2}(1+\lambda^2)}\frac{1}{\sigma_S}\bsigma_S + \frac{\lambda\epsilon}{\sqrt{2(1+\lambda^2)}}\frac{1}{\sigma_S}\bsigma_S - \frac{\lambda\epsilon}{\sqrt{2(1+\lambda^2)}}\frac{1}{\sigma_S}\left( 1+ \frac{1}{\sqrt{1+\lambda^2}} \right) \bsigma_S = \bzero
    \end{align*}
    and,
    \begin{align*}
         \sum_{i \in S} q_i \bw_i + \sum_{i \in S} r_i \bt_i
        &= 
        \sum_{i \in S} \frac{1}{\sqrt{2}m(1+m\eta^2)} m\eta\sqrt{(1-\epsilon^2)}\evec_i - \sum_{i \in S} \frac{\eta^2}{\sqrt{2}(1+m\eta^2)}\frac{\sqrt{(1-\epsilon^2)}}{\eta}\evec_i = \bzero.
    \end{align*}
    Hence, 
    $$
    p_0 \bv_0 + \sum_{i \in S} p_i \bv_i + \sum_{i \in S} q_i \bw_i + \sum_{i \in S} r_i \bt_i = \bzero.
    $$
    \vspace{3mm}{\noindent \bf Second Moment.} We have,
    \begin{align*}
        p_0 \bv_0\bv_0\trsp + \sum_{i \in S} p_i \bv_i\bv_i\trsp
        &=
        \frac{1}{(1+\lambda^2)}\frac{\lambda^2\epsilon^2}{\sigma_S^2}\bsigma_S\bsigma_S\trsp + \sum_{i\in S} \frac{\lambda^2}{(1+\lambda^2)}\left(\frac{\sigma_i}{\sigma_S}\right)^2\frac{\epsilon^2(1+\lambda^2)}{\lambda^2\sigma^2_i} 
         \sigma^2_S \cdot {\evec_i}{\evec_i}\trsp \\
         &- \sum_{i\in S} \frac{\lambda^2}{(1+\lambda^2)}\left(\frac{\sigma_i}{\sigma_S}\right)^2 \frac{\epsilon^2\sqrt{1+\lambda^2}(1+\sqrt{1+\lambda^2})}{\lambda^2\sigma_i}(\evec_i\bsigma_S\trsp+\bsigma_S\evec_i\trsp)\\
         &+ \sum_{i \in S} \frac{\lambda^2}{(1+\lambda^2)}\left(\frac{\sigma_i}{\sigma_S}\right)^2\frac{(1+\sqrt{(1+\lambda^2)})^2}{\lambda^2}\frac{\epsilon^2}{\sigma^2_S} \bsigma_S\bsigma_S\trsp\\\\
         &= \frac{1}{1+\lambda^2}\frac{\lambda^2\epsilon^2}{\sigma_S^2}\bsigma_S\bsigma_S^T + \sum_{i \in S} \epsilon^2{\evec_i\evec_i^T} -2\frac{\epsilon^2}{\sigma_S^2} \left( 1+ \frac{1}{\sqrt{1+\lambda^2}} \right)\bsigma_S\bsigma_S^T\\
         &+ \frac{\epsilon^2}{\sigma_S^2}\left( 1+ \frac{1}{\sqrt{1+\lambda^2}} \right)^2\bsigma_S\bsigma_S^T\\\\
         &= \sum_{i \in S} \epsilon^2 \evec_i\evec_i^T = \epsilon^2I_S
    \end{align*}
    and,
    \begin{align*}
         \sum_{i \in S} q_i \bw_i\bw_i\trsp + \sum_{i \in S} r_i \bt_i\bt_i\trsp
        &= 
        \sum_{i \in S} \frac{1}{m(1+m\eta^2)} m^2\eta^2(1-\epsilon^2)\evec_i\evec_i\trsp + \sum_{i \in S} \frac{\eta^2}{(1+m\eta^2)}\frac{(1-\epsilon^2)}{\eta^2}\evec_i\evec_i\trsp\\
        &= (1-\epsilon^2)I_S.
    \end{align*}
    Hence, 
    $$
    p_0 \bv_0\bv_0\trsp + \sum_{i \in S} p_i \bv_i\bv_i\trsp + \sum_{i \in S} q_i \bw_i\bw_i\trsp + \sum_{i \in S} r_i \bt_i\bt_i\trsp = I_S.
    $$
    \vspace{3mm}{\noindent \bf Conclusion.} Let $\bx$ be a random vector following the distribution in Table~\ref{table:dist} and consider the vector $\bd = \bmu_S + \Sigma^{\frac{1}{2}}_S\bx$, where $\Sigma^{\frac{1}{2}}_S$ is given by the Cholesky decomposition of $\Sigma_S$ i.e., such that $\Sigma^{\frac{1}{2}}_S(\Sigma^{\frac{1}{2}}_S)\trsp = \Sigma_S$ The mean of $\bd$ is $\bmu_S$ and its covariance is $\Sigma_S$. Moreover, 
    every vector in the support of $\bd$ is non-negative, that is because,
    $
    \mu_i \geq \frac{\sqrt{2}(1+\sqrt{1+\lambda^2})}{\lambda}\frac{\epsilon}{\sigma_S}\sigma_i^2
    $
    for every $i \in X$, by definition of $\epsilon$, and 
    $
    \mu_i \geq \frac{\sqrt{2(1-\epsilon^2)}}{\eta} \sigma_i
    $
    for every $i \in X$ by definition of $\eta$. 
    
    Finally, $\bd$ takes $\bmu_S + \Sigma^{\frac{1}{2}}_S\bv_0 = \bmu_S + \frac{\lambda\epsilon\sqrt{2}}{\sigma_S}\bsigma_S^2$ with probability $\Omega(\frac{1}{1+\lambda^2})$. Now for $\lambda \geq 1$ and since $\nu$ is constant, $\sqrt{2}\epsilon$ is larger than some universal constant $C \in (0,1)$. Hence, given $\alpha \geq 1$, taking $\lambda=\frac{\alpha}{C} \geq 1$ gives that $\bd \geq \bmu_S + \frac{\alpha}{\sigma_S} \bsigma_S^2$ w.p. $\Omega(\frac{1}{1+\lambda^2}) = \Omega(\frac{1}{1+\alpha^2})$. This proves the lemma.

\subsection{Proof of Lemma~\ref{lem:AtLeastMu}}

    Let us first show the following triangle inequality property of the fulfillment cost function $\Theta(\bq, \bd)$.
    \begin{claim}
    \label{clm:triangle}
        For every ${\bf u,v,w}$ such that $u_X=v_X=w_X$, we have,
        $$
        \Theta(\bu, \bw) \leq \Theta(\bu, \bv) + \Theta(\bv, \bw)
        $$
    \end{claim}
    {\noindent \em Proof.} Recall that we have,
    $$
    \Theta(\bu, \bv) = \inf_{\bx \geq \bzero} \left\{\sum_{i,k \in X} \ell_{ik} x_{ik} \;\left|\; \sum_{k \in X} x_{ik} = u_i,\; \sum_{i\in X} x_{ik} = v_k \right.\right\},
    $$
    $$
    \Theta(\bv, \bw) = \inf_{\bx \geq \bzero} \left\{\sum_{k,j \in X} \ell_{kj} x_{kj} \;\left|\; \sum_{j \in X} x_{kj} = v_k,\; \sum_{k\in X} x_{kj} = w_j\right.\right\},
    $$
    and,
    $$
    \Theta(\bu, \bw) = \inf_{\bx \geq \bzero} \left\{\sum_{i, j} \ell_{ij} x_{ij} \;\left|\; \sum_{j \in X} x_{ij} = u_i,\; \sum_{i\in X} x_{ij} = w_j\right.\right\}.
    $$
    Let $\bx^1$ and $\bx^2$ be optimal solutions of the first and second LPs respectively. Consider the solution 
    $$
    x_{ij} = \sum_{k: v_k \neq 0} \frac{x^1_{ik} x^2_{kj}}{v_k}
    $$
    We have,
    \begin{align*}
         \sum_j x_{ij} 
         &= \sum_j \sum_{k: v_k \neq 0} \frac{x^1_{ik} x^2_{kj}}{v_k}
         = \sum_{k: v_k \neq 0} \frac{x^1_{ik} \sum_j x^2_{kj}}{v_k} 
         = \sum_{k: v_k \neq 0} x^1_{ik}
         = \sum_{k} x^1_{ik}
         = u_i
    \end{align*}
    where the second to last equality follows from the fact that when $v_k= 0$ it must be the case that $x^1_{ik} = 0$ by feasibility of $\bx^1$.
    Similarly,
    \begin{align*}
         \sum_i x_{ij} 
         &= \sum_i \sum_{k: v_k \neq 0} \frac{x^1_{ik} x^2_{kj}}{v_k} = \sum_{k: v_k \neq 0} \frac{x^2_{kj} \sum_i x^1_{ik}}{v_k} = \sum_{k: v_k \neq 0} x^2_{kj} = \sum_{k} x^2_{kj} = w_j
    \end{align*}
    Hence, $\bx$ is feasible for the third LP. The objective value of this solution is,
    \begin{align*}
         \sum_{i,j} \ell_{ij} x_{ij} 
         &= \sum_{i,j} \sum_{k: v_k \neq 0} \ell_{ij} \frac{x^1_{ik} x^2_{kj}}{v_k}\\
         &\leq \sum_{i,j} \sum_{k: v_k \neq 0} (\ell_{ik} + \ell_{kj}) \frac{x^1_{ik} x^2_{kj}}{v_k}\\
         &= \sum_{i,j} \sum_{k: v_k \neq 0} \ell_{ik} \frac{x^1_{ik} x^2_{kj}}{v_k} + \sum_{i,j} \sum_{k: v_k \neq 0} \ell_{kj} \frac{x^1_{ik} x^2_{kj}}{v_k}\\
         &= \sum_{i} \sum_{k: v_k \neq 0} \ell_{ik} \frac{x^1_{ik} \sum_j x^2_{kj}}{v_k} + \sum_{j} \sum_{k: v_k \neq 0} \ell_{kj} \frac{x^2_{kj} \sum_i x^1_{ik}}{v_k}\\
         &= \sum_{i} \sum_{k: v_k \neq 0} \ell_{ik} x^1_{ik} + \sum_{j} \sum_{k: v_k \neq 0} \ell_{kj} x^2_{kj}\\
         &= \sum_{ik} \ell_{ik} x^1_{ik} + \sum_{kj} \ell_{kj} x^2_{kj} = \Theta(\bu, \bv) + \Theta(\bv, \bw),
    \end{align*}
    where the inequality follows from the triangle inequality $\ell_{ij} \leq \ell_{ik} + \ell_{kj}$. Which implies the claim.

    \hfill $\square$

    Next, fix an optimal solution $\bq^*$ of \eqref{eq:offline}. We define two inventory vectors $\hat{\bq}^*$ and $\tilde{\bq}^*$ as follows,
    \begin{itemize}
        \item If $q^*_X < \mu_X$, then,
        $$\hat{\bq}^* \in \argmin_{\bq \geq \bzero}\{\Theta(\bq, \bq^*) \;|\; q_X = q^*_X,\; q_i \leq \mu_i, \forall i\in X\} \text{ and } \tilde{\bq}^* = \bmu.$$
        \item Otherwise,
        $$\hat{\bq}^* \in \argmin_{\bq \geq \bzero}\{\Theta(\bq, \bq^*) \;|\; q_X = q^*_X,\; q_i \geq \mu_i, \forall i\in X\} \text{ and } \tilde{\bq}^* = \hat{\bq}^*.$$
    \end{itemize}
    In words, $\hat{\bq}^*$ is the ``closest" vector to $\bq^*$ of same total inventory that ``fills" the locations with less than mean inventory, and $\tilde{\bq}^*$ is the ``closest" vector to $\bq^*$ that ``fills" the locations with less than mean by potentially using more inventory than $\bq^*$ if necessary.

    We show the following claim that upper-bounds the cost $C(\tilde{\bq}^*, \bd)$ in function of the cost $C(\bq^*, \bd)$, a cost that accounts for the transport between $\bq^*$ and $\tilde{\bq}^*$, and a cost that accounts for any extra units that $\tilde{\bq}^*$ might contain compared to $\bq^*$. 
    \begin{claim}\label{clm:boundtriangle}
        For every $\bd \geq \bzero$,
        \begin{align*}
            C(\tilde{\bq}^*, \bd) 
            &\leq 3C(\bq^*, \bd) + 3b(\tilde{q}^*_X - q^*_X) + \Theta(\hat{\bq}^*, \bq^*)
        \end{align*}
    \end{claim}
    {\noindent \em Proof.}
    First, 
        \begin{align*}
            C(\tilde{\bq}^*, \bd)
            &= h(\tilde{q}^*_X - d_X)^+ + b(d_X- \tilde{q}^*_X)^+ + \Theta(\tilde{\bq}^*, \bd)\\
            &\leq h(\tilde{q}^*_X - \hat{q}^*_X)^+ + h(\hat{q}^*_X - d_X)^+ + b(d_X- \hat{q}^*_X)^+ + \Theta(\tilde{\bq}^*, \bd)\\
            &\leq h(\tilde{q}^*_X - \hat{q}^*_X)^+ + h(\hat{q}^*_X - d_X)^+ + b(d_X- \hat{q}^*_X)^+ + \Theta(\hat{\bq}^*, \bd) + (b+h)(\tilde{q}^*_X - \hat{q}^*_X)^+\\
            &=C(\hat{\bq}^*, \bd) + (b+2h)(\tilde{q}^*_X - \hat{q}^*_X)\\
            &\leq C(\hat{\bq}^*, \bd) + 3b(\tilde{q}^*_X - \hat{q}^*_X)
        \end{align*}
        where the second inequality follows from the fact that $\hat{\bq}^* \leq \tilde{\bq}^*$ and hence, one feasible fractional matching of $\tilde{\bq}^*$ to $\bd$ is to use the optimal matching of $\hat{\bq}^*$ to $\bd$ (the matching of cost $\Theta(\hat{\bq}^*, \bd)$) then match the rest of the inventory (at most $(\tilde{q}^*_X - \hat{q}^*_X)$) arbitrarily at a cost of at most $h+b$ per unit. 

        Hence, 
        \begin{align}\label{eq:eq1}
            C(\tilde{\bq}^*, \bd) 
            &\leq C(\hat{\bq}^*, \bd) + 3b(\tilde{q}^*_X - \hat{q}^*_X)
        \end{align}

        Second, given vectors $\bu$ and $\bv$ such that $u_X \leq v_X$ we denote by $\bv_{\bu}$ the sub-vector of $\bv$ (i.e. such that $\bv_{\bu} \leq \bv$) to which $\bu$ is matched in an optimal matching $\Theta(\bv, \bu)$ (if there are many optimal matching we fix an arbitrary one). Hence, $\Theta(\bv, \bu) = \Theta(\bv_{\bu}, \bu)$. We have,
        \begin{align*}
            C(\hat{\bq}^*, \bd)
            &= h(\hat{q}^*_X - d_X)^+ + b(d_X- \hat{q}^*_X)^+ + \Theta(\hat{\bq}^*, \bd) = h(q^*_X - d_X)^+ + b(d_X- q^*_X)^+ + \Theta(\hat{\bq}^*, \bd)
        \end{align*}

        We distinguish two cases. First, the case when $d_X \leq \hat{q}^*_X = q^*_X$. In this case, consider an optimal (perfect fractional) matching solution of $\Theta(\hat{\bq}^*, \bq^*)$. In this matching, the sub-vector $\bq^*_{\bd}$ of $\bq^*$ is matched to some sub-vector $\bq'$ of $\hat{\bq}^*$, and we have,
         \begin{align*}
            \Theta(\hat{\bq}^*, \bd)
            &= \Theta(\hat{\bq}^*_{\bd}, \bd) \leq \Theta(\bq', \bd) \leq \Theta(\bq', {\bq}^*_{\bd}) + \Theta({\bq}^*_{\bd}, \bd) \leq \Theta(\hat{\bq}^*, {\bq}^*) + \Theta({\bq}^*, \bd)
        \end{align*}
        where the first inequality follows from the definition of $\hat{\bq}^*_{\bd}$. In fact, by definition of $\hat{\bq}^*_{\bd}$, it must be the case that no sub-vector of $\hat{\bq}^*$ of total inventory equal to $d_X$ can be matched to $\bd$ with cost strictly less than $\Theta(\hat{\bq}^*_{\bd}, \bd)$. The inequality follows from the fact that $\bq'$ is a sub-vector of $\hat{\bq}^*$ of total inventory equal to $d_X$. The second inequality follows from the triangle inequality Claim~\ref{clm:triangle}, and the third from the fact that the sub-matching between $\bq'$ and ${\bq}^*_{\bd}$ of the optimal matching between $\hat{\bq}^*$ and ${\bq}^*$ considered in the definition of $\bq'$  is a feasible matching for $\Theta(\bq', {\bq}^*_{\bd})$. 
        
        Second, the case when 
        $d_X > \hat{q}^*_X = q^*_X$, in this case we have,
         \begin{align*}
            \Theta(\hat{\bq}^*, \bd)
            &= \Theta(\hat{\bq}^*, \bd_{\hat{\bq}^*}) \leq\Theta(\hat{\bq}^*, \bq^*) +\Theta(\bq^*, \bd_{\hat{\bq}^*})  \leq \Theta(\hat{\bq}^*, {\bq}^*) + \Theta({\bq}^*, \bd_{\bq^*}) + \Theta(\bd_{\bq^*}, \bd_{\hat{\bq}^*})\\
            &\leq \Theta(\hat{\bq}^*, {\bq}^*) + \Theta({\bq}^*, \bd) + (b+h)(d_X - q^*_X)
        \end{align*}
        where the second and third inequality follows from the triangle inequality Claim~\ref{clm:triangle},  and the last inequality from the fact that at most a quantity $\sum_{i \in X}(\bd_{\bq^*} - \bd_{\hat{\bq}^*})^+_i$ is moved in an optimal matching of $\Theta(\bd_{\bq^*}, \bd_{\hat{\bq}^*})$, and that the diameter of $(X, \ell)$ is at most $b+h$. In fact, because of the triangle inequality, one can assume without loss of generality that an optimal matching of $\Theta(\bd_{\bq^*}, \bd_{\hat{\bq}^*})$ matches  $\min\{(\bd_{\bq^*})_i, (\bd_{\hat{\bq}^*})_i\}$ of each coordinate $i$ of $\bd_{\bq^*}$ in place at a cost zero, otherwise, if there exists a coordinate $i$ such that an amount $\delta$ of $(\bd_{\bq^*})_{i}$ and of $(\bd_{\hat{\bq}^*})_{i}$ are matched outside of $i$, say to locations $j$ and $j'$, then we can match the amount $\delta$ from $(\bd_{\bq^*})_{i}$ and $(\bd_{\hat{\bq}^*})_{i}$ in-place and match the amount $\delta$ in $j$ and in $j'$ that we have unmatched together. Because $\ell_{jj'} \leq \ell_{ji}+\ell_{ij'}$, the new cost $\delta \cdot 0 + \delta \cdot \ell_{jj'}$ is smaller or equal to the old one $\delta \ell_{ji}+ \delta\ell_{ij'}$.
        We conclude by noticing that,
        \begin{align*}
            \sum_{i \in X}(\bd_{\bq^*} - \bd_{\hat{\bq}^*})^+_i
            &\leq \sum_{i \in X}(\bd - \bd_{\hat{\bq}^*})^+_i =\sum_{i \in X}(\bd - \bd_{\hat{\bq}^*})_i =(d_X - (\bd_{\hat{\bq}^*})_X)  = (d_X-q^*_X)
        \end{align*}
        
        Hence, in both cases, we have,
        $$
        \Theta(\hat{\bq}^*, \bd) \leq \Theta(\hat{\bq}^*, {\bq}^*) + \Theta({\bq}^*, \bd) + (b+h)(d_X - q^*_X)^+
        $$
        implying that,
        \begin{align*}
            C(\hat{\bq}^*, \bd)
            &= h(q^*_X - d_X)^+ + b(d_X- q^*_X)^+ + \Theta(\hat{\bq}^*, \bd)\\
            &\leq h(q^*_X - d_X)^+ + b(d_X- q^*_X)^+ + \Theta(\hat{\bq}^*, {\bq}^*) + \Theta({\bq}^*, \bd) + (b+h)(d_X - q^*_X)^+\\
            &\leq 3C(\bq^*, \bd) + \Theta(\hat{\bq}^*, \bq^*)
        \end{align*}
        as $h \leq b$.
        Hence,
        \begin{align}\label{eq:eq2}
            C(\hat{\bq}^*, \bd)
            \leq 3C(\bq^*, \bd) + \Theta(\hat{\bq}^*, \bq^*)
        \end{align}
        The claim follows from inequalities~\eqref{eq:eq1} and ~\eqref{eq:eq2}.
    
    \hfill $\square$

    Next, we show the following claim relating the term $3b(\tilde{q}^*_X - q^*_X)+ \Theta(\hat{\bq}^*, \bq^*)$ to the optimum of~\eqref{eq:offline}.
    \begin{claim}\label{clm:boundtozstar}
        $$3b(\tilde{q}^*_X - q^*_X)+ \Theta(\hat{\bq}^*, \bq^*) \leq O(1) \cdot z_{\textsf{ODRNM}}$$
    \end{claim}
    {\noindent \em Proof.}
    Let $\cD$ be the distribution we get from Lemma~\ref{lem:chebyshev} for $\alpha = 1$ and $S = X$. In particular, 
        $
        \mathbb{P}_{\bd \sim \cD} \left(\bd \geq \bmu\right) \geq \Omega\left(1\right).
        $
        We have,
        \begin{align*}
            z_{\textsf{ODRNM}} \geq \mathbb{E}_{\bd \sim \cD}(C(\bq^*, \bd))
            &= 
            \mathbb{E}_{\bd \sim \cD}\left(h(q^*_X - d_X)^+ + b(d_X - q^*_X)^+ + \Theta(\bq^*, \bd)\right)\\
            &\geq
            \Omega(1) \cdot \mathbb{E}_{\bd \sim \cD}\left(\left.h(q^*_X - d_X)^+ + b(d_X - q^*_X)^+ + \Theta(\bq^*, \bd) \;\right|\; \bd \geq \bmu\right)\\
            &\geq \Omega(1) \cdot \mathbb{E}_{\bd \sim \cD}\left(\left. b(d_X - q^*_X)^+ + \Theta(\bq^*, \bd) \;\right|\; \bd \geq \bmu\right).
        \end{align*}
        We distinguish two cases. The first case is when $q^*_X < \mu_X$. In this case, $\tilde{\bq}^* = \bmu$ and we have,
        \begin{align*}
            \mathbb{E}_{\bd \sim \cD}\left(\left. b(d_X - q^*_X)^+ + \Theta(\bq^*, \bd) \;\right|\; \bd \geq \bmu\right)
            &= \mathbb{E}_{\bd \sim \cD}\left(\left. b(d_X - q^*_X) + \Theta(\bq^*, \bd) \;\right|\; \bd \geq \bmu\right)\\
            &= \mathbb{E}_{\bd \sim \cD}\left(\left. b(\mu_X - q^*_X) + b(d_X - \mu_X) + \Theta(\bq^*, \bd) \;\right|\; \bd \geq \bmu\right)\\
            &= b(\tilde{q}_X - q^*_X) + \mathbb{E}_{\bd \sim \cD}\left(\left.  b(d_X - \mu_X) + \Theta(\bq^*, \bd) \;\right|\; \bd \geq \bmu\right)\\
            &\geq b(\tilde{q}_X - q^*_X) + \frac{1}{2}\Theta(\bq^*, \bmu)\\
            &= b(\tilde{q}_X - q^*_X) + \frac{1}{2}\Theta(\bq^*, \bmu_{\bq^*})\\
            &\geq b(\tilde{q}_X - q^*_X) + \frac{1}{2}\Theta(\bq^*, \hat{\bq}^*)
        \end{align*}
        To see why the first inequality holds, consider an optimal matching of $\Theta(\bq^*, \bd)$, using this matching, we can construct a feasible matching for $\Theta(\bq^*, \bmu)$ as follows: use the exact same matching used to match the first $\bmu$ units of $\bd$ then match the rest of $\bq^*$ (at most $(d_X - \mu_X)$ total units stay unmatched) using any arbitrary matching. The cost of such matching is at most $\Theta(\bq^*, \bd)$ plus $(b+h)(d_X - \mu_X) \leq 2b (d_X - \mu_X)$. The last inequality holds by definition of $\hat{\bq}^*$ as a minimizer.
        
        Hence, the claim holds in this case.

        The second case is when $q^*_X \geq \mu_X$. In this case, $\tilde{\bq}^* = \hat{\bq}^*$ and in particular $\tilde{q}^*_X = q^*_X$. Let $\bd \geq \bmu$. We distinguish two further sub-cases. The first is when $\mu_X \leq d_X \leq q^*_X$. In this sub-case, we have,
        \begin{align*}
            b(d_X - q^*_X)^+ +  \Theta(\bq^*, \bd)
            &=
            0 + \Theta(\bq^*_{\bd}, \bd) \geq
            \Theta(\bq^*, (\bq^* - \bq^*_{\bd}) + \bd) \geq
            \Theta(\bq^*, \hat{\bq}^*)
        \end{align*}
        the second equality holds by noticing that for every vectors $\bv, \bu$ and $\bw$ it holds that $\Theta(\bu+\bw, \bv+\bw) \leq \Theta(\bu, \bv)$ as one can consider the feasible matching of $\Theta(\bu+\bw, \bv+\bw)$ where $\bw$ units are matched in-place and the rest of the units, i.e., $\bu$ to $\bv$ are matched in the same way as in an optimal matching of $\Theta(\bu, \bv)$. The last inequality follows again by definition of $\hat{\bq}^*$ as a minimizer.

        The second sub-case is when $d_X > q^*_X \geq \mu_X$. Let $\bd'$ be a sub-vector of $\bd$ such that $\bd' \geq \bmu$ and $d'_X = q^*_X$. We have,
        \begin{align*}
            b(d_X - q^*_X)^+ +  \Theta(\bq^*, \bd)
            &\geq
            \frac{1}{2}\Theta(\bq^*, \bd')\geq
            \frac{1}{2}\Theta(\bq^*, \hat{\bq}^*),
        \end{align*}
        where the first inequality holds because one can construct a feasible matching $\tau$ of $\Theta(\bq^*, \bd')$ from an optimal matching $\tau'$ of $\Theta(\bq^*, \bd)$ by using the same matching in $\tau'$ for the first units $\bd'$ of $\bd$ then matching the rest of $\bq^*$ (at most $(d_X-d'_X) = (d_X - q^*_X)$ units remain) arbitrarily at a cost at most $b+h \leq 2b$ per unit. The second inequality follows again by definition of $\hat{\bq}^*$ as a minimizer.

        Hence, in both sub-cases, we have that,
        $$
        b(d_X - q^*_X)^+ +  \Theta(\bq^*, \bd) \geq \frac{1}{2}\Theta(\bq^*, \hat{\bq}^*)
        $$
        implying that,
        $$
        \mathbb{E}_{\bd \sim \cD}\left(\left. b(d_X - q^*_X)^+ + \Theta(\bq^*, \bd) \;\right|\; \bd \geq \bmu\right) \geq \frac{1}{2}\Theta(\bq^*, \hat{\bq}^*)
        $$
        and hence the claim holds in this case as well by noticing that,
        $
        b(\tilde{q}_X - q^*_X) = 0.
        $
    \hfill $\square$

    Now given any distribution $\cD \in \cF(\bmu, \Sigma)$, it follows from the Claims~\ref{clm:boundtriangle} and~\ref{clm:boundtozstar} that,
    \begin{align*}
        \mathbb{E}_{\bd \sim \cD}(C(\tilde{\bq}^*,\bd)) 
        &\leq
        3\mathbb{E}_{\bd \sim \cD}(C(\bq^*,\bd)) + O(1) \cdot z_{\textsf{ODRNM}}\\
        &\leq
        3z_{\textsf{ODRNM}} + O(1) \cdot z_{\textsf{ODRNM}}= O(1) \cdot z_{\textsf{ODRNM}}
    \end{align*}
    which proves the lemma as $\tilde{\bq}^* \geq \bmu$.  

\subsection{Proof of Lemma~\ref{lem:wshp}}

    The constructions are as follows,
    
    \vspace{3mm}{\noindent \bf Uniform:} Let $(X, \ell)$ be a uniform metric space such that $\ell_{ij} = \lambda$ for every $i,j \in X$. Let $\alpha > 1$, $\beta = 1$, $\gamma > \alpha$. Note that $\Delta_1 = \lambda/\alpha$ and $R=2$. The first partition $\cP_1$ consists of a single family containing $|X|$ clusters, each of them is a singleton point. Note that $\cP_1$ is an $(\alpha, \beta)$-well-separated partition of margin $\Delta_1$. In fact, there is $1$ ($\leq \beta$) family, each cluster has diameter $0$ ($< \alpha \Delta_1$), and each pair of clusters are separated by $\lambda$ ($> \Delta_1$). Then, the second partition $\cP_2$ consists of a single family given by a single cluster $X$. Note that $\cP_2$ is an $(\alpha, \beta)$-well-separated partition of margin $\Delta_2=\gamma \lambda/\alpha > \lambda$. In fact, there is $1$ ($\leq \beta$) family, each cluster has diameter $\lambda$ ($< \alpha \Delta_2$), and there is only one cluster so the separation condition is also verified.
    Hence $\cH = \cP_1, \cP_2$ is an $(\alpha,\beta, \gamma)$ well separated hierarchical partition of $(X, \ell)$.

    \vspace{3mm}{\noindent \bf Euclidean:} Let now $(X, \ell)$ be a subspace of the Euclidean space $\mathbb{R}^d$. Let $\alpha > 2\sqrt{d}$, $\beta=2^d$ and $\gamma \geq 2$ an even integer. We take the last partition $\cP_R$ to be a single family given by a single cluster $X$. Note that this is indeed an $(\alpha, \beta)$-well separated partitions of margin $\Delta_R$. In fact, there is $1$ ($\leq \beta$) family, each cluster has diameter at most $\diam(X)$ ($< \alpha\Delta_{R}$; note that $\Delta_R \geq \diam(X)$ by definition of $R$), and there is only one cluster so the separation condition is also verified.
    
    Consider a hypercube of side length $\Delta_R$ that contains all the points of $X$ (which is possible as $\Delta_R \geq \diam(X)$).
    For each of the subsequent partitions $\cP_r$ for $r=R-1, \dots, 1$, we construct $\cP_r$ as follows: we divide each side of the hypercube into $\gamma^{R-r}/2$ equal segments and label these segments as $1, 2, \dots, \gamma^{R-r}/2$ starting from some fixed corner of the hypercube. This partitions the hypercube into a grid of $\gamma^{(R-r)d}/2^d$ equal sub-hypercubes and induces labels of these sub-hypercubes that are vectors $(i_1, \dots, i_d)$ where $i_k \in \{1,\dots, \gamma^{R-r}/2\}$ for every $k$ denotes that the sub-hypercube corresponds to the segment $i_k$ of dimension $k \in [d]$.
    Now, label each dimension $k$ with a label $b_k \in \{0,1\}$. The family of sub-hypercubes whose coordinates $(i_1, \dots, i_d)$ are such that $i_k \equiv b_k [2]$ (i.e., $i_k$ even if $b_k=0$ and odd otherwise)
    are such that: (i) no two sub-hypercubes are adjacent (the labeling of every two sub-hypercubes have the same parity along every dimension, and hence, if the sub-hypercubes are different, there must exist a dimension $k$ such that the difference of their coordinates is at least $2$, which implies that they cannot be adjacent). (ii) the set of all families of sub-hypercubes constructed for different labeling of the dimensions gives the set of all sub-hypercubes (every sub-hypercube $(i_1, \dots, i_d)$ belongs to some family; specifically, the family where $b_k=0$ if $i_k$ is even and $b_k = 1$ otherwise, for every $k$).
    For each such family of sub-hypercubes we consider a family of clusters of $\cP_r$ consisting of the points inside the sub-hypercubes.
    
    Note that by construction, this leads to a coarsing set of partitions (the sub-hypercubes we use to construct each partition are subsets of the sub-hypercubes used to construct the partition at one level higher). We only need to show that for every $r=1, \dots, R-1$, the partition $\cP_r$ is indeed an $(\alpha, \beta)$-well-separated partition of margin $\Delta_r$. In fact, there is at most $2^d$ ($\leq \beta$) families in $\cP_r$ (the number of possible labelings of the dimensions), each cluster has diameter at most $\sqrt{d} \times 2\Delta_R / \gamma^{R-r}$  ($< \frac{\alpha}{\gamma^{R-r}}\Delta_R = \alpha \Delta_{r}$), and since the sub-hypercubes used to construct the clusters of each family are at least $2\Delta_R/\gamma^{R-r}$ far apart from each other (no adjacent sub-hypercubes in the same family), every two clusters of the same family are separated by at least $2\Delta_R/\gamma^{R-r}$ ($> \Delta_r$).
    
    \vspace{3mm}{\noindent \bf General:} Let $\alpha = 6\log(n)+1$, $\beta = \log(n)$, $\gamma > 12\log(n)+2$. Define $R$ and $\{\Delta_r\}_{r \in \{1, \dots, R\}}$ as in Definition~\ref{def:wshp}. We construct an an $(\alpha, \beta, \gamma)$-well-separated hierarchical partition of $(X, \ell)$ starting from the finest well-separated partition to the coarsest one.

    We begin by constructing the initial well-separated partition $\cP_1= \cF_1, \dots, \cF_{K_1}$. Given $\lambda>0$, let 
    $$
    B(x, \lambda) = \{y \in X \;|\; l_{x,y} \leq \lambda\}
    $$
    denote the ball of radius $\lambda$ around $x$. The construction runs in phases, in each phase $k$, the family $\cF_k$ is constructed. In the first phase, an arbitrary point $x_1 \in X$ is selected. Let $s_1 \geq 0$ be the smallest integer such that $$|B(x_1, (s_1+1)\Delta_{1})| < 2|B(x_1, s_1\Delta_{1})|.$$ Note that $s_1 \leq \log(n)$ since by construction,
    $$
    n \geq |B(x_1, s_1\Delta_{1})| \geq 2|B(x_1, (s_1-1)\Delta_{1})| \geq \cdots \geq 2^{s_1}|B(x_1, 0)| \geq 2^{s_1}
    $$
    Now set the first cluster of the family $\cF_1$ to $C_1 = B(x_1, s_1 \Delta_{1})$ and remove all the points $B(x_1, s_1 \Delta_{1})$ from $X$, then mark all the points in the boundary $B(x_1, (s_1+1) \Delta_{1}) \setminus B(x_1, s_1 \Delta_{1})$ and remove them from $X$ as well. Next, select a second point $x_2 \in X$, find again the smallest integer $s_2 \geq 0$ such that $$|B(x_2, (s_2+1)\Delta_{1})| < 2|B(x_2, s_2\Delta_{1})|.$$ Note that $s_2 \leq \log(n)$. Set the second cluster of the family to $C_2 = B(x_2, s_2 \Delta_{1})$ and remove
    all the points $B(x_2, s_2 \Delta_{1})$ from $X$, then mark all the points in the boundary $B(x_2, (s_2+1) \Delta_{1}) \setminus B(x_2, s_2 \Delta_{1})$ and remove them from $X$ as well. We continue until $X$ becomes empty, in which case, the phase ends and a new phase begins where a new family of clusters is constructed using the exact same construction applied to the set of marked points during the current phase. The construction ends when all points are clustered.

    We have,
    \begin{itemize}
        \item $\clusters(\cP_1)$ form a partition of $X$ by construction. 
        \item Every cluster constructed has diameter at most
        $
        2 \times \log(n) \Delta_1 < \alpha \Delta_1
        $
        (two times radius of the ball corresponding to the cluster).
        \item Within each phase, after a cluster is constructed, all the points within $\leq \Delta_1$ of the cluster are removed from $X$, hence, the clusters constructed next must have points that are $>\Delta_1$ away from the points of the current cluster. Hence, the distance between any two constructed clusters in the same family (phase) is
        $
        >\Delta_1
        $.
        \item Finally, note that, the number of marked points at the end of each phase is at most $1/2$ of the total number of points at the beginning of the phase, that is because when a cluster is constructed, we mark the points in the boundary and remove the points inside the cluster without marking. Since the number of points inside the cluster is at least $1/2$ of the total number of removed points, the number of marked points is therefore at most $1/2$ of the total number of points removed when the cluster is constructed. This implies a geometric decrease in the number of points considered at the beginning of each phase and implies a total number of phases (hence families) of at most $\beta=\log(n)$.
    \end{itemize}
    Hence, $\cP_1$ is indeed an $(\alpha, \beta)$-well-separated partition of $(X, \ell)$ of margin $\Delta_1$.

    We now construct the following partitions recursively in a similar manner. In particular, suppose we constructed the first $(\alpha, \beta)$-well-separated partitions $\cP_1, \dots, \cP_{r}$ of margins $\Delta_1, \dots, \Delta_r$ respectively for some $r \in \{1, \dots, R-1\}$, such that each $\cP_{l+1}$ is a coarsening of $\cP_{l}$ for every $l \in \{1, \dots, r-1\}$. We show how to construct the next $(\alpha, \beta)$-well-separated partition $\cP_{r+1} = \cF_1, \dots, \cF_{K_{r+1}}$ of margin $\Delta_{r+1}$ that is a coarsening of $\cP_r$. For each cluster $C \in \clusters(\cP_r)$, let $x_C \in C$ be some arbitrary point in $C$ that we refer to as the representative of $C$. Let $X^r = \{x_C \;|\; C \in \clusters(\cP_r)\}$ be the set of all representatives, and let $x \rightarrow C(x)$ denote the inverse function, i.e., such that $x_{C(x)} = x$ for every $x \in X^r$. Given $\lambda>0$ and $x \in X^r$, let $$B_r(x, \lambda) = \bigcup_{y \in X^r,\; l_{x,y} \leq \lambda} C(y),$$ denote the union of all clusters of $\cP_{r}$ whose representatives are at a distance at most $\lambda$ from $x$. The construction runs in phases, in each phase $k$, the family of clusters $\cF_k$ is constructed. In the first phase, an arbitrary representative $x_1 \in X^r$ is selected. Let $s_1 \geq 0$ be the smallest integer such that $$|B_r(x_1, 3(s_1+1)\Delta_{r+1})| < 2|B_r(x_1, 3s_1\Delta_{r+1})|.$$ Note that $s_1 \leq \log(n)$. Now set the first cluster of the family $\cF_1$ to $C_1 = B_r(x_1, 3s_1 \Delta_{r+1})$ and remove all the representatives $\{y \in X^r \;|\; l_{x_1,y} \leq 3s_1\Delta_{r+1}\}$ of the points $C_1$ from $X^r$, then mark all the representatives in the boundary $\{y \in X^r \;|\; 3s_1\Delta_{r+1} < l_{x_1,y} \leq 3(s_1+1)\Delta_{r+1}\}$ and remove them from $X^r$ as well. Next, select a second representative $x_2 \in X^r$, find again the smallest integer $s_2 \geq 0$ such that $$|B(x_2, 3(s_2+1)\Delta_{r+1})| < 2|B(x_2, 3s_2\Delta_{r+1})|.$$ Note that $s_2 \leq \log(n)$. Set the second cluster of the family to $C_2 = B_r(x_2, 3s_2 \Delta_{r+1})$ and remove the representatives $\{y \in X^r \;|\; l_{x_2,y} \leq 3s_2\Delta_{r+1}\}$ of the points $C_2$ from $X^r$, then mark all the representatives in the boundary $\{y \in X^r \;|\; 3s_2\Delta_{r+1} < l_{x_2,y} \leq 3(s_2+1)\Delta_{r+1}\}$ and remove them from $X^r$ as well. We continue until $X^r$ becomes empty, in which case, the phase ends and new phase begins where a new family of clusters is constructed using the exact same construction applied to the set of marked representative during the current phase. The construction ends when all points are clustered.

    The fact that the clusters of $\cP_{r+1}$ form a partition of $X$ and that $\cP_{r+1}$ coarsens $\cP_{r}$ is trivial from the construction, we only need to show that $\cP_{r+1}$ is an $(\alpha, \beta)$-well-separated partition of margin $\Delta_{r+1}$. During the construction of $\cP_{r+1}$, we have,
    \begin{itemize}
        \item Every cluster has diameter at most $2 \times 3\log(n)\Delta_{r+1}$ (which is the maximum distance between any two representatives in the cluster), plus at most $2\alpha \Delta_{r}$ (which is, by induction hypothesis, two times the maximum diameter of a cluster in $\cP_{r}$). In fact, the maximum distance between two points in the cluster is at most the distance of the two points to their representatives plus the distance between the representatives. This implies a total diameter of at most $$6\log(n)\Delta_{r+1} + 2\alpha \Delta_r = 6\log(n)\Delta_{r+1} + \frac{2\alpha}{\gamma} \Delta_{r+1} < \alpha \Delta_{r+1}$$
        \item In each phase, after a cluster is constructed, all the representatives within $3\Delta_{r+1}$ of the cluster are removed from $X^r$, hence, the clusters constructed next must have representatives that are at least $3\Delta_{r+1}$ away from the representatives in the current cluster. By subtracting twice the diameter of a cluster in $\cP_r$ (twice the maximum distance that a representative in a cluster of $\cP_r$ has to a point of the cluster), the points of the clusters constructed in each family of $\cP_{r+1}$ are at least a distance,
        \begin{align*}
            &3\Delta_{r+1}- 2 \alpha \Delta_r = \Delta_{r+1}(3 - 2\frac{\alpha}{\gamma}) > \Delta_{r+1}
        \end{align*}
        away from each other.
        \item Finally, note that, the number of marked representatives at the end of each phase is at most $1/2$ of the total number of representatives at the beginning of the phase, that is because when a cluster is constructed, we mark the representatives in the boundary and remove the representatives inside the cluster without marking. Since the number of representatives inside the cluster is at least $1/2$ of the total number of removed representatives, the number of marked representatives is therefore at most $1/2$ of the total number of representatives removed when the cluster is constructed. This implies a geometric decrease in the number of representatives considered at the beginning of each phase and implies a total number of phases of at most $\beta=\log(n)$.
    \end{itemize}

\subsection{Proof of Lemma~\ref{lem:UBHierarchicalMatching}}

    Consider a hierarchical fulfillment of the demand $\bd$ from the inventory $\bq$.

    In the fulfillment step inside the clusters of $\cP_1$, we distinguish two cases. The first case is when $\Delta_1 = \frac{1}{\alpha} \min_{i,j \in X: i \neq j} \ell_{ij}$. In this case, the clusters $\clusters(\cP_1)$ are given by singleton locations, and we are just fulfilling a maximum demand inside each location using inventory inside the location with no cost. The second case is when $\Delta_1 = \frac{1}{n\alpha} \sup_{i,j \in X: i \neq j} \ell_{ij}$, in this case, the maximum diameter of a cluster in $\clusters(\cP_1)$ is at most $\frac{\diam(X)}{n}$. Because we first fulfill demand locally from each location before moving to the fulfillment inside the clusters of $\cP_1$, in the worst case, inside each cluster $C \in \clusters(\cP_1)$, the remaining demand
    $
    \sum_{i \in C}(d_i - q_i)^+
    $
    will all be fulfilled inside the cluster at a cost of at most $
    \diam(C) \cdot \sum_{i \in C}(d_i - q_i)^+
    $
    hence a total cost of at most
    $$
    \sum_{C \in \clusters(\cP_1)}\diam(C) \cdot \sum_{i \in C}(d_i-q_i)^+ \leq \frac{\diam(X)}{n}\sum_{i \in X} (d_i - q_i)^+
    $$
    
    Next, in the fulfillment step inside the clusters of $\cP_2$ there is at most 
    $
    (d_C - q_C)^+
    $
    quantity of demand that is not fulfilled in the first step inside each cluster $C \in \clusters(\cP_1)$. In the worst case, all of this demand is fulfilled inside the parent cluster of $C$ in $\cP_2$. Overall, this leads to a total cost of at most
    $
    \sum_{C \in \clusters(\cP_{1})} \diam_p(C) \cdot (d_C - q_C)^+.
    $

    Next, in the fulfillment step inside the clusters of $\cP_3$ there is at most 
    $
    (d_C - q_C)^+
    $
    quantity of demand that is not fulfilled in the second step inside each cluster $C \in \clusters(\cP_2)$. In the worst case, all of this demand is fulfilled inside the parent cluster of $C$ in $\cP_3$. Overall, this leads to a total cost of at most
    $
    \sum_{C \in \clusters(\cP_{2})} \diam_p(C) \cdot (d_C - q_C)^+,
    $
    and so on.

    After the last step $R$ (recall that the last partition $\cP_R$ contains a single family of a single cluster $X$). An underage cost of $b$ per-unit is paid for the remaining unfulfilled demand of total
    $
    (d_X-q_X)^+,
    $
    and an overage cost of $h$ per-unit is paid for the remaining inventory of total,
    $
    (q_X-d_X)^+.
    $
    
    The hierarchical fulfillment gives a feasible solution of $C(\bq, \bd)$. By the above reasoning, the objective value of such solution is at most,
    $$
    C^H(\bq, \bd) = h(q_X-d_X)^+ + b(d_X-q_X)^+ + \sum_{r=1}^{R-1} \sum_{C \in \clusters(\cP_{r})} \diam_p(C) \cdot (d_C - q_C)^+ + \frac{\diam(X)}{n}\sum_{i \in X}(d_i-q_i)^+
    $$
    Hence, $C(\bq, \bd) \leq C^H(\bq, \bd)$.

\subsection{Proof of Lemma~\ref{lem:UB}}

    For every $\bq \in \cQ$ we have,
    \begin{align*}
        \sup_{\cD \in \cF(\bmu, \Sigma)} \mathbb{E}_{\bd \sim \cD}\left(C^H(\bq, \bd)\right)
        &= \sup_{\cD \in \cF(\bmu, \Sigma)} \mathbb{E}_{\bd \sim \cD}\left( h(q_X-d_X)^+ + b(d_X-q_X)^+ + \frac{\diam(X)}{n}\sum_{i \in X}(d_i-q_i)^+\right.\\
        &\left.\quad\quad  \quad \quad \quad \quad \quad \quad +\sum_{r=1}^{R-1} \sum_{C \in \clusters(\cP_{r})} \diam_p(C) \cdot (d_C - q_C)^+\right)\\
        &= h(q_X-\mu_X) + \sup_{\cD \in \cF(\bmu, \Sigma)} \mathbb{E}_{\bd \sim \cD}\left((b+h)(d_X-q_X)^+ + \frac{\diam(X)}{n}\sum_{i \in X}(d_i-q_i)^+\right.\\
        &\left.\quad \quad \quad \quad \quad \quad \quad \quad \quad \quad \quad \quad \quad +\sum_{r=1}^{R-1} \sum_{C \in \clusters(\cP_{r})} \diam_p(C) \cdot (d_C - q_C)^+\right),
    \end{align*}
    where the last equality follows from the fact that $x^+=x+(-x)^+$ for every $x \in \mathbb{R}$.
    Then by observing that the supremum of the sum of functions is smaller than the sum of supremums of the functions, we have,
    \begin{align*}
        &\sup_{\cD \in \cF(\bmu, \Sigma)} \mathbb{E}_{\bd \sim \cD}\left(C^H(\bq, \bd)\right) \\
        &\leq h(q_X-\mu_X) + (b+h)\sup_{\cD \in \cF(\bmu, \Sigma)} \mathbb{E}_{\bd \sim \cD}\left((d_X-q_X)^+\right) + \frac{\diam(X)}{n}\sum_{i \in X} \sup_{\cD \in \cF(\bmu, \Sigma)}\mathbb{E}_{\bd \sim \cD} (d_i-q_i)^+\\
        &\quad + \sum_{r=1}^{R-1} \sum_{C \in \clusters(\cP_{r})} \diam_p(C) \cdot  \sup_{\cD \in \cF(\bmu, \Sigma)}\mathbb{E}\left((d_C - q_C)^+\right)\\
        &\leq h(q_X-\mu_X) + \frac{b+h}{2} \left(\sqrt{\sigma_X^2 + \left(q_X - \mu_X\right)^2} - \left(q_X - \mu_X\right)\right)
        \\
        &\quad + \frac{\diam(X)}{2n}  \sum_{i \in X} \left(\sqrt{\sigma_i^2 + \left(q_i - \mu_i\right)^2} - \left(q_i - \mu_i\right)\right)  + \sum_{r=1}^{R-1} \sum_{C \in \clusters(\cP_{r})} \frac{\diam_p(C)}{2} \cdot \left(\sqrt{\sigma_C^2 + \left(q_C - \mu_C\right)^2} - \left(q_C - \mu_C\right)\right)
    \end{align*}
in the last inequality we use the fact that
$
\mathbb{E}(X^+) = \frac{1}{2}\left(\mathbb{E}(|X|) + \mathbb{E}(X)\right) \leq \frac{1}{2}\left(\sqrt{\mathbb{E}(X^2)} +\mathbb{E}(X)\right)
$
for every random variable $X$.

Let 
$$
f(\bq) = h(q_X-\mu_X),
$$
and
$$
g(\bq) = \frac{b+h}{2} \left(\sqrt{\sigma_X^2 + \left(q_X - \mu_X\right)^2} - \left(q_X - \mu_X\right)\right),
$$
and
$$
h(\bq) = \frac{\diam(X)}{2n}  \sum_{i \in X} \left(\sqrt{\sigma_i^2 + \left(q_i - \mu_i\right)^2} - \left(q_i - \mu_i\right)\right),
$$
and for every $r=1, \dots, R-1$, 
$$
s_r(\bq) = \sum_{C \in \clusters(\cP_{r})} \frac{\diam_p(C)}{2}\cdot \left(\sqrt{\sigma_C^2 + \left(q_C - \mu_C\right)^2} - \left(q_C - \mu_C\right)\right)
$$
we now bound each one of these functions for $\bq=\bq^{\textsf{GSM}}$.

Let
$
\bq'
$ such that $$
q'_i = \frac{\sigma_i}{\sum_i \sigma_i}\sigma_X\left(\sqrt{\frac{b}{h}} - \sqrt{\frac{h}{b}}\right),
$$
note that in particular, $q'_X = \sigma_X\left(\sqrt{\frac{b}{h}} - \sqrt{\frac{h}{b}}\right)$. Then, for $r=1, \dots, R-1$, let
$
\bq^r
$
some vector such that for every $C \in \clusters(\cP_r)$,
\begin{align*}
    q^r_{C} = \left\{
    \begin{matrix}
        \frac{\sigma_C}{2}\left(\sqrt{\frac{b_C}{h}} - \sqrt{\frac{h}{b_C}}\right) & \text{ if } C \in \Gamma\\
        0 & \text{ otherwise }
    \end{matrix}\right.
\end{align*}
The vector
$
\bq = \bmu + \bq' + \sum_{r=1}^{R-1} \bq^r
$
is therefore a feasible solution for~\eqref{eq:gsm}. Hence,
\begin{align*}
    f(\bq^{\textsf{GSM}})\leq h q'_X + h\sum_{r=1}^{R-1} q^r_X
    = h q'_X + \sum_{r=1}^{R-1} h \sum_{C \in \clusters(\cP_r)} q^r_C = \sigma_X\left(\sqrt{bh} - h\sqrt{\frac{h}{b}}\right) + \sum_{C \in \Gamma} \frac{\sigma_C}{2}\left(\sqrt{b_Ch} - h\sqrt{\frac{h}{b_C}}\right) 
\end{align*}
Next, we have,
\begin{align*}
    g(\bq^{\textsf{GSM}}) 
    &= \frac{b+h}{2} \left(\sqrt{\sigma_X^2 + \left(q^{\textsf{GSM}}_X - \mu_X\right)^2} - \left(q^{\textsf{GSM}}_X - \mu_X\right)\right)\\
    &\leq \frac{b+h}{2} \left(\sqrt{\sigma_X^2 + \left( \frac{\sigma_X}{2}\left(\sqrt{\frac{b}{h}} - \sqrt{\frac{h}{b}}\right) \right)^2} - \frac{\sigma_X}{2}\left(\sqrt{\frac{b}{h}} - \sqrt{\frac{h}{b}}\right)\right) = \frac{\sigma_X}{2} \left(\sqrt{bh} + h\sqrt{\frac{h}{b}}\right)
\end{align*}
where the inequality holds by the fact that the function $x \rightarrow \sqrt{a^2 + x^2} - x$ is non-increasing for every $a \in \mathbb{R}$, and that $q^{\textsf{GSM}}_X - \mu_X \geq \sigma_X(\sqrt{\frac{b}{h}} - \sqrt{\frac{h}{b}}) \geq \frac{\sigma_X}{2}(\sqrt{\frac{b}{h}} - \sqrt{\frac{h}{b}})$. Then,
\begin{align*}
    h(\bq^{\textsf{GSM}}) 
    &= \frac{\diam(X)}{2n}  \sum_{i \in X} \left(\sqrt{\sigma_i^2 + \left(q^{\textsf{GSM}}_i - \mu_i\right)^2} - \left(q^{\textsf{GSM}}_i - \mu_i\right)\right)\\
    &\leq \frac{b+h}{2n}  \sum_{i \in X} \left(\sqrt{\sigma_i^2 + \left(\frac{\sigma_i}{\sum_i \sigma_i}\frac{\sigma_X}{2}\left(\sqrt{\frac{b}{h}} - \sqrt{\frac{h}{b}}\right)\right)^2} - \left(\frac{\sigma_i}{\sum_i \sigma_i}\frac{\sigma_X}{2}\left(\sqrt{\frac{b}{h}} - \sqrt{\frac{h}{b}}\right)\right)\right)\\
    &= \frac{b+h}{2n}  \sum_{i \in X} \frac{\sigma_i}{\sum_i \sigma_i}\left(\sqrt{\left(\sum_i \sigma_i\right)^2 + \left(\frac{\sigma_X}{2}\left(\sqrt{\frac{b}{h}} - \sqrt{\frac{h}{b}}\right)\right)^2} - \left(\frac{\sigma_X}{2}\left(\sqrt{\frac{b}{h}} - \sqrt{\frac{h}{b}}\right)\right)\right)\\
    &\leq \frac{b+h}{2n} \left(\sqrt{n\sigma_X^2 + \left(\frac{\sigma_X}{2}\left(\sqrt{\frac{b}{h}} - \sqrt{\frac{h}{b}}\right)\right)^2} - \left(\frac{\sigma_X}{2}\left(\sqrt{\frac{b}{h}} - \sqrt{\frac{h}{b}}\right)\right)\right),
\end{align*}
where the first inequality holds by $\diam(X) \leq b+h$ and by the fact that the function $x \rightarrow \sqrt{a^2 + x^2} - x$ is non-increasing for every $a \in \mathbb{R}$. The last inequality holds by Cauchy Schwartz inequality implying that 
$
\sum_i \sigma_i \leq \sqrt{n}\sigma_X.
$
Now notice that for every $a, b >0$ we have,
\begin{align*}
    \sqrt{na^2 + b^2} - b 
    &= \frac{na^2}{\sqrt{na^2+b^2}+b}\leq \frac{na^2}{\sqrt{a^2+b^2}+b}= n(\sqrt{a^2 + b^2} - b)
\end{align*}
hence,
\begin{align*}
    h(\bq^{\textsf{GSM}}) 
    &\leq \frac{b+h}{2n} \left(\sqrt{n\sigma_X^2 + \left(\frac{\sigma_X}{2}\left(\sqrt{\frac{b}{h}} - \sqrt{\frac{h}{b}}\right)\right)^2} - \left(\frac{\sigma_X}{2}\left(\sqrt{\frac{b}{h}} - \sqrt{\frac{h}{b}}\right)\right)\right)\\
    &\leq \frac{b+h}{2} \left(\sqrt{\sigma_X^2 + \left(\frac{\sigma_X}{2}\left(\sqrt{\frac{b}{h}} - \sqrt{\frac{h}{b}}\right)\right)^2} - \left(\frac{\sigma_X}{2}\left(\sqrt{\frac{b}{h}} - \sqrt{\frac{h}{b}}\right)\right)\right)\\
    &\leq \frac{\sigma_X}{2} \left(\sqrt{bh} + h\sqrt{\frac{h}{b}}\right)
\end{align*}
Finally, for $s_r$ we have,
\begin{align*}
    &s_r(\bq^{\textsf{GSM}}) \\
    &=\sum_{C \in \clusters(\cP_{r})} \frac{\diam_p(C)}{2}  \left(\sqrt{\sigma_C^2 + \left(q^{\textsf{GSM}}_C - \mu_C\right)^2} - \left(q^{\textsf{GSM}}_C - \mu_C\right)\right)\\
    &\leq \sum_{C \in \clusters(\cP_{r})\cap \Gamma} \frac{\diam_p(C)}{2} \left(\sqrt{\sigma_C^2 + \left(\frac{\sigma_C}{2}\left(\sqrt{\frac{b_C}{h}} - \sqrt{\frac{h}{b_C}}\right)\right)^2} - \frac{\sigma_C}{2}\left(\sqrt{\frac{b_C}{h}} - \sqrt{\frac{h}{b_C}}\right)\right) + \sum_{C \in \clusters(\cP_{r})\setminus \Gamma} \frac{\diam_p(C)}{2}  \sigma_C\\
    &= \sum_{C \in \clusters(\cP_{r})\cap \Gamma} \frac{b_C + h}{2} \left(\sqrt{\sigma_C^2 + \left(\frac{\sigma_C}{2}\left(\sqrt{\frac{b_C}{h}} - \sqrt{\frac{h}{b_C}}\right)\right)^2} - \frac{\sigma_C}{2}\left(\sqrt{\frac{b_C}{h}} - \sqrt{\frac{h}{b_C}}\right)\right) + \sum_{C \in \clusters(\cP_{r})\setminus \Gamma} \frac{\diam_p(C)}{2}  \sigma_C\\
    &= \sum_{C \in \clusters(\cP_{r})\cap \Gamma} \frac{\sigma_C}{2} \left(\sqrt{b_C h} + h \sqrt{\frac{h}{b_C}}\right) +\sum_{C \in \clusters(\cP_{r})\setminus \Gamma} \frac{\diam_p(C)}{2}  \sigma_C
\end{align*}
where the inequality follows again by the fact that $x \rightarrow \sqrt{a^2 + x^2} - x$ is non-increasing for every $a \in \mathbb{R}$. Hence,
$$
\sum_{r=1}^{R-1} s_r(\bq^{\textsf{GSM}}) \leq \sum_{C \in \Gamma} \frac{\sigma_C}{2}\left(\sqrt{b_C h} + h \sqrt{\frac{h}{b_C}}\right) + \sum_{C \in \clusters(\cH) \setminus \Gamma} \frac{\diam_p(C)}{2}  \sigma_C
$$

Hence, combining all of the above inequalities gives,
\begin{align*}
    \sup_{\cD \in \cF(\bmu, \Sigma)} \mathbb{E}_{\bd \sim \cD}\left(C^H(\bq, \bd)\right)
    &\leq 2\sigma_X \sqrt{bh} + \sum_{C \in \Gamma}\sigma_C \sqrt{b_Ch} + \sum_{C \in \clusters(H) \setminus \Gamma} \frac{\diam_p(C)}{2}\sigma_C
\end{align*}

\subsection{Proof of Lemma~\ref{lem:LBqd}}

The first inequality is trivial form the expression of $C(\bq, \bd)$ by noticing that $\Theta(\bq, \bd) \geq 0$ and that $h\left(q_X - d_X\right)^+ =  h\left(q_X - d_X\right) + h(d_X-q_X)^+$.

For the second, note that for every $C \in \cF$, the locations inside $C$ have at most $\overline{q}_{C}$ inventory available for them at a distance $\frac{\Delta_r}{2}$. Past this threshold, demand must be fulfilled from inventory at a distance of at least $\frac{\Delta_r}{2}$ or stay unfulfilled. Since the total unfulfilled demand in any optimal fulfillment is at most $(d_X-q_X)^+$, and since the clusters within a family are disjointed, in any optimal fulfillment, the amount of demand that is fulfilled from inventory at a distance at least $\frac{\Delta_r}{2}$ must be at least,
$
\sum_{C \in \cF}(d_C - \overline{q}_{C})^+ - (d_X-q_X)^+.
$
Hence,
\begin{align*}
    C(\bq, \bd) 
    &\geq h(q_X-d_X)^++b(d_X-q_X)^++\frac{\Delta_r}{2}\sum_{C \in \cF}(d_C - \overline{q}_{C})^+ - \frac{\Delta_r}{2}(d_X-q_X)^+\\
    &\geq h(q_X-d_X)+\frac{\Delta_r}{2}\sum_{C \in \cF}(d_C - \overline{q}_{C})^+ + \left(b+h-\frac{\Delta_r}{2}\right)(d_X-q_X)^+\\
    &\geq h(q_X-d_X)+\frac{\Delta_r}{2}\sum_{C \in \cF}(d_C - \overline{q}_{C})^+
\end{align*}
where the last inequality follows from the fact that $\Delta_r \leq \Delta_{R-1} = \frac{\Delta_R}{\gamma}\leq \diam(X) \leq (b+h)$ (note that $\Delta_R < \gamma\diam(X)$ by definition of $R$).

\subsection{Proof of Lemma~\ref{lem:ExpBound}}

We use Lemma~\ref{lem:chebyshev} with 
$$
\alpha = \frac{\sqrt{\sigma_S^2 + \left(Q - \mu_S\right)^2} + (Q- \mu_S)}{\sigma_S}.
$$
Note that because $Q \geq \mu_S$, $\alpha \geq 1$. Hence, by Lemma~\ref{lem:chebyshev}, there exists a distribution $\cD \in \cF_{S}(\bmu, \Sigma)$ such that,
$$
\mathbb{P}_{\bd \sim \cD}\left(\bd \geq {\bmu_S} + \frac{\alpha}{{\sigma_S}} {\bsigma_S^2}\right) = \Omega\left(\frac{1}{1+\alpha^2}\right)
$$
Next, let $\alpha' = 1$. Then, by the same lemma, there exists a distribution $\cD' \in \cF_{T \setminus S}(\bmu, \Sigma)$ such that,
$$
\mathbb{P}_{\bd' \sim \cD'}\left({\bd'} \geq {\bmu_{T\setminus S}} + \frac{\alpha'}{{\sigma_{T\setminus S}}} {\bsigma_{T\setminus S}^2}\right) = \Omega\left(1\right),
$$
implying that,
$$
\mathbb{P}_{\bd' \sim \cD'}\left({\bd'} \geq {\bmu_{T\setminus S}} \right) = \Omega(1).
$$
Now let $\cD'' \in \cF_{T}(\bmu, \Sigma)$ such that sampling $\bd'' \sim \cD''$ consists of sampling $\bd \sim \cD$ and $\bd' \sim \cD'$ independently and letting $\bd''$ take $\bd$ on $S$ and take $\bd'$ on $T\setminus S$. Under such distribution, we have,
\begin{align*}
    &\mathbb{E}_{\bd'' \sim \cD''}\left((\sum_{i \in S} d''_i - Q)^+ \cdot \mathbf{1}_{\{d_i \geq \mu, \forall i \in T\}}\right) \\
        &= \Omega\left( \frac{1}{1+\alpha^2} \left({\mu_S} + \frac{\alpha}{{\sigma_S}}\cdot {\sigma_S^2} - Q\right)^+\right)\\ &=\Omega\left(\frac{{\sigma_S^2}}{{\sigma_S^2} + \left(\sqrt{{\sigma_S^2}+(Q-{\mu_S})^2} + (Q-{\mu_S})\right)^2} \cdot \left({\mu_S} + \sqrt{{\sigma_S^2}+(Q-{\mu_S})^2} + (Q-{\mu_S}) - Q\right)^+\right)\\
        &=\Omega\left( \frac{{\sigma_S^2} \sqrt{{\sigma_S^2}+(Q-{\mu_S})^2}}{{\sigma_S^2} + \left(\sqrt{{\sigma_S^2}+(Q-{\mu_S})^2} + (Q-{\mu_S})\right)^2}\right) =\Omega\left(\sqrt{{\sigma_S^2}+(Q-{\mu_S})^2} - (Q-{\mu_S})\right)
\end{align*}

\subsection{Proof of Lemma~\ref{lem:LBzStar}}

    Let $\bq \in \cQ$. We have, 
    \begin{align*}
        \sup_{\cD \in \cF(\bmu, \Sigma)} \mathbb{E}_{\bd\sim \cD}\; C(\bq, \bd) 
        &\geq h(q_X - \mu_X) + (b+h) \sup_{\cD \in \cF(\bmu, \Sigma)} \mathbb{E}_{\bd \sim \cD}\left( d_X - q_X\right)^+\\
        &\geq \Omega\left(h(q_X - \mu_X) + \frac{b+h}{2} \left(\sqrt{\sigma_X^2 + \left(q_X - \mu_X\right)^2} - \left(q_X - \mu_X\right)\right)\right)\\
        &\geq \Omega\left(\sigma_X\sqrt{bh}\right)
    \end{align*}
    where the first inequality follows from Lemma~\ref{lem:LBqd}, the second from Lemma~\ref{lem:ExpBound} applied to $S=X$ and $T=X$, and the third is simply by taking the minimum of the scalar function over $q_X$.
    Next, for every $r \in \{1, \dots, R-1\}$ and $\cF \in \families(\cP_r)$, we have,
    \begin{align*}
        \sup_{\cD \in \cF(\bmu, \Sigma)} \mathbb{E}_{\bd\sim \cD}\; C(\bq, \bd) 
        &\geq h(q_X - \mu_X) + \frac{\Delta_r}{2} \sup_{\cD \in \cF(\bmu, \Sigma)} \mathbb{E}_{\bd \sim \cD} \left(\sum_{C \in \cF} \left(d_{C} - \overline{q}_{C}\right)^+\right)\\
        &\geq h(q_X - \mu_X) + \frac{\Delta_r}{2} \sum_{C \in \cF} \sup_{\cD \in \cF_{\overline{C}}(\bmu, \Sigma)} \; \mathbb{E}_{\bd \sim \cD} \left(\left(d_{C} - \overline{q}_{C}\right)^+\right)\\
        &\geq h(q_X - \mu_X) + \frac{\Delta_r}{2} \sum_{C \in \cF} \sup_{\cD \in \cF_{\overline{C}}(\bmu, \Sigma)} \; \mathbb{E}_{\bd \sim \cD} \left(\left(d_{C} - \overline{q}_{C}\right)^+ \bm{1}_{\{d_i \geq \mu_i, \forall i \in \overline{C}\}}\right)\\
        &\geq h(q_X - \mu_X) + \frac{\Delta_r}{2} \sum_{C \in \cF} \sup_{\cD \in \cF_{\overline{C}}(\bmu, \Sigma)} \; \mathbb{E}_{\bd \sim \cD} \left(\left(d_{C} - q_{\overline{C}} + \mu_{\overline{C} \setminus C}\right)^+ \bm{1}_{\{d_i \geq \mu_i, \forall i \in \overline{C}\}}\right)\\
        &\geq \Omega\left(h(q_X - \mu_X) + \frac{\Delta_r}{2} \sum_{C \in \cF} \left(\sqrt{\sigma_C^2 + \left(q_{\overline{C}} - \mu_{\overline{C}}\right)^2} - \left(q_{\overline{C}} - \mu_{\overline{C}}\right)\right)\right)\\
        &\geq \sum_{C\in \cF} \Omega\left( h\left( q_{\overline{C}} - \mu_{\overline{C}}\right) + \frac{\Delta_r}{2} \left(\sqrt{\sigma_C^2 + \left(q_{\overline{C}} - \mu_{\overline{C}}\right)^2} - \left(q_{\overline{C}} - \mu_{\overline{C}}\right)\right)\right)\\
        &\geq \sum_{C\in \cF} \Omega\left( h\left( q_{\overline{C}} - \mu_{\overline{C}}\right) + \frac{\Delta_{r+1}}{2\gamma} \left(\sqrt{\sigma_C^2 + \left(q_{\overline{C}} - \mu_{\overline{C}}\right)^2} - \left(q_{\overline{C}} - \mu_{\overline{C}}\right)\right)\right)\\
        &\geq \sum_{C\in \cF} \Omega\left( h\left( q_{\overline{C}} - \mu_{\overline{C}}\right) + \frac{\diam_p(C)}{2\alpha\gamma} \left(\sqrt{\sigma_C^2 + \left(q_{\overline{C}} - \mu_{\overline{C}}\right)^2} - \left(q_{\overline{C}} - \mu_{\overline{C}}\right)\right)\right)\\
        &\geq \frac{1}{\alpha\gamma}\sum_{C\in \cF} \Omega\left( h\left( q_{\overline{C}} - \mu_{\overline{C}}\right) + \frac{\diam_p(C)}{2} \left(\sqrt{\sigma_C^2 + \left(q_{\overline{C}} - \mu_{\overline{C}}\right)^2} - \left(q_{\overline{C}} - \mu_{\overline{C}}\right)\right)\right)\\
        &\geq \frac{1}{\alpha\gamma}\Omega\left(\sum_{C \in \cF \cap \Gamma} \sigma_C \sqrt{b_C h} + \sum_{C \in \cF \setminus \Gamma} \sigma_C \diam_p(C)\right)
    \end{align*}
    The first inequality follows from Lemma~\ref{lem:LBqd}. The second inequality follows by the fact that the sets $\overline{C}$ are disjointed for $C \in \cF$ and hence, given a maximal distribution $\cD_{\overline{C}} \in \cF_{\overline{C}}(\mu, \Sigma)$ for each $C$, one can construct a distribution $\cD_X \in \cF_{X}(\mu, \Sigma)$ such that  
    $$
    \mathbb{E}_{\bd \sim \cD_X} \left(\sum_{C \in \cF}\left(d_{C} - \overline{q}_{C}\right)^+\right) = \sum_{C \in \cF} \sup_{\cD \in \cF_{\overline{C}}(\bmu, \Sigma)} \; \mathbb{E}_{\bd \sim \cD} \left(\left(d_{C} - \overline{q}_{C}\right)^+\right)
    $$
    as follows: sampling $\bd \sim \cD_X$ samples independently $\bd_{\overline C} \sim \cD_{\overline C}$ for each set $\overline{C}$ (and $\bd' \sim \cD'$ for any arbitrary $\cD' \in \cF_{X \setminus \cup_{C \in \cF} C}(\bmu, \Sigma)$) then $\bd$ takes $\bd_{\overline C}$ on each set $\overline C$ (and takes $\bd'$ for the remaining points). The fourth inequality is by using that $d_i \geq \mu_i$ for every $i \in \overline{C} \setminus C$. The next inequality is by Lemma~\ref{lem:ExpBound} applied to $S=C$ and $T=\overline{C}$ for each cluster $C \in \cF$. The following inequality is by noticing that $q_{X} - \mu_{X} \geq \sum_{C \in \cF} (q_{\overline C} - \mu_{\overline C})$ as $\overline{C}$'s are disjointed. The last inequality is by minimizing each scalar function over $q_{\overline{C}} \geq \mu_{\overline{C}}$, where the minimum is $q_{\overline{C}} = \mu_{\overline C}$ when $C \notin \Gamma$ and the minimum is $q_{\overline{C}} = \frac{\sigma_C}{2} \left(\sqrt{\frac{\diam_p(C)-h}{h}} - \sqrt{\frac{h}{\diam_p(C)-h}}\right)$ otherwise.
    
    The lemma follows by taking the average over all $\cF \in \cP_r$ (at most $\beta$ many families) and over all $r \in \{1, \dots, R-1\}$.

\section{The Hierarchical Balance (HB) Policy: Omitted Proofs}

\subsection{Proof of Lemma~\ref{lem:InputSize}}

First, note that because we are selecting $\delta d$ to be the minimum of the remaining demand at location $i$ and the maximum demand that can be fulfilled before some location $j$ runs out of inventory, at each iteration of the while loop, either a location runs out of inventory or the demand of the current time round $t$ is completely fulfilled at some location. Hence, the while loop does at most $2n$ iterations. 

We now show that the algorithm can be implemented so that $\bd^{\textsf{curr}}$ and $\bq^{\textsf{curr}}$ have a size that is polynomial in the input size at every step of the algorithm, where a step is a period between two times the algorithm is at line~\ref{linestep} (we adopt the following conventions: the first step starts from the first time the algorithm is at line~\ref{linestep}, and, after each iteration of the while loop we go back to line~\ref{linestep} only if the condition of the loop holds). Note that every other computation performed by the algorithm is either via simple elementary operations from these values (e.g. computing $\delta d$) or is trivially polynomial (e.g. computing the $k_j$ values). 

To show this, let $M$ and $N$ be any integers such that $M d^t_i$ is integer for every $t$ and $i$ and $M q^{\textsf{HB}}_i$ is integer for every $i$ and such that both these integers are less or equal than $N$. Let $l=1, \dots, L$ denote the steps of the algorithm. For every variable, we will add a super script $l$ to denote the value of the variable at the beginning of step $l$ (e.g. $d^{l, \textsf{curr}}_{X}$ denotes the value of $d^{\textsf{curr}}_{X}$ at line~\ref{linestep} at the beginning of step $l$). Note that $L \leq 2n T$ where $T$ denotes the number of time rounds as there are at most $2n$ iterations of the while loop in each time round. We show by induction that for each step $l$, it holds that
$
(n!)^{l n} M d^{l, \textsf{curr}}_i
\text{ and }
(n!)^{l n} M q^{l, \textsf{curr}}_i
$
are integers smaller or equal than $(n!)^{ln} N$.

For $l=1$, this is trivial as 
$(n!)^{n} M d^{1,\textsf{curr}}_i = (n!)^{n} M d^1_i$ and $(n!)^{n} M q^{1,\textsf{curr}}_i = (n!)^{n} M q^{\textsf{HB}}_i$ are both integers less or equal than $\leq (n!)^{n} N$ by definition of $N$ and $M$.

Suppose now this holds for some $1 \leq l \leq L-1$ we show it holds for $l+1$. 

For $\bd^{l+1, \textsf{curr}}$, we have: (i) if between step $l$ and $l+1$ we left the while loop, then 
$\bd^{l+1, \textsf{curr}} = \bd^{t^{l+1}}$, hence, 
$
(n!)^{(l+1)n}Md^{l+1, \textsf{curr}}_j = (n!)^{(l+1)n}Md^{t^{l+1}}_j
$
for every $j$, which is an integer less or equal than $(n!)^{(l+1)n}N$ by definition of $N$ and $M$. (ii) if we did not leave the while loop between $l$ and $l+1$ then we distinguish two sub-cases: (ii)-1 the if condition line~\ref{lineif} holds at step $l$ then, 
$d^{l+1, \textsf{curr}}_j = d^{l, \textsf{curr}}_j$ for $j \neq i^{l+1}$ hence by hypothesis,
\begin{align*}
    (n!)^{(l+1)n}Md^{l+1, \textsf{curr}}_{j} 
    = (n!)^n(n!)^{ln}Md^{l, \textsf{curr}}_{j}\leq (n!)^{(l+1)n}N
\end{align*}
for $j \neq i^{l+1}$ and is integer, and, $$d^{l+1, \textsf{curr}}_{i^{l+1}} = d^{l, \textsf{curr}}_{i^{l+1}} - \min(d^{l, \textsf{curr}}_{i^{l+1}}, q^{l, \textsf{curr}}_{i^{l+1}})$$
hence by hypothesis,
\begin{align*}
    (n!)^{(l+1)n}Md^{l+1, \textsf{curr}}_{i^{l+1}} 
    &= (n!)^{(l+1)n}Md^{l, \textsf{curr}}_{i^{l+1}} - \min\{(n!)^{(l+1)n}d^{l, \textsf{curr}}_{i^{l+1}}, (n!)^{(l+1)n}q^{l, \textsf{curr}}_{i^{l+1}}\}\\
    &\leq (n!)^n(n!)^{ln}Md^{l, \textsf{curr}}_{i^{l+1}}\\
    &\leq (n!)^{(l+1)n}N
\end{align*}
and is integer.
(ii)-2 if we did not leave the loop between $l$ and $l+1$ and the if condition line~\ref{lineif} does not hold at step $l$ then,
$
 d^{l+1, \textsf{curr}}_j = d^{l, \textsf{curr}}_j$
for every $j \neq i^{l+1}$ hence by hypothesis,
\begin{align*}
    (n!)^{(l+1)n}Md^{l+1, \textsf{curr}}_{j} 
    = (n!)^n(n!)^{ln}Md^{l, \textsf{curr}}_{j}\leq (n!)^{(l+1)n}N
\end{align*}
for every $j \neq i^{l+1}$ and is integer,
and 
\begin{align*}
    d^{l+1, \textsf{curr}}_{i^{l+1}} 
    = d^{l, \textsf{curr}}_{i^{l+1}} - \min\{d^{l, \textsf{curr}}_{i^{l+1}}, \; \min_{j \in C^{l+1}} k^{l+1}_{j} \cdot q^{l, \textsf{curr}}_j\}.
\end{align*}
Hence,
\begin{align*}
    (n!)^{(l+1)n}Md^{l+1, \textsf{curr}}_{i^{l+1}} 
    &= (n!)^{(l+1)n}Md^{l, \textsf{curr}}_{i^{l+1}} - \min\{(n!)^{(l+1)n}Md^{l, \textsf{curr}}_{i^{l+1}}, \; \min_{j \in C^{l+1}} \{(n!)^{(l+1)n}M k^{l+1}_{j} \cdot q^{l, \textsf{curr}}_j\}\}\\
   & \leq (n!)^n (n!)^{ln}Md^{l, \textsf{curr}}_{i^{l+1}}\\
   & \leq (n!)^{(l+1)n}N 
\end{align*}
and is integer.

Now for $\bq^{l+1, \textsf{curr}}$, we have: (i) if all locations run out of inventory between step $l$ and $l+1$, then 
$\bq^{l+1, \textsf{curr}} = \bzero$, hence, 
$
(n!)^{(l+1)n}Mq^{l+1, \textsf{curr}}_j = 0
$
for every $j$, which is an integer less or equal than $(n!)^{(l+1)n}N$. (ii) if there is still inventory at some location at the beginning of step $l+1$, we distinguish two sub-cases: (ii)-1 the if condition line~\ref{lineif} holds at step $l$ then, 
$q^{l+1, \textsf{curr}}_j = q^{l, \textsf{curr}}_j$ for $j \neq i^{l+1}$ hence by hypothesis,
\begin{align*}
    (n!)^{(l+1)n}Mq^{l+1, \textsf{curr}}_{i^{l+1}} 
    = (n!)^n(n!)^{ln}Mq^{l, \textsf{curr}}_{i^{l+1}}\leq (n!)^{(l+1)n}N
\end{align*}
and is integer, and, $$q^{l+1, \textsf{curr}}_{i^{l+1}} = q^{l, \textsf{curr}}_{i^{l+1}} - \min(d^{l, \textsf{curr}}_{i^{l+1}}, q^{l, \textsf{curr}}_{i^{l+1}})$$
hence by hypothesis,
\begin{align*}
    (n!)^{(l+1)n}Mq^{l+1, \textsf{curr}}_{i^{l+1}} 
    &= (n!)^{(l+1)n}Mq^{l, \textsf{curr}}_{i^{l+1}} - \min\{(n!)^{(l+1)n}d^{l, \textsf{curr}}_{i^{l+1}}, (n!)^{(l+1)n}q^{l, \textsf{curr}}_{i^{l+1}}\}\\
    &\leq (n!)^n(n!)^{ln}Mq^{l, \textsf{curr}}_{i^{l+1}}\\
    &\leq (n!)^{(l+1)n}N
\end{align*}
and is integer.
(ii)-2 the if condition line~\ref{lineif} does not hold at step $l$ then,
$q^{l+1, \textsf{curr}}_j = q^{l, \textsf{curr}}_j$
for every $j \notin C^{l+1}$ or $j \in C^{l+1}$ s.t. $q^{l, \textsf{curr}}_j = 0$, hence by hypothesis,
\begin{align*}
    (n!)^{(l+1)n}Mq^{l+1, \textsf{curr}}_{i^{l+1}} 
    = (n!)^n(n!)^{ln}Mq^{l, \textsf{curr}}_{i^{l+1}} \leq (n!)^{(l+1)n}N
\end{align*}
for $j \notin C^{l+1}$ or $j \in C^{l+1}$ s.t. $q^{l, \textsf{curr}}_j = 0$ and is integer,
and 
\begin{align*}
    q^{l+1, \textsf{curr}}_{j} 
    = q^{l, \textsf{curr}}_{j} - \frac{1}{k^{l+1}_j} \min\{d^{l, \textsf{curr}}_{i^{l+1}}, \; \min_{j' \in C^{l+1}} k^{l+1}_{j'} \cdot q^{l, \textsf{curr}}_{j'}\}
\end{align*}
for $j \in C^{l+1}$ s.t. $q^{l, \textsf{curr}}_j > 0$, hence by hypothesis,
\begin{align*}
    (n!)^{(l+1)n}Mq^{l+1, \textsf{curr}}_{j} 
    &= (n!)^{(l+1)n}Mq^{l, \textsf{curr}}_{j} - \frac{(n!)^n}{k^{l+1}_{j}} \min\{(n!)^{ln}Md^{l, \textsf{curr}}_{i^{l+1}}, \; \min_{j' \in C^{l+1}} \{(n!)^{ln}M k^{l+1}_{j'} \cdot q^{l, \textsf{curr}}_{j'}\}\}\\
   & \leq (n!)^n (n!)^{ln}Mq^{l, \textsf{curr}}_{j}\\
   & \leq (n!)^{(l+1)n}N 
\end{align*}
for $j \in C^{l+1}$ s.t. $q^{l, \textsf{curr}}_j > 0$ and is integer. Note that $k^{l+1}_{j}$ divides $(n!)^n$. This proves our claim for $l+1$. Now, if we let $M$ be the smallest common multiple of the denominators of all $d^t_i$ for all $t$ and $i$ and $q^{\textsf{HB}}_i$ for all $i$ and we let $N$ be the maximum of $Md^t_i$ for all $t$ and $i$ and $Mq^{\textsf{HB}}_i$ for all $i$ then $M$ and $N$ are polynomial sized in the input since $\bq^{\textsf{HB}}$ is polynomial sized as a (basic feasible) solution of an LP with polynomial sized data. Hence, the bit size of $d^{\textsf{curr}}_i$ and $q^{\textsf{curr}}_i$ when represented in their reduced form (the nominator and denominator are co-prime and hence the denominator must divide $(n!)^{Ln}M$ and nominator must be at most $(n!)^{Ln}N$) for every $i$ at every step of the algorithm is at most $O(\log((n!)^{Ln}M)+\log((n!)^{Ln}N)) = O(n^3 T \log(n)(\log(MN)))$ which is polynomial in the input size. Now we can make sure to always store $d^{\textsf{curr}}_i$ and $q^{\textsf{curr}}_i$ in their reduced form (after each update of these quantities, if they are not in the reduced form, their reduced form can be computed in polynomial time using Euclid's algorithm) and get a polynomial time implementation of the algorithm.

\subsection{Proof of Lemma~\ref{lem:reductionToC}}

    Assume
    $
    C^H_{\cA^{\textsf{HB}}}(\bq^{\textsf{HB}}, \pi) \leq O(\log(n)) \cdot C^H(\bq^{\textsf{HB}}, \bd) 
    $
    for every $\bd$ such that $d_X = q^{\textsf{HB}}_X$ and $\pi \in \Pi(\bd)$. We show that for every $\bd \geq \bzero$ and $\pi \in \Pi(\bd)$, the inequality
    $
    C^H_{\cA^{\textsf{HB}}}(\bq^{\textsf{HB}}, \pi) \leq O(\log(n)) \cdot C^H(\bq^{\textsf{HB}}, \bd) 
    $
    still holds. This would imply that,
    \begin{align*}
        \sup_{\cD \in \cF(\bmu, \Sigma)} \mathbb{E}_{\bd \sim \cD}\; \left(\sup_{\pi \in \Pi(\bd)}
        C_{\cA^{\textsf{HB}}}(\bq^{\textsf{HB}}, \pi)
        \right)
        &\leq \sup_{\cD \in \cF(\bmu, \Sigma)} \mathbb{E}_{\bd \sim \cD}\; \left(\sup_{\pi \in \Pi(\bd)}
        C^H_{\cA^{\textsf{HB}}}(\bq^{\textsf{HB}}, \pi)
        \right)\\
        &\leq 
        O(\log(n))\sup_{\cD \in \cF(\bmu, \Sigma)} \mathbb{E}_{\bd \sim \cD}\; \left(
        C^H(\bq^{\textsf{HB}}, \bd)
        \right)\\
        &\leq O(\log(n)) \left(2\sigma_X \sqrt{bh} + \sum_{C \in \Gamma}\sigma_C \sqrt{b_Ch} + \sum_{C \in \clusters(H) \setminus \Gamma} \frac{\sigma_C\diam_p(C)}{2}\right)\\
        &\leq 
        O(\alpha\beta\gamma R \log(n)) \cdot z_{\textsf{ODRNM}}
    \end{align*}
    where the third inequality follows from Lemma~\ref{lem:UB} and the last one from inequality~\eqref{eq:LB}. Hence that $\bq^{\textsf{HB}}$ and $\cA^{\textsf{HB}}$ achieve the bounds in Lemma~\ref{lem:main2}.

    Let $\bd \geq \bzero$ and $\pi \in \Pi(\bd)$. If $d_X > q^{\textsf{HB}}_X$, note that by construction $\cA^{\textsf{HB}}$ will fulfill a total of exactly $q^{\textsf{HB}}_X$ of demand and the rest is left unfulfilled. Let $\bd'$ denote the sub-vector of $\bd$ of demand that is fulfilled by $\cA^{\textsf{HB}}$, let $t$ denote the time step that the last demand part $\bd^t$  has been fulfilled fully or partially by the algorithm (the last part might be fulfilled only partially), let $\bd^{' t}$ denote the sub-portion of $\bd^{t}$ that has been fulfilled by the algorithm at time $t$, then the algorithm will take the exact same decisions under the sequence of arrivals $\pi' = (\bd^1, \dots, \bd^{t-1}, \bd^{' t})$ and we have, 
\begin{align*}
    C^H_{\cA^{\textsf{HB}}}(\bq^{\textsf{HB}}, \pi)
    &= b(d_X-d'_X) + C^H_{\cA^{\textsf{HB}}}(\bq^{\textsf{HB}}, \pi')\\
    &= h(q^{\textsf{HB}}_X-d_X)^+ + b(d_X-q^{\textsf{HB}}_X)^+ + C^H_{\cA^{\textsf{HB}}}(\bq^{\textsf{HB}}, \pi')\\
    &\leq O(\log(n)) \left(h(q^{\textsf{HB}}_X-d_X)^+ + b(d_X-q^{\textsf{HB}}_X)^+ + C^{H}(\bq^{\textsf{HB}}, \bd')\right)\\
    &= O(\log(n)) \left(h(q^{\textsf{HB}}_X-d_X)^+ + b(d_X-q^{\textsf{HB}}_X)^+\right.
    \\
    &\quad \left.+\sum_{r=1}^{R-1} \sum_{C \in \clusters(\cP_{r})} \diam_p(C) \cdot (d'_C - q^{\textsf{HB}}_C)^+ + \frac{\diam(X)}{n}\sum_{i \in X}(d'_i-q^{\textsf{HB}}_i)^+\right)\\
    &\leq O(\log(n)) \left(h(q^{\textsf{HB}}_X-d_X)^+ + b(d_X-q^{\textsf{HB}}_X)^+\right.
    \\
    &\quad \left.+\sum_{r=1}^{R-1} \sum_{C \in \clusters(\cP_{r})} \diam_p(C) \cdot (d_C - q^{\textsf{HB}}_C)^+ + \frac{\diam(X)}{n}\sum_{i \in X}(d_i-q^{\textsf{HB}}_i)^+\right)\\
    &= O(\log(n))\cdot C^H(\bq^{\textsf{HB}}, \bd)
\end{align*}
Now, if $d_X < q^{\textsf{HB}}_X$, it must be the case that there exists a cluster $C_{R-1} \in \clusters(\cP_{R-1})$ that has strictly more inventory than demand. Then again, there must exist a cluster $C_{R-2} \in \clusters(C_{R-1})$ that has strictly more inventory than demand, and so on until we reach a cluster $C_{1} \in \clusters(C_{2})$ with strictly more inventory then demand which implies in turn that one of the locations $i \in C_1$ has strictly more inventory than demand. Now reducing the inventory at location $i$ by a small amount $\delta$ can only influence the terms $(d_{C_{R-1}}-q^{\textsf{HB}}_{C_{R-1}})^+, \dots, (d_{C_{1}}-q^{\textsf{HB}}_{C_{1}})^+$, $(d_{i}-q^{\textsf{HB}}_{i})^+$, $(d_X-q^{\textsf{HB}}_X)^+$ and $(q^{\textsf{HB}}_X-d_X)^+$ in the expression of $C^H(\bq^{\textsf{HB}}, \bd)$. The first $R+2$ terms are zero and will remain zero if the inventory is reduced by a small enough amount, the last term will decrease by $\delta$. Hence, the new inventory vector obtained $\bq^\delta$ is such that $C^H(\bq^\delta,\bd) = C^H(\bq^{\textsf{HB}}, \bd) - h \delta$. Iterating this argument, this implies the existence of a vector $\bq' \leq \bq^{\textsf{HB}}$ such that $q'_X = d_X$ and $
C^H(\bq',\bd) = C^H(\bq^{\textsf{HB}}, \bd) - h (q^{\textsf{HB}}_X - q'_X).$
Let $\bd' = \bd + (\bq^{\textsf{HB}} - \bq')$ and let $\pi'$ be the sequence $\pi$ to which we add a demand part $\bq^{\textsf{HB}} - \bq'$ arriving last. 
We have,
\begin{align*}
    C^H_{\cA^{\textsf{HB}}}(\bq^{\textsf{HB}}, \pi)
    &\leq h(q^{\textsf{HB}}_X - q'_X) + C^H_{\cA^{\textsf{HB}}}(\bq^{\textsf{HB}}, \pi') \\
    &\leq O(\log(n)) (h(q^{\textsf{HB}}_X - q'_X) + C^H(\bq^{\textsf{HB}}, \bd'))\\
    &= O(\log(n))(h(q^{\textsf{HB}}_X - q'_X) + C^H(\bq', \bd))\\
    &= O(\log(n)) \cdot C^H(\bq^{\textsf{HB}}, \bd).
\end{align*}
The first inequality follows from the fact that $C^H_{\cA^{\textsf{HB}}}(\bq^{\textsf{HB}}, \pi)$ is equal to the overage cost of $h(q^{\textsf{HB}}_X - q'_X)$ at the end of the online algorithm plus the cost (under $\ell^H$) of fulfilling demand online, but since the sequence $\pi'$ has the sequence $\pi$ as a prefix, the online fulfillment policy will make the exact same decisions when fulfilling the first sub-sequence $\pi$ of $\pi'$ and hence will pay at least the same fulfillment cost. The second inequality is by hypothesis as $q^{\textsf{HB}}_X = d'_{X}$, the next equality follows from the fact that for every $\bd, \bq$ and $\bu$ it is holds that $C^H(\bq + \bu, \bd + \bu) = C^H(\bq, \bd)$ ($C^H(\bq, \bd)$ is a function of the vector $(d_i - q_i)_{i \in X}$ which is invariant by translation of $\bd$ and $\bq$ by the same vector).

\subsection{Proof of Lemma~\ref{lem:comparepotential}}

For every $C' \in V_C$, let 
$
\theta^{\textsf{online}, C}_{[1:l-1]}(C') 
$
denote the total amount of demand from $\delta d^1, \dots, \delta d^{l-1}$ that arrived inside $C'$ and that the online algorithm fulfilled using inventory outside of $C'$, where by convention $
\theta^{\textsf{online}, C}_{[1:0]}(C') = 0
$.

Next, consider the inventory remaining after the online algorithm fulfilled $\delta d^1, \dots, \delta d^{l-1}$. Consider an offline algorithm that fulfills the remaining portions $\delta d^{l}, \dots, \delta d^{L}$ inside $C$ using this remaining inventory in a way that fulfills a maximum amount of demand to the remaining inventory inside each $C' \in V_C$, then fulfills the remaining demand arbitrary inside $C$. Then let,
$
\theta^{\textsf{offline}, C}_{[l:L]}(C') 
$
denote the amount of demand from $\delta d^{l}, \dots, \delta d^{L}$ that the offline algorithm fulfiled outside of $C' \in V_C$, where by convention $
\theta^{\textsf{offline}, C}_{[L+1:L]}(C') = 0
$. Note that for every $l \in \{1, \dots, L+1\}$,
$$
\theta^{\textsf{offline}, C}_{[l:L]}(C') = \left(\sum_{k:\; k\geq l,\; i^{k} \in C'} \delta d^k - q^{l}_{C'}\right)^+
$$
where $\bq^l$ denotes the remaining inventory vector after the online algorithm fulfilled the first $\delta d^1, \dots, \delta d^{l-1}$.

Let $C(0), C(1), \dots, C(U)$ denote the order by which the elements of $V_C$ run out of inventory (the virtual cluster $C(0)$ has no inventory and hence its first in the list). Consider the following potential function
$$
\Phi_l =  \sum_{u=0}^U \theta^{\textsf{online}, C}_{[1:l-1]}(C(u)) + \sum_{u=0}^U 3\ln(U-u+1) \cdot \theta^{\textsf{offline}, C}_{[l:L]}(C(u)).
$$
For every $l \in \{1,\dots, L+1\}$. We show that $\Phi_l$ is non-increasing in $l$. The lemma follows from the fact that 
\begin{align*}
    \Phi_{L+1} &= \sum_{C' \in V_C} \theta^{\textsf{online}, C}_{[1:L]}(C') = \sum_{C' \in V_C} \theta^{\textsf{online}, C}(C') \leq \Phi_{1} \leq 3\ln(n) \sum_{C' \in V_C} \theta^{\textsf{offline}, C}_{[1:L]}(C') = 3\ln(n) \sum_{C' \in V_C} \theta^{\textsf{offline}, C}(C').
\end{align*}

Fix $l \in \{1, \dots, L\}$
and suppose $\delta d^l$ arrived inside $C(u)$ for some $u \in \{0, \dots, U\}$. We now show that 
$
\Phi_{l+1} \leq \Phi_l.
$

Consider the step at which $\delta d^l$ was considered by the algorithm. There are two cases,

\vspace{3mm}{\noindent \bf Case 1.} There was available inventory to fulfill $\delta d^l$ inside $C(u)$. In this case, the online algorithm fulfills $\delta d^l$ using inventory inside $C(u)$, hence,
$$
\theta^{\textsf{online}, C}_{[1:l]}(C(u')) = \theta^{\textsf{online}, C}_{[1:l-1]}(C(u'))
$$
for every $u' \in \{0, \dots, U\}$. This also implies that the inventory within every $C(u')$ such that $u' \neq u$ is unchanged, i.e., $q^{l+1}_{C(u')} = q^l_{C(u')}$ for every $u' \neq u$, and the inventory within $C(u)$ decreases by $\delta d^l$ after the online algorithm fulfills $\delta d^l$, hence, $q^{l+1}_{C(u)} = q^l_{C(u)} - \delta d^l$. This further implies that,
\begin{align*}
    \theta^{\textsf{offline}, C}_{[l:L]}(C(u')) &= \left(\sum_{k:\; k\geq l,\; i^{k} \in C(u')} \delta d^k - q^{l}_{C(u')}\right)^+ = \left(\sum_{k:\; k\geq l+1,\; i^{k} \in C(u')} \delta d^k - q^{l+1}_{C(u')}\right)^+ = \theta^{\textsf{offline}, C}_{[l+1:L]}(C(u'))
\end{align*}
for $u' \neq u$ as $i^l \notin C(u')$, and 
\begin{align*}
    \theta^{\textsf{offline}, C}_{[l:L]}(C(u)) &= \left(\sum_{k:\; k\geq l,\; i^{k} \in C(u)} \delta d^k - q^{l}_{C(u)}\right)^+\\
    &= \left(\sum_{k:\; k\geq l+1,\; i^{k} \in C(u)} \delta d^k - q^{l+1}_{C(u)} + \delta d^{l} - (q^{l}_{C(u)}-q^{l+1}_{C(u)})\right)^+\\
    &= \left(\sum_{k:\; k\geq l+1,\; i^{k} \in C(u')} \delta d^k - q^{l+1}_{C(u')}\right)^+\\
    &= \theta^{\textsf{offline}, C}_{[l+1:L]}(C(u)).
\end{align*}
Hence,
$$
\theta^{\textsf{offline}, C}_{[l+1:L]}(C(u')) = \theta^{\textsf{offline}, C}_{[l:L]}(C(u'))
$$
for every $u' \in \{0, \dots, U\}$.
Hence, $\Phi_{l+1}=\Phi_l$ in this case.

\vspace{3mm}{\noindent \bf Case 2.} 
There is no available inventory to fulfill $\delta d^l$ inside $C(u)$. Note that this necessarily means that $C(u)$ has no left inventory, as $\delta d^l$ is an infinitesimal portion (or part of an infinitesimal portion if $u=0$). In this case, both the online algorithm and the offline algorithm will both have to fulfill $\delta d^l$ from outside $C(u)$. We have, 
$$
\theta^{\textsf{online}, C}_{[1:l]}(C(u)) = \theta^{\textsf{online}, C}_{[1:l-1]}(C(u)) + \delta d^l
$$
and 
$$
\theta^{\textsf{online}, C}_{[1:l]}(C(u')) = \theta^{\textsf{online}, C}_{[1:l-1]}(C(u'))
$$
for every other $u' \neq u$.

Next, let $u_l, u_{l}+1, \dots, U$ be the indices of the elements of $V_C$ that still have inventory at the moment $\delta d^l$ arrived (these are the clusters or locations from which the online algorithm will take an equal amount of inventory). Note that $u_l > u$ (since $C(u)$ has run out of inventory already). Note that,
\begin{itemize}
    \item For every $u' < u_l$ such that $u' \neq u$, the inventory inside $C(u')$ has not changed such that $q^{l+1}_{C(u')} = q^{l}_{C(u')} = 0$, and we have,
    \begin{align*}
        \theta^{\textsf{offline}, C}_{[l:L]}(C(u')) 
        &= \left(\sum_{k:\; k\geq l,\; i^{k} \in C(u')} \delta d^k - q^{l}_{C(u')}\right)^+  =\left(\sum_{k:\; k\geq l+1,\; i^{k} \in C(u')} \delta d^k - q^{l+1}_{C(u')}\right)^+ =\theta^{\textsf{offline}, C}_{[l+1:L]}(C(u')) 
    \end{align*}
    \item The inventory inside $C(u)$ also has not changed such that $q^{l+1}_{C(u)} = q^{l}_{C(u)} = 0$, and we have,
    \begin{align*}
        \theta^{\textsf{offline}, C}_{[l:L]}(C(u)) 
        &= \left(\sum_{k:\; k\geq l,\; i^{k} \in C(u)} \delta d^k - q^{l}_{C(u)}\right)^+ \\
        &= \sum_{k:\; k\geq l,\; i^{k} \in C(u)} \delta d^k\\
        &= \sum_{k:\; k\geq l+1,\; i^{k} \in C(u)} \delta d^k + \delta d^l\\
        &= \left(\sum_{k:\; k\geq l+1,\; i^{k} \in C(u)} \delta d^k - q^{l+1}_{C(u)}\right)^+ + \delta d^l\\
        &=\theta^{\textsf{offline}, C}_{[l+1:L]}(C(u')) + \delta d^l
    \end{align*}
    \item The inventory inside each $C(u')$ for $u' = u_l, u_l+1, \dots, U$ has decreased by $\frac{\delta d^l}{U-u_l+1}$, i.e., $q^{l+1}_{C(u')} = q^{l}_{C(u')} - \frac{\delta d^l}{U-u_l+1}$, and
    \begin{align*}
        \theta^{\textsf{offline}, C}_{[l:L]}(C(u')) 
        &= \left(\sum_{k:\; k\geq l,\; i^{k} \in C(u')} \delta d^k - q^{l}_{C(u')}\right)^+ \\
        &= \left(\sum_{k:\; k\geq l+1,\; i^{k} \in C(u')} \delta d^k - q^{l+1}_{C(u')} - \frac{\delta d^l}{U-u_l+1}\right)^+ \\
        &\geq \left(\sum_{k:\; k\geq l+1,\; i^{k} \in C(u')} \delta d^k - q^{l+1}_{C(u')}\right)^+ - \frac{\delta d^l}{U-u_l+1}\\
        &=\theta^{\textsf{offline}, C}_{[l+1:L]}(C(u')) - \frac{\delta d^l}{U-u_l+1}
    \end{align*}
\end{itemize}
Hence,
\begin{align*}
    &\Phi_{l+1}-\Phi_l\\
    &=\sum_{u'=0}^U \left(\theta^{\textsf{online}, C}_{[1:l]}(C(u')) - \theta^{\textsf{online}, C}_{[1:l-1]}(C(u'))\right) + 3\ln(U+1-u') \left(\theta^{\textsf{offline}, C}_{[l+1:L]}(C(u')) - \theta^{\textsf{offline}, C}_{[l:L]}(C(u'))\right)\\
    &\leq  \delta d^l - \delta d^l \cdot 3\ln(U+1-u) + \frac{\delta d^l}{U-u_l+1} \sum_{u'=u_l}^{U}  3\ln(U+1-u')\\
    &\leq \delta d^l - \delta d^l \cdot 3\ln(U+2-u_l) + \frac{\delta d^l}{U-u_l+1} \sum_{u'=u_l}^{U}  3\ln(U+1-u')\\
    &\leq \delta d^l - \delta d^l \cdot 3\ln(U+2-u_l) + \frac{\delta d^l}{U-u_l+1} \sum_{u'=1}^{U-u_l+1}  3\ln(u')
\end{align*}
where the second inequality follows from the fact that $u \leq u_l-1$. When $u_l = U$, then the above upper-bound is negative. Otherwise, let $x=U+2-u_l$ (which is hence $\geq 3$), we have,
\begin{align*}
    -3\ln(x) + \frac{1}{x-1}\sum_{u'=1}^{x-1} 3\ln(u')
    &\leq -3\ln(x) + \frac{1}{x-1}\int_{u'=1}^{x} 3\ln(u')\\
    &\leq -3\ln(x) + \frac{1}{x-1} \left(3x\ln(x) - 3x + 3\right) \leq -1
\end{align*}
where the last inequality follows from the fact that $x \geq 3$ and the function on $x$ in decreasing over $(3, \infty)$ with value smaller than $-1$ at $x=3$. This implies that $\Phi_{l+1}-\Phi_l \leq 0$.

\subsection{Proof of Lemma~\ref{lem:onlineUB}}
    Consider some amount $\delta$ of inventory that the online algorithm moved between locations $i$ and $j$ when fulfilling demand of an infinitesimal portion $\delta d$. Let $r$ be the smallest level such that $i$ and $j$ belong to the same cluster $C \in \clusters(\cP_r)$, and hence such that $\ell^H_{ij} = \diam(C)$. If $r\geq 2$, let $C^{j}$ be the cluster in $\cP_{r-1}$ to which $j$ belongs, otherwise let $C^j=j$. Then, on the one hand, the contribution of $\delta$ to the left-hand-side is $\ell^H_{ij} \cdot \delta = \delta \diam(C)$. On the other hand, it must be the case that $\delta d$ arrived inside $C^j$ and was fulfilled using inventory outside of $C^j$ (inside $C$) and hence $\delta$ is part of $\theta^{\textsf{online}, C}(C^j)$ with coefficient $\diam(C)$, hence a contribution to the right-hand-side of $\delta \diam(C)$.

\subsection{Proof of Lemma~\ref{lem:onlineLB}}
We have,
\begin{align*}
    &\sum_{C \in \clusters(\cH)} \diam(C) \sum_{C' \in V_C}\theta^{\textsf{offline}, C}(C')\\
    &= \sum_{r=2}^R \sum_{C \in \clusters(\cP_r)}\diam(C) \left(\sum_{C' \in \clusters(C)}\theta^{\textsf{offline}, C}(C') + \theta^{\textsf{offline}, C}(C(0))\right) \\
    &\quad + \sum_{C \in \clusters(\cP_1)} \diam(C) \left( \sum_{i \in C}\theta^{\textsf{offline}, C}(i) + \theta^{\textsf{offline}, C}(C(0))\right)\\\\
    &\leq \sum_{r=2}^R \sum_{C \in \clusters(\cP_r)}\diam(C) \left(\sum_{C' \in \clusters(C)}(d_{C'} - q^{\textsf{HB}}_{C'})^+ + \theta^{\textsf{offline}, C}(C(0))\right) \\
    &\quad + \sum_{C \in \clusters(\cP_1)} \diam(C) \left( \sum_{i \in C}(d_i - q^{\textsf{HB}}_{i})^+ + \theta^{\textsf{offline}, C}(C(0))\right)\\\\
    &= \sum_{r=2}^R \sum_{C \in \clusters(\cP_r)}\diam(C) \sum_{C' \in \clusters(C)}(d_{C'} - q^{\textsf{HB}}_{C'})^+ + \sum_{C \in \clusters(\cP_1)} \diam(C) \sum_{i \in C}(d_i - q^{\textsf{HB}}_{i})^+\\
    &\quad + \sum_{C \in \clusters(\cH)}\diam(C) \theta^{\textsf{offline}, C}(C(0))\\\\
    &\leq 2C^H(\bq^{\textsf{HB}}, \bd) + \sum_{C \in \clusters(\cH)}\diam(C) \theta^{\textsf{offline}, C}(C(0))
\end{align*}
The first inequality follows from the definition of $\theta^{\textsf{offline}, C}(C')$, in fact, when $C'$ is a cluster or a location, the expression is given by the positive part of some of the demand that arrived at $C'$ minus the inventory inside $C'$ and hence can be upper-bounded by $(d_{C'}-q^{\textsf{HB}}_{C'})^+$. For the second inequality recall that if $\Delta_1 = \frac{1}{\alpha}\min_{i,j \in X: i \neq j} \ell_{ij}$, then the clusters in the first partition are necessarily singletons and hence,
\begin{align*}
    &\sum_{r=1}^{R-1} \sum_{C \in \clusters(\cP_r)}\diam_p(C) (d_{C} - q^{\textsf{HB}}_{C})^+ + \sum_{C \in \clusters(\cP_1)} \diam(C) \sum_{i \in C}(d_i - q^{\textsf{HB}}_{i})^+\\
    &\leq 2 \sum_{r=1}^{R-1} \sum_{C \in \clusters(\cP_r)}\diam_p(C) (d_{C} - q^{\textsf{HB}}_{C})^+\\
    &\leq 2 C^H(\bq^{\textsf{HB}}, \bd)
\end{align*}
otherwise, the diameter of the clusters at level $1$ is at most $\diam(X)/n$ and hence,
\begin{align*}
    &\sum_{r=1}^{R-1} \sum_{C \in \clusters(\cP_r)}\diam_p(C) (d_{C} - q^{\textsf{HB}}_{C})^+ + \sum_{C \in \clusters(\cP_1)} \diam(C) \sum_{i \in C}(d_i - q^{\textsf{HB}}_{i})\\
    &\leq \sum_{r=1}^{R-1} \sum_{C \in \clusters(\cP_r)}\diam_p(C) (d_{C} - q^{\textsf{HB}}_{C})^+ + \frac{\diam(X)}{n}\sum_{i \in X}(d_i - q^{\textsf{HB}}_{i})\\
    &\leq C^H(\bq^{\textsf{HB}}, \bd).
\end{align*}
We now show that,
$$\sum_{C \in \clusters(\cH)}\diam(C) \theta^{\textsf{offline}, C}(C(0)) \leq \frac{2}{21\log(n)}\sum_{C \in \clusters(\cH)} \diam(C) \sum_{C' \in V_C}\theta^{\textsf{online}, C}(C')
$$
to conclude.

To prove this inequality, note that for a given level $r \in \{1, \dots, R-1\}$, we have,
\begin{align*}
    \sum_{C \in \clusters(\cP_r)} \diam(C)\theta^{\textsf{offline}, C}(C(0))
    &\leq 
    \alpha \Delta_r\sum_{C \in \clusters(\cP_r)} \theta^{\textsf{offline}, C}(C(0))\\
    &\leq 
    \alpha \Delta_r \sum_{r'=r+1}^R\sum_{\substack{C \in \clusters(\cP_{r'})\\ \text{s.t } |\clusters(C)| \geq 2}} \sum_{C' \in V_C}\theta^{\textsf{online}, C}(C')\\
    &\leq 
    \frac{1}{21\log(n)}\sum_{r'=r+1}^R \frac{1}{\gamma^{r'-r - 1}}\sum_{\substack{C \in \clusters(\cP_{r'})\\ \text{s.t } |\clusters(C)| \geq 2}} \Delta_{r'} \sum_{C' \in V_C}\theta^{\textsf{online}, C}(C')\\
    &\leq 
    \frac{1}{21\log(n)}\sum_{r'=r+1}^R \frac{1}{\gamma^{r'-r - 1}} \sum_{\substack{C \in \clusters(\cP_{r'})}} \diam(C) \sum_{C' \in V_C}\theta^{\textsf{online}, C}(C')
\end{align*}
where the first inequality follows by definition. For the second inequality, note that the sum $\sum_{C \in \clusters(\cP_r)} \theta^{\textsf{offline}, C}(C(0))$ is exactly the amount of demand that has been fulfilled by the online algorithm from inventory inside the clusters of $\cP_r$ but arrived outside of these clusters (i.e. demand arrived at some cluster of $\cP_r$ and fulfilled from another). Consider a portion $\delta$ of such demand, then there must exist a cluster $\hat{C}$ at a higher level $k \geq r+1$ such that $\delta$ arrived at one cluster $C^1$ in $\clusters(\hat{C})$ and has been fulfilled from another cluster $C^2$ in $\clusters(\hat{C})$ implying that $\delta$ is part of $\theta^{\textsf{online}, \hat{C}}(C^1)$. The third inequality follows from the fact that $\Delta_{r'}=\gamma^{r'-r}\Delta_r$ and by our choice that $21\log(n)\alpha < \gamma$. The fourth inequality follows from the fact that when $C$ of level $r'$ contains at least $2$ subclusters at level $r'-1$, its diameter must be at least as large as the distance (margin) between these two subclusters clusters i.e. $\diam(C) \geq \Delta_{r'}$.

We now sum the above inequality over all $r=1, \dots, R-1$ to get,
\begin{align*}
    \sum_{C \in \clusters(\cH)} \diam(C)\theta^{\textsf{offline}, C}(C(0))
    &\leq 
    \frac{1}{21\log(n)}\sum_{r=1}^{R-1}\sum_{r'=r+1}^R \frac{1}{\gamma^{r'-r - 1}} \sum_{\substack{C \in \clusters(\cP_{r'})}} \diam(C) \sum_{C' \in V_C}\theta^{\textsf{online}, C}(C')\\
    &= 
    \frac{1}{21\log(n)}\sum_{r'=2}^R\sum_{r=1}^{r'-1} \frac{1}{\gamma^{r'-r - 1}} \sum_{\substack{C \in \clusters(\cP_{r'})}} \diam(C) \sum_{C' \in V_C}\theta^{\textsf{online}, C}(C')\\
    &\leq 
    \frac{1}{1-\frac{1}{\gamma}}\frac{1}{21\log(n)}\sum_{r'=2}^R \sum_{\substack{C \in \clusters(\cP_{r'})}} \diam(C) \sum_{C' \in V_C}\theta^{\textsf{online}, C}(C')\\
    &\leq 
    \frac{2}{21\log(n)} \sum_{\substack{C \in \clusters(\cH)}} \diam(C) \sum_{C' \in V_C}\theta^{\textsf{online}, C}(C').
\end{align*}
The last inequality follows from $\gamma \geq 21\log(n) \geq 2$ for $n \geq 2$, implying that $1/(1-1/\gamma) \leq 2$.

\section{Derivation of~\eqref{SDP}}
\label{apx:derivation-sdp}
Consider a metric space as constructed in the the numerical experiments Section~\ref{sec:num}. The cost function $C(\bq, \bd)$ can be formulated as follows,
\begin{align*}
    C(\bq, \bd)
    &= h(q_X - d_X)^+ + b(d_X - q_X)^+ 
    \\
    &+ 2c \sum_{r=1}^{K-1} \lambda^{r} \left(\sum_{v \in L_r} (d_v - q_v)^+ - \sum_{w \in L_{r+1}} (d_w - q_w)^+\right)\\
    &= h(q_X - d_X) + (b+h-2c\lambda^{K-1}) (d_X - q_X)^+ 
    \\
    &+ 2c\lambda\sum_{v \in L_1} (d_v-q_v)^+ + 2c\sum_{r=2}^{K-1} (\lambda^r - \lambda^{r-1})\sum_{v \in L_{r}}(d_v-q_v)^+\\
    &= \max_{\bm{\epsilon}_1 \in \Pi_1, \dots, \bm{\epsilon}_K \in \Pi_K} \left(-h \evec\trsp + (b+h -2c\lambda^{K-1})\bm{\epsilon}_K\trsp + 2c\lambda \bm{\epsilon}_1\trsp 
    + 2c\sum_{r=2}^{K-1} (\lambda^r -\lambda^{r-1})\bm{\epsilon}_r\trsp \right)(\bd - \bq)\\
    &= \max_{\bm{\epsilon}_1 \in \Pi_1, \dots, \bm{\epsilon}_K \in \Pi_K} f(\bm{\epsilon}_1, \dots, \bm{\epsilon}_K)\trsp (\bd - \bq)
\end{align*}
% }
Such that $\Pi_{r} = \{\bm{\epsilon} \in \{0,1\}^n \;|\; \epsilon_i = \epsilon_j \; \forall i,j \in v, \forall v \in L_r\}$ for every $r =1, \dots, K$, and where we abuse notation and denote by $v \in T$ both the vertex $v$ of the tree and the set leaves in the subtree rooted at $v$.

Now, recall the expression of the offline problem~\eqref{eq:offline},
\begin{align*}
    \inf_{\bq \geq \bzero}\quad \left\{\sup_{\cD \in \cF(\bmu, \Sigma)} \mathbb{E}_{\bd \sim \cD}\; C(\bq, \bd)\right\},
\end{align*}
which can also be written as,
\begin{align*}
    \inf_{\bq \geq \bzero}\quad 
    \sup_{f: \mathbb{R}_+^n \rightarrow \mathbb{R}} \quad 
    & \int_{\mathbb{R}^n_+} C(\bq, \bd) f(\bd)\; \partial\bd  \\
    s.t \quad 
    & \int_{\mathbb{R}^n_+} \bd f(\bd)\; \partial\bd = \bmu, \\
    & \int_{\mathbb{R}^n_+} \bd\bd\trsp f(\bd)\; \partial\bd = \bmu\bmu\trsp + \Sigma, \\
    & \int_{\mathbb{R}^n_+} f(\bd)\; \partial\bd = 1, \\
    & f(\bd) \geq 0, \quad \forall \bd \in \mathbb{R}^n_+.
\end{align*}
Taking the dual of the inner moment maximization problem, we get,
\begin{align*}
   \min_{\bq \geq \bzero, Y, \by, y_0} \quad 
    &(\bmu\bmu\trsp + \Sigma)\cdot Y + \bmu\trsp \by + y_0 \\
    s.t. \quad 
    & \bd\trsp Y \bd + \bd\trsp \by + y_0 \geq C(\bq, \bd), \quad \forall \bd \geq \bzero
\end{align*}
where $A \cdot B$ denotes the Frobenius inner product of $A$ and $B$. Strong duality in the above problem holds as $\Sigma \succ 0$~\cite{smith1995generalized}, and hence the above formulation is equivalent to~\eqref{eq:offline}. Relaxing non-negativity of the demand vectors and replacing $C(\bq, \bd)$ by its expression,
{\fontsize{11.1px}{10px}
\begin{align*}
    \min_{\bq \geq \bzero, Y, \by, y_0} \quad 
    &(\bmu\bmu\trsp + \Sigma)\cdot Y + \bmu\trsp \by + y_0 \\
    s.t. \quad 
    & \bd\trsp Y \bd + \bd\trsp \by + y_0 \geq f(\bm{\epsilon}_1, \dots, \bm{\epsilon}_K)\trsp(\bd - \bq) \\
    &\forall \bd \in \mathbb{R}^n, \forall \bm{\epsilon}_1 \in \Pi_1, \dots, \bm{\epsilon}_K \in \Pi_K
\end{align*}
}
which can be formulated as the following exponentially sized SDP,
{\fontsize{10.5px}{10px}
\begin{align}
    \min_{\bq \geq \bzero, Y, \by, y_0} \quad 
    &\sum_{i,j \in X} (\mu_i \mu_j + \Sigma_{ij}) Y_{ij} + \sum_{i \in X} \mu_i y_i + y_0 \notag\\
    s.t. \quad 
    & 
    \begin{bmatrix}
        Y & \frac{\by - f( \bm{\epsilon}_1, \dots, \bm{\epsilon}_K)}{2}\\
        \frac{\by\trsp - f( \bm{\epsilon}_1, \dots, \bm{\epsilon}_K)\trsp}{2} & y_0 + f( \bm{\epsilon}_1, \dots, \bm{\epsilon}_K)\trsp \bq
    \end{bmatrix}\succeq 0, \notag
    \\
    &\forall \bm{\epsilon}_1 \in \Pi_1, \dots, \forall \bm{\epsilon}_K \in \Pi_K \notag.
\end{align}
}

\end{document}